\documentclass[12pt]{article}
\usepackage{amssymb,amsthm,amscd,amsmath,graphicx,color,xy}

\newcommand{\R}{{\mathbb R}}
\newcommand{\C}{{\mathbb C}}
\newcommand{\Z}{{\mathbb Z}}
\newcommand{\N}{{\mathbb N}}
\newtheorem{theorem}{Theorem}[section]
\newtheorem{corollary}[theorem]{Corollary}
\newtheorem{lemma}[theorem]{Lemma}
\newtheorem{proposition}[theorem]{Proposition}
\newtheorem{proposition-definition}[theorem]{Proposition - Definition}
\newtheorem{definition}[theorem]{Definition}
\newtheorem{remark}[theorem]{Remark}
\newtheorem{example}[theorem]{Example}

\newtheorem{assumption}[theorem]{Assumption}

\begin{document}

\thispagestyle{empty}
\begin{flushright}pre-print\end{flushright}
\par\bigskip\par
\vfill
\begin{center}
{\bfseries\Large
Quantum corrections to the holomorphic structure of the mirror bundle
along the caustic and the bifurcation locus}
\par\addvspace{20pt}
{\sc G. Marelli}
\par\medskip
Scuola Internazionale Superiore di Studi Avanzati (SISSA),\\
Via Beirut 4, 34013 Trieste, Italy
\end{center}
\vfill
\begin{quote} \footnotesize {\sc Abstract.}
Given, in the Lagrangian torus fibration $R^4\rightarrow R^2$, a Lagrangian
submanifold $L$, endowed with a trivial flat connection, the corresponding mirror object
is constructed on
the dual fibration by means of a family of Morse homologies
associated to the generating function of $L$, and it is provided with 
a holomorphic structure. 
Morse homology, however, is not defined
along the caustic $C$ of $L$ or along the bifurcation locus $B$, where the family
does not satisfy the
Morse-Smale condition. The holomorphic structure is extended
to the subset $C\cup B$, except cusps, 
yielding the so called quantum corrections to the mirror object.
\end{quote}
\vfill
\leftline{\hbox to8cm{\hrulefill}}\par
{\footnotesize
\par
\noindent\emph{E-Mail address:} {\tt marelli@kusm.kyoto-u.ac.jp}
}


\section{Introduction}
One of the reasons for which it may be desirable to embark on the study of mirror symmetry
for dual torus fibrations is that Calabi-Yau threefolds represents, at least
in String Theory, a case of remarkable interest. 
For the kind of problems this paper is concerned with, 
first steps in this direction were undertaken in
\cite{F1}, \cite{BMP1} and \cite{BMP2}: it was provided,
for the torus fibration $T^{2n}\rightarrow T^n$,
or, in general, for a smooth trivial family of Lagrangian tori,
a correspondence between Lagrangian submanifolds endowed with a flat connection
on one side and holomorphic bundles on the other. Some restrictions are necessary:
the most substantial, besides the absence of singular fibres in the fibrations, 
is that Lagrangian
submanifolds are assumed to exhibit no caustic. The constructions of the
correspondence 
are different: by
means of families of Floer homologies or by a kind of Fourier-Mukai
transform; however, at least in this simple setting, they are equivalent.
An attempt to allow for more general Lagrangian submanifolds, that is to include
the caustic, is contained in \cite{F2}: some ``quantum corrections'' must be added
to the construction and it is argued that these should be provided by
pseudoholomorphic disks; it is also conjectured that Floer homology and
pseudoholomorphic disks can be replaced, in an appropriate sense, by Morse
homology and gradient lines. The analysis of caustic and bifurcation locus, particularly
for the case of the perturbed elliptic umbilic, was carried on, in dimension 2,
in \cite{M1} and \cite{M2}, and an attempt of study of quantum corrections was
developed in \cite{M3}. 

The present paper tries to generalize the content of the three
previous works, by analysing the behaviour of gradient lines 
for given caustic and bifurcation locus and proposing quantum corrections
to get a whole defined holomorphic mirror object:
this is carried out by assigning, under suitable hypothesis, submanifolds $C$ and $B$, 
acting respectively as caustic and bifurcation locus,
and a class of orbit equivalent gradient vector fields for each subset $U_i$ determined by $C$ and $B$, 
followed by a study of bifurcations relating the phase portraits in nearby 
$U_i$ and $U_j$.
However, while 
in \cite{M3}, for the specific case
of the perturbed elliptic umbilic, the monodromy of the holomorphic structure
of the mirror bundle around cusps was considered, this paper does not deal the generalization
of this aspect. 
The fibration $T^4\rightarrow T^2$ is considered , though, since quantum corrections
are defined locally, the fibration $\R^4\rightarrow \R^2$ is kept in mind.
What follows is a summary of
the content of this paper.

In chapter 2 caustic and bifurcation locus associated to a Lagrangian
submanifold are introduced and their features are exposed: in particular,
it is studied when a codimension 1 subset of $\R^2$ represents caustic and
bifurcation locus of some Lagrangian submanifold.

Chapter 3 outlines the construction of the mirror object by Morse homology: it is not defined
along caustic and bifurcation locus, which form a codimension 1 subset of $\R^2$.

In chapter 4 it is studied the behaviour of gradient lines and of phase portraits
in a neighbourhood of folds of the caustic, 
not containing the bifurcation locus, 
and, as a consequence, quantum
corrections are defined, in order to extend the holomorphic structure of the
mirror object through such points of the caustic.

In chapter 5 gradient lines and phase portraits are analysed near codimension 1
points of the bifurcation locus, leading to the definition of quantum
corrections for glueing the holomorphic structure of the mirror object along 
these points.

Chapter 6 is concerned with the relative position of caustic and bifurcation locus
and with intersection of bifurcation lines.

Chapter 7 is devoted to show that the holomorphic mirror object has no
monodromy around the codimension 2 subset of points for which quantum
corrections, introduced in chapter 5 and 6, are not defined, and so it can be
extended across such points. However, as already said, cusps are not considered.
Theorem \ref{theo} sums up the achieved results. 

\newpage

\section{$CB$-diagrams}
Consider $\R^4$, endowed with its
canonical symplectic structure and standard Euclidean metric, and
the natural projection $\R^4\rightarrow\R^2$.
To any smooth function $f:\R^2\rightarrow\R$ (see \cite{AGZV} for details) 
it is associated
a 2-dimensional Lagrangian submanifold $L$, of which $f$ is its
generating function; the set of critical values of the projection
$L\hookrightarrow\R^4\rightarrow\R^2$ is called caustic of $L$ (or of $f$)
and denoted by $C$. Any Lagrangian submanifolds can be described locally by
some generating functions. Generically, in neighbourhood of any point of $L$, but
for a discrete set, $f$ can be assumed to be a function on the coordinates along the fibres.

\begin{proposition}
The caustic $C$ is a codimension 1 immersed
submanifold of $\R^2$ with singularities.
\end{proposition}

$C$ is a stratified submanifold: generically, folds form the stratum of codimension 1 and
cusps, the singularities of $C$, the stratum of codimension 2. Different branches
of $C$ can intersect transversely one with another, generically at folds.

In \cite{M1}, the family of functions $f_x:\R^2\rightarrow\R$, with $x\in\R^2$,
$f_x(y)=f(y)-xy$, is associated to $L$.  
The solutions of the gradient system
$$\frac{dy}{dt}=\nabla f_x(y)$$
are named gradient lines.
Observe that $f_x$ is a Morse function for $x\notin C$.
The subset of $\R^2\setminus C$, 
at whose points the gradient vector field $\nabla f_x$
is not Morse-Smale, is called bifurcation locus and denoted by $B$: 
at these points a saddle-to-saddle
separatrix occurs in the phase portrait of $\nabla f_x$.

\begin{proposition}
The bifurcation locus $B$ is a codimension 1 immersed
submanifold of $\R^2$.
\end{proposition}

$B$ is a stratified submanifold: generically, at codimension 1 points one saddle-to-saddle separatrix
occurs, at codimension 2 points two saddle-to-saddle separatrixes occur. Codimension 2
points are the intersections of the codimension 1 stratum.

Thus $C$ and $B$ determine in $\R^2$ a diagram.

\begin{definition}
\label{cbl}
\rm
The $CB$-diagram associated to a Lagrangian submanifold $L$ of $\R^4$,
is the partition of $\R^2$ 
determined by the caustic
$C$ and by the bifurcation locus $B$ of $L$. 
\end{definition}

Definition \ref{cbl} is extended as follows:

\begin{definition}
\rm
Given two codimension 1 submanifolds $C$ and $B$ of $\R^2$, holding the features outlined
above of, respectively,
the caustic and the bifurcation locus, 
a $CB$-diagram generated by $C$ and $B$ is
the partition $U_i$, $i\in I$, $U_i$ connected, determined in $\R^2$ by $C$ and $B$, 
together with an assignment, for each $i\in I$, of a phase portrait $P_i$.
It is denoted by $(C,B,(U_i,P_i)_{i\in I})$.
\end{definition}

%
%
The question now is: when is a $CB$-diagram, generated by $C$ and $B$ and
with family of phase portraits $P_i$,
the $CB$-diagram of some Lagrangian submanifold $L$, whose caustic
and bifurcation locus are respectively $C$ and $B$, and, for $x\in U_i$, $\nabla f_x$
is orbit isotopic (see the discussion preceding definition \ref{adm})
to $P_i$, where $f$ is a generating function of $L$? 
The answer 
depends first of all on the family of phase portraits $P_i$: 
in fact, $P_i$ must be the phase portrait of some gradient vector
field. Into this direction, a result is provided by the following lemma:

\begin{lemma}
A gradient vector field in $\R^2$ has no periodic 
orbits, moreover, if structurally stable, has only hyperbolic 
critical points and the intersection of the stable and unstable 
submanifolds $W^s(p)$ and $W^u(q)$ of any two critical points $p$ and $q$  
is always transverse.
\end{lemma}

So to any $x\notin C\cup B$ it must be associated a vector field
exhibiting only hyperbolic critical points, which in $\R^2$ turn out to be
either stable nodes or unstable nodes or saddles, and with no 
saddle-to-saddle separatrix. In this case we can prove the following
proposition:

\begin{proposition}
\label{grad}
Given a phase portrait exhibiting a finite number of hyperbolic critical points,
no periodic orbits and no saddle-to-saddle separatrixes, then there exists
a function $f:\R^2\rightarrow\R$ such that $\nabla f$ has the
given phase portrait.
\end{proposition}
\begin{proof}
Let $p_j$, for $j\in J$, where is $J$ finite,
be the critical points in the given phase portrait, and consider, 
for each $j$, a ball $B_j$ centred in $p_j$ and such that $B_j\cap
B_k=\varnothing$ for $j\neq k$. As explained, each $p_j$ is expected to be either a node,
stable or unstable, or a saddle: since there exist gradient vector field
with such singular points, choose bounded functions $f_j$ such that 
$\nabla f_j$ exhibits in $B_j$
the critical point $p_j$ of the type as prescribed by the given phase portrait.
Assume $B_j$ itself as domain of $f_j$.
Choose now a point $p_m$ among the critical 
points $p_j$ and for each $j\neq m$ choose a path $\gamma_j$ from
$p_m$ to $p_j$. Choose also a tubular neighbourhood $N_j$
of $\gamma_j$ such that $(N_i\cap N_j)\setminus B_m=\varnothing$. 
To extend $f_m$ along $N_j$, assume $\partial B_j\cap N_j$ to 
be a level set of $f_j$ as well as a fibre of $N_j$, 
and take the fibres of $N_j$ as the level
curves of the extended functions: since each $f_j$ is defined
up to a constant, $f_m$ can be matched with $f_j$. Proceed
till obtaining a bounded function $f$ defined on the
contractible subset $B=(\cup_j B_j)\cup(\cup_j N_j)$. Now
$f$ can be extended to a function defined on the whole $\R^2$: indeed, 
as $B$ is contractible,
the problem is equivalent to extending a bounded function from a ball
$B$ to the whole $R^2$. Note that this can be carried out without introducing new critical 
points.
\end{proof}

The proof relies on the fact that there exist gradient vector fields exhibiting a stable node or an unstable
node or a saddle. Since there are also gradient vector fields, though non-structurally stable, exhibiting a saddle-node (over folds of the caustic,
when referring to the setting considered in this paper), or a point, which can be called ``saddle-node-saddle'',
given by two saddles and a node glued together  
(over cusps), or saddle-to-saddle separatrixes (over
points of the bifurcation locus), the following corollary generalizes proposition \ref{grad}:

\begin{corollary}
\label{grad2}
Given a phase portrait exhibiting a finite number of critical points (hyperbolic or saddle-nodes
or saddle-node-saddles) and
no periodic orbits, then there exists
a function $f:\R^2\rightarrow\R$ such that $\nabla f$ has the
given phase portrait.
\end{corollary}

Let $U_i$, with $i\in I$ and where $I$ is either finite or $\N$, be
the connected components of the $CB$-diagram associated to $L$. For every $i\in I$,
$x_1,x_2\in U_i$ and path $\gamma(t)\subset U_i$, with $t\in[0,1]$,
$\gamma(0)=x_1$ and $\gamma(1)=x_2$, there exists an orbit isotopy 
from $x_1$ to $x_2$, that is 
a smooth family $\Phi_t$ of diffeomorphisms of $\R^2$, with $t\in[0,1]$, 
such that $\Phi_0=Id$ and
$\Phi_t$ provides an orbit equivalence between $\nabla f_{x_1}$ and $\nabla f_{\gamma(t)}$
for every $t\in[0,1]$. Therefore,
another property that the family $(P_i)$ of a $CB$-diagram generated by $C$ and $B$ 
must satisfy, in order to be the $CB$-diagram of some Lagrangian submanifold, 
is that, intuitively,
if $U_k$ and $U_l$ are separated by $C$ or $B$, then it must be possible
to switch from $P_k$ to $P_l$ by adding or removing a pair of critical points
(forming at $C$ a degenerate critical point) or by exchanging the separatrixes
of two saddles (forming at $B$ a saddle-to-saddle separatrix). The last 
considerations can be resumed
rigorously in the following definition:

\begin{definition}
\label{adm}
\rm
A $CB$-diagram $(C,B,(U_i,P_i)_{i\in I})$ is admissibile if and only if 
$I$: it is finite; $U_i$ is open for every $i\in I$;
each $P_j$
exhibits a finite number of only hyperbolic critical points
and no closed orbits;
if $\partial U_i\cap\partial U_j\neq\varnothing$ then 
for every path $\gamma:[0,1]\rightarrow \bar{U_i}\cup \bar{U_j}$
such that $\gamma([-1,0))\subset U_i$, $\gamma((0,1])\subset U_j$ and $\gamma(0)\in \partial U_i\cap\partial U_j$,
there is a smooth family $X_t$ of vector fields, with $t\in[-1,1]$, 
such that $X_t$ has phase portrait orbit isotopic to $P_i$ for
$t\in[-1,0)$, and to $P_j$ for $t\in(0,1]$, and such that
$X_0$ exhibits either
a degenerate critical point or a saddle-to-saddle separatrix depending on whether
$\gamma(0)$ belongs respectively to $C$ or $B$; moreover, any two family $X_t$ and $X'_t$
as above are orbit isotopic. 
\end{definition}

\begin{theorem}
\label{mgrad}
Any admissible $CB$-diagram $(C,B,(U_i,P_i)_{i\in I})$
is the $CB$- diagram
of some Lagrangian submanifolfd $L$ at least on any compact subset of $\R^2$, 
in the sense that $L$ has caustic and 
bifurcation locus diffeomorphic, respectively, to $C$ and $B$, determining
a partition $(W_i)_{i\in I}$ of $\R^2$, with $W_i$ diffeomorphic to $U_i$,
and, for each $x\in W_i$ such that the projection of $L$ over $W_i$ is non-empty, 
the vector fields $\nabla f_x$ has phase portrait orbit
isotopic to $P_i$,
where $f$ is a local generating function of $L$.
\end{theorem}
\begin{proof}
Observe that, since the critical points of $\nabla f_x$ correspond to the intersection points of $L$ with
the fibre over $x$, $L$ will be defined only over those $U_i$ endowed with a phase portrait $P_i$ having at least a critical point.
So, for every $i\in I$ such that $P_i$ is as described above, choose a point $p_i\in U_i$. 
By proposition \ref{grad}, choose a gradient vector field $X_i$ with
phase portrait orbit equivalent to $P_i$ and such that $X_i(0)=0$. Let $\tilde{f}_i$ be a function in the variable $y=(y_1,y_2)$ such that 
$\nabla \tilde{f}_i=X_i$ and let $f_i=\tilde{f}_i+p_i\cdot y$: observe that $\nabla (f_i)_{p_i}=X_i$, and, since $p_i\notin C\cup B$,
there is a subset $V_i'$ of $U_i$ such that for all $x\in V_i'$ the vector field $\nabla (f_i)_x$ is orbit equivalent to $X_i$.
Observe also that $(\partial f_i/\partial y)(0)=(\partial\tilde{f}_i/\partial y)(0)+p_i=X_i(0)+p_i=p_i$, so 
the equation $x=\partial f/\partial y$ defines a Lagrangian submanifold $L_i$ over a neighbourhood $V_i^1$ of $p_i$ which
can be assumed contained in $V_i'$. 

$L_i$ can be extended over an open subset $V^2_i$ of $U_i$ diffeomorphic to $U_i$: clearly $V^2_i=V_i^1$  when
$\pi_1(U_i)=0$; otherwise, consider for simplicity the case where $\pi_1(U_i)=\Z$: since $p_i\notin C$ and
$V^1_i$ is an open ball, $L_i$ can be generated by a finite set of functions $g_i^j$, one for each sheet $L_i^j$ of $L_i$ over $V^1_i$ and defined
on $V^1_i$ (this corresponds to the fact that each sheet $L_i^j$ of $L_i$ can be seen as the graph of a closed 1-form $\sigma_i^j$, which, being
$V^1_i$ an open ball, is exact, that is, $\sigma_i^j=dg_i^j$ for some function $g_i^j$); it is enough now to extend each $g_i^j$ to a function defined on an open subset
$V^2_i$ diffeomorphic to $U_i$. This argument also shows how to extend a Lagrangian submanifolds, defined over two disjoint open
subsets of the base of the fibration, where it has the same number of sheets and no critical points, 
onto a new subset containing the two subsets.

The admisibility of the given $CB$-diagram implies that to each point of $C$ and $B$ it is associated a vector field,
which by corollary \ref{grad2} can be assumed to be a gradient vector field. Since $C$ and $B$ generically have two strata, choose,

Choose a point $q_{ik}$ on
each connected component $C_{ik}$, $B_{ik}$ of the codimension 1 stratum of $C$, respectively $B$,
where ``intersection points'' between $C$ and $C$ are removed (the issue of intersection points will be analyzed in chapter \ref{int}). 
Each $C_{ik}$, $B_{ik}$ will bound two subsets $U_i$ and $U_k$ in the partition of the given $CB$-diagram. 
Let $\gamma_{ik}$ be a path
from $p_i$ to $p_k$ as in definition \ref{adm}, with $\gamma_{ik}(0)=q_{ik}$ and associated
family of vector fields $X_{ik}^t$, and  let $f_{ik}$ be a function such that
$(\partial f_{ik}/\partial y)(0)=q_{ik}$ and
$\nabla (f_{ik})_{\gamma(t)}$ is orbit equivalent to $X_{ik}^t$, for $t$ in a neighbourhood $N_{ik}$ of 0: this is possible
because $X_{ik}^0$ exhibits a saddle-node or a saddle-to-saddle separatrix, particularly it is not stable, and after a small
perturbation, that is for $t$ in a neighbourhood $N_{ik}$ of 0, it is orbit equivalent, by the admissibility of the given $CB$-diagram,
to the vector fields $\nabla (f_i)_{p_i}=X_i$ and $\nabla (f_k)_{p_k}=X_k$ in, respectively, $\gamma_{ik}(N_{ik})\cap U_i$
and $\gamma_{ik}(N_{ik})\cap U_k$. The function $f_{ik}$ defines a Lagrangian submanifold $L_{ik}$ in a neighbourhood $V_{ik}$
of $q_{ik}$; $V_{ik}$ can be suppossed to not intersect any of the subsets $V^2_j$ constructed above. For a generic choice of $f_{ik}$, 
$V_{ik}$ is a ball such that, along one of its diameters, the vector field $\nabla (f_{ik})_{x}$ is orbit equivalent to $X_{ik}^0$,
while, in the two half-disks determined by such a diameter, it is orbit equivalent to respectively $X_i$ and $X_k$. This construction
can be performed also for every point of the codimension 2 stratum of $C$ and $B$ and for the ``intersection points'' between $C$ and $B$.

The Lagrangian submanifolds defined above over the open sets $V^2_j$, $V_{ik}$ and in open neighbourhoods of cusps, of 
intersection points of bifurcation lines and of intersection points between $C$ and $B$, 
can be glued together, by extending them along every path $\gamma_{ik}$, in a new Lagrangian submanifold $L$: indeed, away
from points of $C$ and $B$, which form a codimension 1 subset, generating functions,
one for each sheets of the Lagrangian submanifolds, can be considered, and these functions, as
already explained above, can be extended along every paths $\gamma_{ik}$.

The Lagrangian submanifold $L$ so obtained can be extended now along the remaining points of the caustic $C$.
Suppose indeed to have two Lagrangian submanifolds $L_1$ and $L_2$ defined over two disjoint open balls $W_1$ and $W_2$
in the $(x_1,x_2)$-plane, such that they exhibits a caustic, formed only by folds, along a diameter $C_i$ of $W_i$, for
$i=1,2$. For simplicity, suppose that $L_i$ has two sheets over one of the two connected components determined by $C_i$
and no sheet over the other. Let $b_i$ one of the two points in $\partial W_i\cap \bar{C}_i$ and consider a path $C:[1,2]\rightarrow\R^2$
such that $C((1,2))\cap W_i=\varnothing$, $C(i)=b_i$, for $i=1,2$ and such that it extends the paths $C_i$ (it can be $W_i=\varnothing$
for $i=1$ or $i=2$ or $W_1=W_2$ but $b_1\neq b_2$). By admissibility of the given $CB$-diagram, it can be assumed that $L_1$ and $L_2$ have the same number of sheets
in the component of $W_i$ lying on the same side with respect to the path $C_1\cup C\cup C_2$. Choose coordinates $t$ along $C_1\cup C\cup C_2$
and $u$ such that, if $(t,u,y_t,y_u)$ are canonical coordinates, 
$C_1\cup C\cup C_2$ lies on the $t$-axis and such that, since $C_i$ contains only folds, $L_i$ has equation
\begin{displaymath}
\left\{ \begin{array}{ccc}
u & = & y_u^2 \\
t & = & y_t
\end{array} \right.
\end{displaymath}
in $W_i$. This equation gives also the wanted extension along $C$, when the coordinate $t$ corresponds to points of $C$.
Since the map $(x_1,x_2,y_1,y_2)\rightarrow(t,u,y_t,y_u)$ is a Lagrangian equivalence of the Lagrangian bundle $\R^4\rightarrow\R^2$,
the extension of $L$ along $C$ is obtained.
 
Finally, the extension of $L$ to the whole $\R^2$ is carried out as already done above, since the subset to which now $L$ is extended does not
contain any point of $C$.
Because the vector fields in $U_i$ and $U_k$, when these have $B$ has common boundary, are not orbit equivalent and since
the given $CB$-diagram is admissible, it follows that for each path from $U_i$ to $U_k$ there is a point along this path
where the corrsponding vector field is equivalent to the one chosen over $q_{ik}$.
In principle, this point is not unique,
however
if further bifurcation points
appear, they must appear 
in pairs, that is, each pair will mark the apperance of the samme saddle-to-saddle separatrix, so that the two 
bifurcations cancel each other and 
the admissibility of the given $CB$-diagram is preserved; after a perturbation, each pair of points can be removed, at least
on a compact subset. That the bifurcation locus of $L$ is diffeomorphic to the given $B$ follows from the fact that
$V^2_i$ has been constructed diffeomorphic to $U_i$.
\end{proof}

\begin{remark}
\rm
%
Note that if two Lagrangian submanifolds have diffeomorphic $CB$-diagrams,
in the sense of theorem \ref{mgrad}, this does not imply that they are Lagrangian
equivalent: in fact, for example, two Lagrangian equivalent submanifolds 
have diffeomorphic caustic, however the converse is not true.
\end{remark}

\section{The mirror bundle}
This chapter wants to be only a summary of the idea of the construction of the mirror object
using families of Floer homologies, or, as in this paper, families of Morse
homologies. Details are, in fact, already exposed in \cite{F1}, \cite{F2}
and \cite{M3}. 

In the Lagrangian torus fibration $T^4\rightarrow T^2$, consider
a 2-dimensional Lagrangian submanifold $L\hookrightarrow T^4$,
endowed with a flat connection $\nabla$, and let $f:\R^2\rightarrow\R$ be
a generating function of $L$. Floer homology for families of Lagrangian submanifolds
is treated in \cite{F3} and its application to mirror symmetry in the construction
of the mirror object on the dual fibration is in \cite{F1} and \cite{F2}:
the fibre of the mirror object, an element of $DCoh(\hat{X})$, 
over a point $(x,w)$, where $x\in T^2$ and $w\in\hat{F}_x$,
is given by the
intersection Floer homology $HF((L,\nabla),(F_x,w))$, where $F_x$ is the fibre over $x$
and $\hat{F}_x$ the dual fibre
(all the problems concerning the definition or the existence of $HF(L,F_x)$
are not discussed here, see rather the monograph \cite{FOOO}). A holomorphic
frame is then defined (see \cite{F1} and \cite{F2}), glueing the fibres in a complex of
holomorphic bundles.
In particular,
in \cite{F2} and \cite{FO} it is conjectured that near the caustic the moduli space of
pseudoholomorphic disks is isotopic, after perturbation, 
to the moduli space of gradient lines of the generating function $f$. 
This conjecture is used in this paper: the fibre of the mirror object over $x$
is defined as the Morse homology $HM(f_x)$, when $x\notin C\cup B$: in fact
in this case, $f_x$ is a Morse function and the Morse-Smale condition is satisfied
(of course, all conditions on $f$ necessary to ensure the existence of Morse 
homology are assumed: to this purpose,
for everything concerning Morse homology, see \cite{J}
and, above all, the monograph \cite{S}). A holomorphic frame is then defined, yielding
a complex of holomorphic bundles, away from $x\in C\cup B$:
writing $\nabla=d+A$,
a section $e(x)$ of the mirror object turnes out to be holomorphic and descends on the torus fibres when multiplied by
the weight 
$$exp\Big[2\pi\Big(\frac{h(x)}{2}-\frac{A(x)}{4\pi}+i\frac{\partial h}{\partial x}\cdot w\Big)\Big]$$
where $h$ is a multi-valued function on the base such that each sheet of $L$ is locally the
graph of $dh$ (in other words, $h$ is a set of local generating functions, defined in the
coordinates of the base, one for each sheet of $L$).
The way to extend the holomorphic
structure through the subset $C\cup B$ is provided by ``quantum corrections'', that is
morphisms glueing the mirror object along this subset. 
This is the purpose of the present paper.
Since quantum corrections are defined locally, it is enough to consider
the Lagrangian fibration
$\R^4\rightarrow\R^2$.

\begin{remark}
\rm
Here is a kind of road map showing the way leading to the extension of the
holomorphic structure of the mirror obejct across $C\cup B$:
\begin{itemize}
\item
in an admissible $CB$-diagram $(C,B,(U_i,P_i)_{i\in I})$, confront the phase portraits $P_i$ and $P_j$
for nearby $U_i$ and $U_j$;
\item
show that the Morse homologies $HM(f_x)$, for $x\ U_i$ are isomorphic to those for $x\in U_j$;
\item
pick up an isomoprhism: the choice depends on which kind of points form the common boundary of $U_i$ and $U_j$,
that is, folds not limit points of the bifurcation locus or codimension 1 bifurcation points; it is defined a map
at the chains level, that is, on generators of the Morse complex, that is, on critical points of $\nabla f_x$, inducing
an isomorphism (the quantum correction) in homology;
\item
with such a glueing, check that there is no monodromy when going around the set of the remaining points, that is,
folds which are limit of the bifurcation locus and codimension 2 bifurcation points (cusps are not considered, as 
already said, in this paper), which form a codimension 2 subset; this means that the bundle whose fibres are
$HM(f_x)$ is endowed now with non-vanishing sections which can be extended to any point of $C\cup B$ (but
for cusps);
\item
observe that, since a holomorphic section is obtained from a section of $HM(f_x)$ multiplied by a weight, it follows
that the glueing isomorphisms introduced in the previous steps induce a glueing at the level of holomorphic sections, allowing to
extend them across $C\cup B$ (but
for cusps);
\item
$HM(f_x)$ is now endowed with a holomorphic structure which can be extended to $C\cup B$ (but
for cusps).
\end{itemize}
\end{remark}
 
\section{The caustic}
Points of the caustic are characterized by a degenerate critical point
of $\nabla f_x$, which, if $x$ is a fold, is related to a so-called birth-death pair.

\begin{definition}
\rm
A smooth 2-parameters family of vector fields, defined on an open subset $U$
of the plane, exhibits a birth-death pair if there exists a curve $C$,
decomposing $U$ into two connected components $U_1$ and $U_2$,
such that the family has $k+2$
critical points $p_1(s)$, ... ,$p_k(s)$, $p_{k+1}(s)$, $p_{k+2}(s)$ in $U_1$
and $k$ critical points $p_1(s)$, ... ,$p_k(s)$ in $U_2$, where $s$ is the parameter,
and the two critical points $p_{k+1}$, $p_{k+2}$ 
converge, for $s$ converging to $c\in C$, to a degenerate critical
point $p(c)$, called birth-death point. The pair of critical points 
$(p_{k+1}(s),p_{k+2}(s))$ is called a birth-death pair.
\end{definition}

\begin{proposition}
For a generic $f$, at any fold $x\in C$, 
the function $f_x$ exhibits a birth-death point.
\end{proposition}
\begin{proof}
For simplicity, the fold $x$ can be assumed to be the origin $(0,0)$.
After a Lagrangian equivalence, since $f$ is generic, the local generating 
function $f$ of $L$ can be written, for example, as $f(y_1,y_2)=y_1^3$.
This shows that $C$ determines in any neighbourhood of $(0,0)$ 
two subsets, characterized by the fact that the intersection of the
fibre $F_{x'}$, for $x'$ in these two subsets, and $L$ is either empty or  contains
two points: these form a birth-death pair and, for 
$x'\rightarrow (0,0)$, glue together into a birth-death point.
Equivalently, this means that the vector field $\nabla f_{x'}$
exhibits two critical points in one of these components, glueing
together in a degenerate critical point at $(0,0)$, and no critical
points in the other component.
\end{proof}

In the example considered in \cite{M1} and \cite{M2}, concerned with the cusp
and the elliptic umbilic, the birth-death pairs were formed
by a saddle and an unstable node. Instead, for the hyperbolic
umbilic, the birth-death pairs are given by a saddle and a stable node.
These are also the only two cases that can be met.

Denote by $\mu(p_i)$ the Morse index of a critical point $p_i$.

\begin{proposition}
If $(p_1,p_2)$ is a birth-death pair corresponding to some fold, then
$|\mu(p_1)-\mu(p_2)|=1$.
\end{proposition}
\begin{proof}
At a fold $x$, by definition, $rk(H(f_x))(c(x))=n-1$, where $c(x)$ is the 
birth-death point: this means that both $p_1$ and $p_2$ has an eigenvalue $e(p_i)$ 
which vanishes when $p_1$ and $p_2$ glue together in $c(x)$; if $p_1$ and $p_2$ glue over a fold then
$e(p_1)$ and $e(p_2)$ have opposite sign. The remaining eigenvalues of 
$p_1$ and $p_2$ can not vanish and so have the same signs.
\end{proof}

The quantity $\mu(p_1)-\mu(p_2)$ is named relative Morse index of $p_1$ and $p_2$.

\begin{lemma}
\label{exgrl}
Given a birth-death pair $(p_1,p_2)$ with $\mu(p_1)>\mu(p_2)$, 
then there exists a unique gradient line $\gamma_{p_1,p_2}$.
\end{lemma}
\begin{proof}
The local phase portrait is topologically equivalent to that
given by the generating function $f(y_1,y_2)=y_1^3\pm y_2^2$ (see \cite{GH}), for which
the existence of a gradient line $\gamma_{p_1,p_2}$ can be proved by direct computation.
If $W^u(p_1)\cap W^s(p_2)$ is not empty, then 
$\mu(p_1)-\mu(p_2)=dim(W^u(p_1)\cap W^s(p_2))=1$. Unicity follows now from the
fact that $W^u(p_1)\cap W^s(p_2)$ is connected: if not, the birth-death
point would exhibit a homoclinic orbit, which can not occur for a gradient vector field.
\end{proof}

So, as explained, the caustic divides a small open neighbourhood $U$
of $x$ into two open subsets $U_1$ and $U_2$: for example, for $x'\in U_1$,
$\nabla f_{x'}$ has critical points $p_1(x')$, ... , $p_k(x')$, $p_{k+1}(x')$, $p_{k+2}(x')$, 
where $(p_{k+1}(x'),p_{k+2}(x'))$
is the birth-death pair, and so $p_1(x')$, ... , $p_k(x')$
are the critical points of $\nabla f_{x'}$ for $x'\in U_2$.

\begin{definition}
\rm
A square in a phase portrait is a set 
$$(un,s_1,s_2,sn;\gamma_{un,s_1},\gamma_{s_1,sn},\gamma_{un,s_2},\gamma_{s_2,sn})$$
whose elements are an unstable node $un$, two saddles
$s_1$ and $s_2$, a stable node $sn$, and for each saddle
a pair of separatrixes connecting them to the nodes.
If $s_2=s_1$ and either $\gamma_{un,s_1}=\gamma_{un,s_2}$ or
$\gamma_{s_1,sn}=\gamma_{s_2,sn}$ the square is said to be degenerate.
\end{definition}

Theorem \ref{jost} is quoted from \cite{J} : it is essential in
defining the Morse complex, 
explaining
the structure of the boundary of the moduli space ${\cal M}(p,q)=(W^u(p)\cap W^s(q))/\R$ 
of gradient lines between two critical
points $p$ and $q$ with relative Morse index $\mu(p,q)=2$.
Given such points $p$ and $q$, if there exists a critical point $r$, with $\mu(p,r)=\mu(r,q)=1$,
and gradient lines $\gamma_{p,r}$ and $\gamma_{r,q}$, the triple 
$(\gamma_{p,r},r,\gamma_{r,q})$ is called a broken gradient line from $p$ to $q$, and
denoted also by $\gamma_{p,r}\sharp\gamma_{r,q}$. Observe that the gradient lines
of a square form two broken gradient lines $\gamma_{un,s_i}\sharp\gamma_{s_i,sn}$
for $i=1,2$.
Finally, $\gamma_1\sharp\gamma_2\neq\gamma^{'}_1\sharp\gamma^{'}_2$ if and only if
$\gamma_1\neq\gamma^{'}_1$ or $\gamma_2\neq\gamma^{'}_2$.

Some hypothesis on the function $f$ are needed in order to ensure a good
behaviour of ${\cal M}(p,q)$ and so to define a differential $\partial$ such that
$\partial^2=0$. 

\begin{definition}
\rm
A function $f$ of class $C^1$ satisfies the Palais-Smale condition if
every sequence $(x_n)$, such that $|f(x_n)|$ is bounded and
$|df(x_n)|\rightarrow 0$ for $n\rightarrow\infty$, admits a convergent 
subsequence.
\end{definition}

\begin{theorem}
\label{jost}
Suppose $f$ is a function of class $C^3$, having only non degenerate
critical points, satisfying the Palais-Smale condition and the
Morse-Smale condition. Let $p$ and $q$ be two critical points
of $f$, connected by the flow, and such that $\mu(p)-\mu(q)=2$.
Suppose that the space of gradient lines ${\cal M}(p,q)$
from $p$ to $q$ is contained in a flow-invariant compact set.
Then each connected component of  ${\cal M}(p,q)$
either is compact after including $p$ and $q$, and so 
diffeomorphic to the 2-sphere, or its boundary consists of two
different broken gradient lines from $p$ to $q$.

Conversely each broken gradient line from $p$ to $q$ is contained in
the boundary of precisely one component of ${\cal M}(p,q)$.
\end{theorem}

In other words, each connected component of ${\cal M}(p,q)$, 
if non-compact after including $p$ and $q$, determines by means of its boundary a square.
Theorem \ref{jost}, quoted from \cite{J}, is proved by assuming, for a matter of
convergence, a compactness hypothesis:
$W^u(p)\cap W^s(q)$
is contained in a flow-invariant compact subset. 
In the case we are considering, that is $\R^2$,
this could be a problem:
there could be an unstable node $un$ and stable node $sn$ such that
$W^u(un)\cap W^s(sn)$ is not bounded and so the above compactness hypothesis
can not be applied. Actually, what is important for the purposes of Morse theory
is to show, when $un$ is connected to some saddle, 
that every connected component of $W^u(un)\cap W^s(sn)$ is
bounded by two different broken gradient lines (or, in other words, it forms
a square): this allows to prove that 
$\partial^2=0$.

Suppose, for example, that in $\R^2$ $W^u(un)\cap W^s(sn)$
is not bounded, consider the compactification $S^2$ of $\R^2$ and suppose also that the phase
portrait can be extended to $S^2$ without
adding new critical points: as $\infty$ is not a critical point
and because of unicity of solutions, the
gradient line through $\infty$, in the compactification $S^2$, 
either belongs to $(W^u(un)\cap W^s(sn))/\R^2$, if it connects $un$ to $sn$, 
or bounds $(W^u(un)\cap W^s(sn))/\R^2$, 
if it connects $un$ to a saddle or a saddle to $sn$.
So, in the first case, except when compact after including $un$ and $sn$,
$W^u(un)\cap W^s(sn)$ has in $\R^2$ two non-bounded connected components,
which, however, in $S^2$ form a unique connected bounded component 
of $W^u(un)\cap W^s(sn)$ (see figure \ref{fig1}.1);
in the second case, 
the gradient line through $\infty$ is part of one of the two broken gradient lines
bounding $W^u(un)\cap W^s(sn)$ in $S^2$. So, even though the compactness hypothesis 
fails in $\R^2$, the equation $\partial^2$ still holds true.

\begin{center}
\setlength{\unitlength}{1cm}
\begin{picture}(4,4)
\label{fig1}
\thinlines

\put(3.5,2){\circle*{.1}}
\put(2,.5){\circle*{.1}}
\put(.5,2){\circle*{.1}}
\put(2,3.5){\circle*{.1}}
\put(1.5,2){\circle*{.1}}
\put(2.5,2){\circle*{.1}}

\scriptsize
\put(2,3.6){$\displaystyle s_1$}
\put(0.2,1.7){$\displaystyle un$}
\put(2,.2){$\displaystyle s_2$}
\put(3.5,2.1){$\displaystyle sn$}
\put(2.6,2){$\displaystyle un_1$}
\put(1,2){$\displaystyle sn_1$}
\put(0,2.1){$\displaystyle \infty$}
\put(3.9,1.7){$\displaystyle \infty$}
\normalsize

\qbezier(3.5,2)(2.75,1.25)(2,.5)
\qbezier(3.5,2)(2.75,2.75)(2,3.5)
\qbezier(.5,2)(1.25,1.25)(2,.5)
\qbezier(.5,2)(1.25,2.75)(2,3.5)
\qbezier(1.5,2)(1.75,2.75)(2,3.5)
\qbezier(2.5,2)(2.25,2.75)(2,3.5)
\qbezier(1.5,2)(1.75,1.25)(2,.5)
\qbezier(2.5,2)(2.25,1.25)(2,.5)
\qbezier(3.5,2)(3.75,2)(4,2)
\qbezier(.5,2)(.25,2)(0,2)

\qbezier(1.75,2.75)(1.8,2.775)(1.85,2.8)
\qbezier(1.75,2.75)(1.725,2.8)(1.7,2.85)

\qbezier(2.25,2.75)(2.3,2.725)(2.35,2.7)
\qbezier(2.25,2.75)(2.225,2.7)(2.2,2.65)

\qbezier(1.75,1.25)(1.8,1.225)(1.85,1.2)
\qbezier(1.75,1.25)(1.725,1.2)(1.7,1.15)

\qbezier(2.25,1.25)(2.3,1.275)(2.35,1.3)
\qbezier(2.25,1.25)(2.225,1.3)(2.2,1.35)

\qbezier(3.75,2)(3.8,2.05)(3.85,2.1)
\qbezier(3.75,2)(3.8,1.95)(3.85,1.9)

\qbezier(.25,2)(.3,2.05)(0.35,2.1)
\qbezier(.25,2)(.3,1.95)(0.35,1.9)

\qbezier(2.75,2.75)(2.7,2.75)(2.65,2.75)
\qbezier(2.75,2.75)(2.75,2.8)(2.75,2.85)

\qbezier(1.25,2.75)(1.2,2.75)(1.15,2.75)
\qbezier(1.25,2.75)(1.25,2.74)(1.25,2.65)

\qbezier(2.75,1.25)(2.7,1.25)(2.65,1.25)
\qbezier(2.75,1.25)(2.75,1.2)(2.75,1.15)

\qbezier(1.25,1.25)(1.2,1.25)(1.15,1.25)
\qbezier(1.25,1.25)(1.25,1.3)(1.25,1.35)


\end{picture}
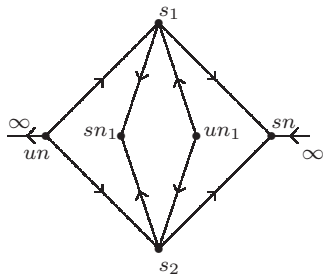
$Fig.~\ref{fig1}.1:~W^u(un)\cap W^s(sn)~and~the~gradient~line~through~\infty$
\end{center}

Suppose now that both $W^u(un)\cap W^s(sn)$
is not bounded and the phase portrait can not be extended to the compactification
$S^2$ of $\R^2$ unless
adding further critical points.
Let ${\cal M}$ be an unbounded connected component of $W^u(un)\cap W^s(sn)$ which
is turned into a bounded component ${\cal M'}$ by the addition of a new
critical point in the phase portrait: then $\partial{\cal M'}\neq\partial{\cal M}$,
on the other hand, since theorem  \ref{jost} implies that $\partial{\cal M'}$
consists of two distinct broken gradient lines, it follows that $\partial{\cal M}$ 
consists instead of a single broken gradient lines.
This means that the Morse differential does not satisfy $\partial^2\neq0$
for the given phase portrait in $\R^2$.
Thus the following assumption is made:

\begin{assumption}
\label{assumption}
Each component of 
$W^u(un)\cap W^s(sn)$ satisfies the compactness hypothesis of theorem \ref{jost},
that is, it is contained in a flow-invariant compact subset either
of $\R^2$ or, provided the phase portrait
can be extended to $S^2$ without the addition of further
critical points, of the compactification $S^2$ of $\R^2$.
\end{assumption}

\begin{lemma}
\label{php}
Under the hypothesis of theorem \ref{jost} and 
assumption \ref{assumption},
given a saddle $s_1$, a stable node $sn$, an unstable nodes $un$ and
gradient lines $\gamma_{un,s_1}$ and $\gamma_{s_1,sn}$
then there exists a saddle $s_2$ ($s_2=s_1$ is a possibility), with separatrixes
$\gamma_{un,s_2}$ and $\gamma_{s_2,sn}$, thus
forming a square (degenerate if $s_2=s_1$), 
in $R^2$ or eventually in its compactification $S^2$, 
together with $s_1$, $sn$ and $un$.
\end{lemma}
\begin{proof}
Theorem \ref{jost}, eventually applied on $S^2$, implies that
the boundary of $W^u(un)\cap W^s(sn)$ contains, besides
the broken gradient line $\gamma_{un,s_1}\sharp\gamma_{s_1,sn}$,
a second broken gradient line $(\gamma_{un,s_2},s_2,\gamma_{s_2,sn})$
for some saddle $s_2$.
\end{proof}


Here are few examples.

\begin{example}
\label{ex1}
\rm
Given a square $(un,s_1,s_2,sn)$, the two separatrixes of each saddle $s_i$ not forming
the sides of the squares lie on the same side with respect to the
broken gradient lines $\gamma_{un,s_i}\sharp\gamma_{s_i,sn}$. Suppose these separatrixes
are in the unbounded region $R_1$ determined by the square:
then, if, in the bounded region $R_2$, there are no other critical points,
$R_2$ is a component of ${\cal M}(un,sn)$, whose boundary is formed by
the broken gradient lines
$\gamma_{un,s_1}\sharp\gamma_{s_1,sn}$ and $\gamma_{un,s_2}\sharp\gamma_{s_2,sn}$.
Suppose instead that the two separatrixes lie in
$R_2$: then at least a stable and an unstable node, are
contained in it (see figure \ref{fig1}.1).
\end{example}

\begin{example}
\rm
If more than one square has $un$ and $sn$
among their vertexes, then ${\cal M}(un,sn)$ has more than one connected component.
\end{example}

\begin{example}
\rm
A case where ${\cal M}(un,sn)$ is not bounded is shown in
in figure \ref{fig1}.1):
in $\R^2$, ${\cal M}(un,sn)$ is the union of two unbounded connected sets $R_1^1$ and $R_1^2$,
such that $\partial R_1^i$ consists of a broken gradien line 
$\gamma_{un,s_i}\sharp\gamma_{s_i,sn}$
and of two gradient lines $\gamma_{un}$ and $\gamma_{sn}$, which, in the compactification
$S^2$ of $\R^2$, connect, respectively, $un$ and $sn$ to $\infty$;
so, in $S^2$,
the boundary of ${\cal M}(un,sn)$ consists of the two broken gradient lines
$\gamma_{un,s_1}\sharp\gamma_{s_1,sn}$ and $\gamma_{un,s_2}\sharp\gamma_{s_2,sn}$,
while $\gamma_{un}$ and $\gamma_{sn}$ are part of a single gradient line from $un$ to $sn$,
belonging to ${\cal M}(un,sn)$.
\end{example}

\begin{example}
\rm
Figure \ref{fig2}.2 shows two types of degenerate squares:\\ $(un_1,s,s,sn)$ bounded by 
$\gamma_1\sharp\gamma_3$ and $\gamma^{'}_1\sharp\gamma^{'}_3$
with $\gamma_1=\gamma^{'}_1$, and $(un_2,s,s,sn)$ bounded by 
$\gamma_2\sharp\gamma_3$ and $\gamma^{'}_2\sharp\gamma^{'}_3$
with $\gamma_2=\gamma^{'}_2$. In the first case $W^u(un_1)\cap W^s(sn)$
is bounded and connected, in the second case $W^u(un_2)\cap W^s(sn)$ is
not bounded and has two connected component, but it is bounded and connected
in the compactification $S^2$ of $\R^2$. Note that the square $(un_2,s,s,sn)$
requires the existence of the unstable node $un_1$ and so the existence of the
square $(un_1,s,s,sn)$, but not the 
converse. This lacking of symmetry is understood in $S^2$:
the square $(un_1,s,s,sn)$ requires the unstable node $un_2$. A similar example
is obtained exchanging the roles of stable and unstable nodes.
\end{example}

\begin{center}
\setlength{\unitlength}{1cm}
\begin{picture}(4,4)
\label{fig2}
\thinlines

\put(3,2){\circle*{.1}}
\put(2,2){\circle*{.1}}
\put(1,2){\circle*{.1}}
\put(.25,2){\circle*{.1}}

\scriptsize
\put(.9,2.1){$\displaystyle s$}
\put(2.1,2.1){$\displaystyle un_1$}
\put(.05,2.2){$\displaystyle un_2$}
\put(3.1,2.1){$\displaystyle sn$}
\put(1.3,1.65){$\displaystyle \gamma_1=\gamma^{'}_1$}
\put(0,1.65){$\displaystyle \gamma_2=\gamma^{'}_2$}
\put(1.9,3.2){$\displaystyle \gamma_3$}
\put(1.9,.5){$\displaystyle \gamma^{'}_3$}
\normalsize

\qbezier(1,2)(.5,2)(0,2)
\qbezier(1,2)(1.5,2)(2,2)
\qbezier(1,2)(2,4)(3,2)
\qbezier(1,2)(2,0)(3,2)
\qbezier(3,2)(3.5,2)(4,2)

\qbezier(1.5,2)(1.55,2.05)(1.6,2.1)
\qbezier(1.5,2)(1.55,1.95)(1.6,1.9)

\qbezier(3.75,2)(3.8,2.05)(3.85,2.1)
\qbezier(3.75,2)(3.8,1.95)(3.85,1.9)

\qbezier(.55,2.1)(.6,2.05)(0.65,2)
\qbezier(.55,1.9)(.6,1.95)(0.65,2)

\qbezier(.1,2.1)(.05,2.05)(0,2)
\qbezier(.1,1.9)(.05,1.95)(0,2)

\qbezier(2,3)(1.95,3.05)(1.9,3.1)
\qbezier(2,3)(1.95,2.95)(1.9,2.9)

\qbezier(2,1)(1.95,1.05)(1.9,1.1)
\qbezier(2,1)(1.95,0.95)(1.9,0.9)


\end{picture}

$Fig.~\ref{fig2}.2:~Degenerate~squares$
\end{center}

\begin{remark}
\rm
If $f$ is bounded or in a bounded subset of $\R^2$,
$f$ exhibits a finite number of critical points (in particular, this is true
if $f$ is defined on a compact manifold). 
In the sequel, $f$ is assumed to have only a finite
number of critical points.
\end{remark}

Under the hypothesis of theorem \ref{jost}, it can be constructed a complex, the Morse complex,
whose homology groups are used to define the mirror object in Mirror Symmetry.
Very briefly, given a Morse function satisfying the hypothesis of theorem \ref{jost} 
and a metric, the Morse complex is defined by the free module (in our case over $\C$) generated by
critical points of the Morse function. The grading is given by
the Morse index. The Morse differential is defined by counting gradient lines between critical points:
\begin{displaymath}
0\rightarrow\C[un_i]\rightarrow^{\partial_0}\C[s_j]\rightarrow^{\partial_1}\C[sn_k]\rightarrow0
\end{displaymath}
here $un_i$ denotes the unstable nodes, $s_j$ the saddles and $sn_k$ the stable nodes.
The construction of the Morse differential, when we consider coefficients in $\C$,
depends on the choice of an orientation of gradient lines between critical points.
That this choice can be done in a compatible way, that is, the orientations
of gradient lines between two critical points $p$ and $q$ such that
$|\mu(p)-\mu(q)|$=2 and the induced orientations on the gradient lines
(between critical points of relative Morse index 1)
in the boundary
of $W^s(p)\cap W^u(q)$ are consistent, allowing to talk of a
``coherent orientation'', can be proved for finite dimensional oriented manifolds
(see \cite{S}).
In this paper the sign of a gradient lines $\gamma$
is denoted by $n(\gamma)$.

In an admissible $CB$-diagram $(C,B,(U_i,P_i)_{i\in I})$, to each point $x\notin C\cup B$, it is associated a gradient vector
field $\nabla f_x$ with phase portrait $P_i$: this enables to define a Morse complex
$$\C[p_{1}(x)]\oplus ... \oplus\C[p_{k}(x)]$$
where $p_i(x)$ are the critical points of $\nabla f_x$, and Morse homology groups. Since for all $x\in U_i$, the vector
fields $\nabla f_x$ are orbit equivalent, it follows that the Morse complexes associated to them are isomorphic: this allows
to write
$$\C[p_{1}^{U_i}]\oplus ... \oplus\C[p_{k}^{U_i}]$$
and to talk about the Morse complex and Morse homology over $U_i$.

The following lemma, besides stating the isomorphism between Morse homology over $U_1$ and $U_2$, when
their boundaries are folds not limit points of $B$, suggests also how to pick up an isomorphism, leading to
the definition of quantum corrections.

\begin{lemma}
\label{causdiag}
Let $U_1$ and $U_2$ be open subsets such that $U_i\cap(C\cup B)=\varnothing$
and $\partial U_1\cap\partial U_2\subset C$ consists only of folds not limit
points of $B$; let $p_1(x')$, ... , $p_k(x')$, $p_{k+1}(x')$, $p_{k+2}(x')$ be
the critical points of $\nabla f_{x'}$ for $x'\in U_1$,
and $p_1(x')$, ... , $p_k(x')$
be the critical points of $\nabla f_{x'}$ for $x'\in U_2$,
so that $(p_{k+1}(x'),p_{k+2}(x'))$
is a birth-death pair.
Then the homology groups of the Morse complexes 
$$\C[p_{1}^{U_2}]\oplus ... \oplus\C[p_{k}^{U_2}]$$
$$\C[p_{1}^{U_1}]\oplus ... \oplus\C[p_k^{U_1}]
\oplus\C[p_{k+1}^{U_1}]\oplus\C[p_{k+2}^{U_1}]$$
are isomorphic.
\end{lemma}
\begin{proof}
\begin{enumerate}
\item We confront first $HM_1(U_1)$ with $HM_1(U_2)$, by computing $Im\partial_0$ and 
$Ker\partial_1$ in $U_1$ and $U_2$.

\begin{itemize}
\item
Suppose, first, that the birth-death pair $(p_{k+1},p_{k+2})$ 
is of type $(n,s)$, where $n$ is an unstable node and $s$ a saddle.
It will be proved that $dimIm\partial_0^{U_1}=dimIm\partial_0^{U_2}+1$ and
$dimKer\partial_1^{U_1}=dimKer\partial_1^{U_2}+1$.

Consider the phase portraits, looking at figure \ref{fig1}.3.
For $x\in U_1$, by lemma \ref{exgrl}, there exists
a gradient line $\gamma_{n,s}$ from $n$ to $s$, which is, therefore, 
one of the two components of $W^s(s)$. This implies that
$s$ may be connected, besides to $n$, to at most a second unstable node $un$, by means of
the second component of $W^s(s)$. Moreover, it may also be connected to at most two
unstable nodes, by means of  
its remaining two separatrixes forming $W^u(s)$.
As to $n$, the gradient lines forming $W^u(n)$ may connect $n$ to further
saddles or stable nodes. For $x\in C$, $n$ and $s$ glue together
in a birth-death point $d$: this exhibits an unstable manifold $W^u(d)$
and a center manifold $W^c(d)$ (see \cite{GH}), corresponding to the zero eigenvalue; observe
that $d$ is connected to $un$ by a gradient line in $W^c(d)$, to all saddles and
stable nodes connected to $n$ in $U_1$, by means of other gradient lines in $W^c(d)$, and to
the stable nodes connected to $s$ in $U_1$, by means of the two components of $W^u(d)$. 
In $U_2$, the center manifold $W^c(d)$ breaks, forming gradient lines
from $un$ to all the stable nodes and saddles connected to $n$ and $s$.

\begin{center}
\setlength{\unitlength}{1cm}
\begin{picture}(11,4)
\label{fig3}
\thinlines

\put(2,2){\circle*{.1}}
\put(3,2){\circle*{.1}}
\put(6,2){\circle*{.1}}

\scriptsize
\put(2,2.1){$\displaystyle n$}
\put(3.1,2.1){$\displaystyle s$}
\put(6.1,2.1){$\displaystyle d$}

\put(3.1,2.8){$\displaystyle (sn)$}
\put(1,2.4){$\displaystyle (sn)$}
\put(1,1.4){$\displaystyle (s)$}
\put(3.6,2.1){$\displaystyle (un)$}
\put(2.5,.5){$\displaystyle U_1$}

\put(6.1,2.8){$\displaystyle (sn)$}
\put(5,2.4){$\displaystyle (sn)$}
\put(5,1.4){$\displaystyle (s)$}
\put(6.6,2.1){$\displaystyle (un)$}
\put(6,.5){$\displaystyle C$}

\put(8,2.1){$\displaystyle (sn)$}
\put(8,1.7){$\displaystyle (s)$}
\put(9.6,2.3){$\displaystyle (un)$}
\put(9,.5){$\displaystyle U_2$}

\normalsize

\qbezier(1,2)(3,2)(4,2)
\qbezier(3,1)(3,2)(3,3)
\qbezier(2,2)(2.5,2)(2.5,3)
\qbezier(2,2)(2.5,2)(2.5,1)
\qbezier(2,2)(1.5,2)(1.5,1)
\qbezier(2,2)(1.5,2)(1.5,3)

\qbezier(5,2)(6,2)(7,2)
\qbezier(6,1)(6,2)(6,3)
\qbezier(6,2)(5.5,2)(5.5,1)
\qbezier(6,2)(5.5,2)(5.5,3)

\qbezier(8,2)(9,2)(10,2)
\qbezier(8,3)(8,2.1)(10,2.1)
\qbezier(8,1)(8,1.9)(10,1.9)

\qbezier(1.5,2)(1.55,2.05)(1.6,2.1)
\qbezier(1.5,2)(1.55,1.95)(1.6,1.9)

\qbezier(3.5,2)(3.55,2.05)(3.6,2.1)
\qbezier(3.5,2)(3.55,1.95)(3.6,1.9)

\qbezier(2.5,2)(2.45,2.05)(2.4,2.1)
\qbezier(2.5,2)(2.45,1.95)(2.4,1.9)

\qbezier(1.6,2.3)(1.65,2.3)(1.7,2.3)
\qbezier(1.6,2.3)(1.6,2.25)(1.6,2.2)

\qbezier(1.6,1.7)(1.65,1.7)(1.7,1.7)
\qbezier(1.6,1.7)(1.6,1.75)(1.6,1.8)

\qbezier(2.4,2.3)(2.35,2.3)(2.3,2.3)
\qbezier(2.4,2.3)(2.4,2.25)(2.4,2.2)

\qbezier(2.4,1.7)(2.35,1.7)(2.3,1.7)
\qbezier(2.4,1.7)(2.4,1.75)(2.4,1.8)

\qbezier(3,1.5)(3.05,1.55)(3.1,1.6)
\qbezier(3,1.5)(2.95,1.55)(2.9,1.6)

\qbezier(3,2.5)(3.05,2.45)(3.1,2.4)
\qbezier(3,2.5)(2.95,2.45)(2.9,2.4)

\qbezier(5.5,2)(5.55,2.05)(5.6,2.1)
\qbezier(5.5,2)(5.55,1.95)(5.6,1.9)

\qbezier(5.6,2.3)(5.65,2.3)(5.7,2.3)
\qbezier(5.6,2.3)(5.6,2.25)(5.6,2.2)

\qbezier(5.6,1.7)(5.65,1.7)(5.7,1.7)
\qbezier(5.6,1.7)(5.6,1.75)(5.6,1.8)

\qbezier(6,1.5)(6.05,1.55)(6.1,1.6)
\qbezier(6,1.5)(5.95,1.55)(5.9,1.6)

\qbezier(6,2.5)(6.05,2.45)(6.1,2.4)
\qbezier(6,2.5)(5.95,2.45)(5.9,2.4)

\qbezier(6.5,2)(6.55,2.05)(6.6,2.1)
\qbezier(6.5,2)(6.55,1.95)(6.6,1.9)

\qbezier(9,2)(9.05,2.05)(9.1,2.1)
\qbezier(9,2)(9.05,1.95)(9.1,1.9)

\qbezier(8.6,2.3)(8.65,2.325)(8.7,2.35)
\qbezier(8.6,2.3)(8.6,2.25)(8.6,2.2)

\qbezier(8.6,1.7)(8.65,1.675)(8.7,1.65)
\qbezier(8.6,1.7)(8.6,1.75)(8.6,1.8)


\end{picture}
$Fig.~\ref{fig1}.3:~Phase~portraits~near~a~fold$
\end{center}

Observe, therefore, that if there are, in $U_1$, a second unstable node $un$ connected to
$s$ and saddles $s_i$ connected to $n$, then in $U_2$ the new gradient lines from $un$ to
$s_i$, not appearing in $U_1$, implies a change in $Im\partial_0^{U_2}$
with respect to $Im\partial_0^{U_1}$.

\item[*]
So consider $Im\partial_0$.
Denote by
$un_1$, ..., $un_n$ the unstable nodes appearing in the phase portrait over $U_2$.
We distinguish two cases:

\item[**]
If no unstable node, except $n$, is connected to $s$ in $U_1$,
then
$$Im\partial_0^{U_1}=i(Im\partial_0^{U_2})\oplus<\partial_0^{U_1}(n)>$$
where $i$ denotes the natural injection, defined by the continuity
of the family $f_x$ in the parameter $x$,
$$\C[p_{1}^{U_2}]\oplus ... \oplus\C[p_{k}^{U_2}]\hookrightarrow
\C[p_{1}^{U_1}]\oplus ... \oplus\C[p_k^{U_1}]
\oplus\C[p_{k+1}^{U_1}]\oplus\C[p_{k+2}^{U_1}]$$
indeed, $\partial_0^{U_1}(n)\nsubseteq i(Im\partial_0^{U_2})$
because writing $\partial_0^{U_1}(n)=s_{i_1}+...+s_{i_l}+s$, for some saddles
$s_{i_1}$, ..., $s_{i_l}$, by hypothesis, $s\notin i(Im\partial_0^{U_2})$.
Note that it should rather be written $n(\gamma_{n,s_{i_1}})s_{i_1}+...+n(\gamma_{n,s_{i_l}})s_{i_l}
+n(\gamma_{n,s})s$, where $n(\gamma_{n,s_{i_j}})$, $n(\gamma_{n,s}) \in\{1,-1\}$ 
depend on the choice of an orientation of gradient lines; however, for
simplicity, it will often assumed, when possible, $n(\gamma_{n,s_{i_j}})=n(\gamma_{n,s})=1$. 

\item[**]
Suppose now $un_k$ is a second unstable node (denoted by $un$ in figure \ref{fig1}.3) connected to $s$
by the gradient line $\gamma_{un_k,s}$:
if $\partial_0^{U_1}(un_k)=<s_{k_1}+...+s_{k_m}-s>$ for some
saddles $s_{k_1}$, ..., $s_{k_m}$ (note here the choice of signs: in fact, the two
components $\gamma_{un_k,s}$ and $\gamma_{n,s}$ of $W^s(s)$ have opposite
orientations), and, as above, $\partial_0^{U_1}(n)=s_{i_1}+...+s_{i_l}+s$, then
$\partial_0^{U_2}(un_k)=<s_{k_1}+...+s_{k_m}+s_{i_1}+...+s_{i_l}>$
and so
\begin{eqnarray}
Im\partial_0^{U_1} &= & i(\partial_0^{U_2}(<un_1,...,un_{k-1},un_{k+1},...,un_n>))
\oplus
\nonumber\\
& & \oplus<\partial_0^{U_1}(un_k)>\oplus<\partial_0^{U_1}(n)>=
\nonumber\\
& = & i(\partial_0^{U_2}(<un_1,...,un_{k-1},un_{k+1},...,un_n>))\oplus
\nonumber\\
& & \oplus<s_{k_1}+...+s_{k_m}-s>\oplus<s_{i_1}+...+s_{i_l}+s>=
\nonumber\\
& = & i(\partial_0^{U_2}(<un_1,...,un_{k-1},un_{k+1},...,un_n>))\oplus
\nonumber\\
& & \oplus<s_{k_1}+...+s_{k_m}+s_{i_1}+...+s_{i_l}>\oplus
\nonumber\\
& & \oplus<s_{i_1}+...+s_{i_l}+s>=
\nonumber\\
& = & i(\partial_0^{U_2}(<un_1,...,un_{k-1},un_{k+1},...,un_n>))\oplus 
\nonumber\\
& & \oplus i(\partial_0^{U_2}(<un_k>))
\oplus<\partial_0^{U_1}(n)>=
\nonumber\\
& = & i(Im\partial_0^{U_2})\oplus<\partial_0^{U_1}(n)>
\nonumber
\end{eqnarray}
Note that if it were $\partial_0^{U_1}(n)=<s_{i_1}+...+s_{i_l}+s>
=<s_{k_1}+...+s_{k_m}-s>=\partial_0^{U_1}(un_k)$,
then $l=m$ and $s_{i_r}=s_{k_r}$ for all $1\leq r\leq l=m$,
thus $\partial_0^{U_2}(un_k)=0$: this means that $un_k$ is connected to
all the saddles $s_{i_r}=s_{k_r}$ by both the separatrixes forming $W^s(s_{i_r})$
(this is shown, reversing the roles of
stable and unstable nodes, in figure \ref{fig2}.2 representing a degenerate square). 
Thus the relation $Im\partial_0^{U_1}=i(Im\partial_0^{U_2})\oplus<\partial_0^{U_1}(n)>$
still holds.
To better llustrate this case, consider only two saddles,
$s$ and $s_g$, such that
$\partial_0^{U_1}(n)=<s_g+s>
=$ $<s_g-s>=\partial_0^{U_1}(un_k)$, and let $\gamma_{n,s}$, $\gamma_{n,s_g}$
and $\gamma_{un_k,s}$, $\gamma_{un_k,s_g}$ be the gradient lines connecting $s$ and $s_g$,
respectively, to $n$ and $un_k$, as in figure \ref{fig1}.4.

\begin{center}
\setlength{\unitlength}{1cm}
\begin{picture}(4,4)
\label{fig4}
\thinlines

\put(3.5,2){\circle*{.1}}
\put(2,.5){\circle*{.1}}
\put(.5,2){\circle*{.1}}
\put(2,3.5){\circle*{.1}}

\scriptsize
\put(2,3.6){$\displaystyle s$}
\put(0,1.8){$\displaystyle un_k$}
\put(2,.2){$\displaystyle s_g$}
\put(3.6,2){$\displaystyle n$}
\normalsize

\qbezier(3.5,2)(2.75,1.25)(2,.5)
\qbezier(3.5,2)(2.75,2.75)(2,3.5)
\qbezier(.5,2)(1.25,1.25)(2,.5)
\qbezier(.5,2)(1.25,2.75)(2,3.5)

\qbezier(2.75,2.75)(2.76,2.75)(2.85,2.75)
\qbezier(2.75,2.75)(2.75,2.74)(2.75,2.65)

\qbezier(1.25,2.75)(1.2,2.75)(1.15,2.75)
\qbezier(1.25,2.75)(1.25,2.74)(1.25,2.65)

\qbezier(2.75,1.25)(2.76,1.25)(2.85,1.25)
\qbezier(2.75,1.25)(2.75,1.3)(2.75,1.35)

\qbezier(1.25,1.25)(1.2,1.25)(1.15,1.25)
\qbezier(1.25,1.25)(1.25,1.3)(1.25,1.35)


\end{picture}

$Fig.~\ref{fig1}.4:~\partial_0^{U_1}(n)=\partial_0^{U_1}(un_k)$
\end{center}

One of the two components of both $W^u(s)$ and $W^u(s_g)$ lies in the bounded
region $R_2$ determined by $\gamma_{n,s}$, 
$\gamma_{n,s_g}$, $\gamma_{un_k,s}$ and $\gamma_{un_k,s_g}$, so,
since the vector fields are gradient,
there must be at least another critical point in $R_2$: this must be a stable
node $sn_g$, and, for simplicity, assume this is the only critical point, as
shown in figure \ref{fig1}.5).

\begin{center}
\setlength{\unitlength}{1cm}
\begin{picture}(4,4)
\label{fig5}
\thinlines

\put(3.5,2){\circle*{.1}}
\put(2,.5){\circle*{.1}}
\put(.5,2){\circle*{.1}}
\put(2,3.5){\circle*{.1}}
\put(2,2){\circle*{.1}}

\scriptsize
\put(2,3.6){$\displaystyle s$}
\put(0,1.8){$\displaystyle un_k$}
\put(2,.2){$\displaystyle s_g$}
\put(3.6,2){$\displaystyle n$}
\put(2.1,2.1){$\displaystyle sn_g$}
\normalsize

\qbezier(3.5,2)(2.75,1.25)(2,.5)
\qbezier(3.5,2)(2.75,2.75)(2,3.5)
\qbezier(.5,2)(1.25,1.25)(2,.5)
\qbezier(.5,2)(1.25,2.75)(2,3.5)

\qbezier(.5,2)(1,2)(2,2)
\qbezier(3.5,2)(3,2)(2,2)
\qbezier(2,.5)(2,1)(2,2)
\qbezier(2,3.5)(2,3)(2,2)

\qbezier(2.75,2.75)(2.76,2.75)(2.85,2.75)
\qbezier(2.75,2.75)(2.75,2.74)(2.75,2.65)

\qbezier(1.25,2.75)(1.2,2.75)(1.15,2.75)
\qbezier(1.25,2.75)(1.25,2.74)(1.25,2.65)

\qbezier(2.75,1.25)(2.76,1.25)(2.85,1.25)
\qbezier(2.75,1.25)(2.75,1.3)(2.75,1.35)

\qbezier(1.25,1.25)(1.2,1.25)(1.15,1.25)
\qbezier(1.25,1.25)(1.25,1.3)(1.25,1.35)

\qbezier(2.75,2)(2.8,2.05)(2.85,2.1)
\qbezier(2.75,2)(2.8,1.95)(2.85,1.9)

\qbezier(1.25,2)(1.2,2.05)(1.15,2.1)
\qbezier(1.25,2)(1.2,1.95)(1.15,1.9)

\qbezier(2,1.25)(2.05,1.2)(2.1,1.15)
\qbezier(2,1.25)(1.95,1.2)(1.9,1.15)

\qbezier(2,2.75)(2.05,2.8)(2.1,2.85)
\qbezier(2,2.75)(1.95,2.8)(1.9,2.85)


\end{picture}

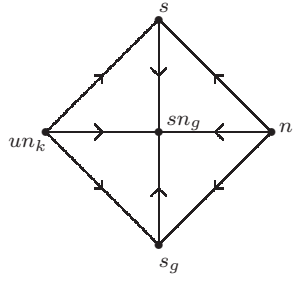
$Fig.~\ref{fig1}.5:~\partial_0^{U_1}(n)=\partial_0^{U_1}(un_k):~the~phase~portrait~in~U_1$
\end{center}

in $U_2$, since no other saddles were assumed to be
in $R^2$, by lemma \ref{php},
a degenerate square with vertexes $un_k$, $s_g$ and $sn_g$
appears in the phase portrait, as shown in figure \ref{fig1}.6.

\begin{center}
\setlength{\unitlength}{1cm}
\begin{picture}(4,4)
\label{fig6}
\thinlines

\put(2,.5){\circle*{.1}}
\put(.5,2){\circle*{.1}}
\put(2,2){\circle*{.1}}

\scriptsize
\put(0,1.8){$\displaystyle un_k$}
\put(2,.2){$\displaystyle s_g$}
\put(2.1,2.1){$\displaystyle sn_g$}
\normalsize

\qbezier(3.5,2)(2.75,1.25)(2,.5)
\qbezier(3.5,2)(2.75,2.75)(2,3.5)
\qbezier(.5,2)(1.25,1.25)(2,.5)
\qbezier(.5,2)(1.25,2.75)(2,3.5)

\qbezier(.5,2)(1,2)(2,2)
\qbezier(3.1,2)(3,2)(2,2)
\qbezier(2,.5)(2,1)(2,2)
\qbezier(2,3.1)(2,3)(2,2)
\qbezier(2,3.1)(1.25,2.55)(.5,2)
\qbezier(2,3.3)(1.25,2.65)(.5,2)
\qbezier(2,3.3)(2.55,2.65)(3.1,2)

\qbezier(2.75,2.75)(2.7,2.75)(2.65,2.75)
\qbezier(2.75,2.75)(2.75,2.8)(2.75,2.85)

\qbezier(1.25,2.75)(1.2,2.75)(1.15,2.75)
\qbezier(1.25,2.75)(1.25,2.74)(1.25,2.65)

\qbezier(2.75,1.25)(2.8,1.25)(2.85,1.25)
\qbezier(2.75,1.25)(2.75,1.3)(2.75,1.35)

\qbezier(1.25,1.25)(1.2,1.25)(1.15,1.25)
\qbezier(1.25,1.25)(1.25,1.3)(1.25,1.35)

\qbezier(2.75,2)(2.8,2.05)(2.85,2.1)
\qbezier(2.75,2)(2.8,1.95)(2.85,1.9)

\qbezier(1.25,2)(1.2,2.05)(1.15,2.1)
\qbezier(1.25,2)(1.2,1.95)(1.15,1.9)

\qbezier(2,1.25)(2.05,1.2)(2.1,1.15)
\qbezier(2,1.25)(1.95,1.2)(1.9,1.15)

\qbezier(2,2.75)(2.05,2.8)(2.1,2.85)
\qbezier(2,2.75)(1.95,2.8)(1.9,2.85)


\end{picture}

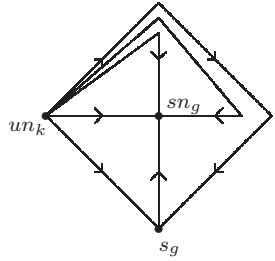
$Fig.~\ref{fig1}.6:~\partial_0^{U_1}(n)=\partial_0^{U_1}(un_k):~the~phase~portrait~in~U_2$
\end{center}

\item[*]
It remains to prove that $dimKer\partial_1^{U_1}=dimKer\partial_1^{U_2}+1$.

\item[**]
If $\partial_1^{U_1}(s)=0$, that is $s$ is not connected to any 
stable node, then 
$$Ker\partial_1^{U_1}=i(Ker\partial_1^{U_2})\oplus<s>$$
since $s\notin i(Ker\partial_1^{U_2})$. 

\item[**]
If, instead, $\partial_1^{U_1}(s)$ does not
vanish, then either there exists a stable node $sn_g$ such that 
$\partial_1^{U_1}(s)=sn_g$, or there exist two stable nodes $sn_{g}$
and $sn_{h}$ such that $\partial_1^{U_1}(s)=sn_{g}+ sn_{k}$
(as already explained, the signs in $\partial_1^{U_i}$ depends on the orientation of gradient lines,
however, it is done this choice to simplify the notation). 

\item[***]
If $\partial_1^{U_2}(s)=sn_g$ then, by lemma \ref{php}, there exists a saddle $s_g$
forming a square together with $n$, $s$ and $sn_g$.
Assume that the square
is non-degenerate, that is $s_g\neq s$. Suppose also, for
the moment, that $s_g$ is the
only saddle, besides $s$, connected to $sn_g$. 

\item[****]
Suppose that only an unstable node, that is, $n$, is connected to $s$

\item[*****]
If
$s_g$ is not connected to a second stable node, then  
$$Ker\partial_1^{U_1}=i(Ker\partial_1^{U_2})\oplus<s\pm s_g>$$ 
where the sign 
depends on the orientation of gradient lines,
particularly on the orientation of the square $(n,s,s_g,sn_g)$. 

\item[*****]
If, instead,
$s_g$ is connected to a second stable node $sn_{g_1}$,
as figure \ref{fig1}.7 shows, then by lemma \ref{php}
there exists a saddle $s_{g_1}$, which suppose, for the moment, distinct
from $s$ and $s_g$, forming a square together with $n$, $s_g$
and $sn_{g_1}$; repeating this argument, being finite the number of critical points,
it follows that there is at most a finite number $k$ of such squares having
as vertexes the node $n$, the saddles $s_{g_i}$ and $s_{g_{i+1}}$, and
the stable node $sn_{g_{i+1}}$, as shown in figure \ref{fig1}.7.

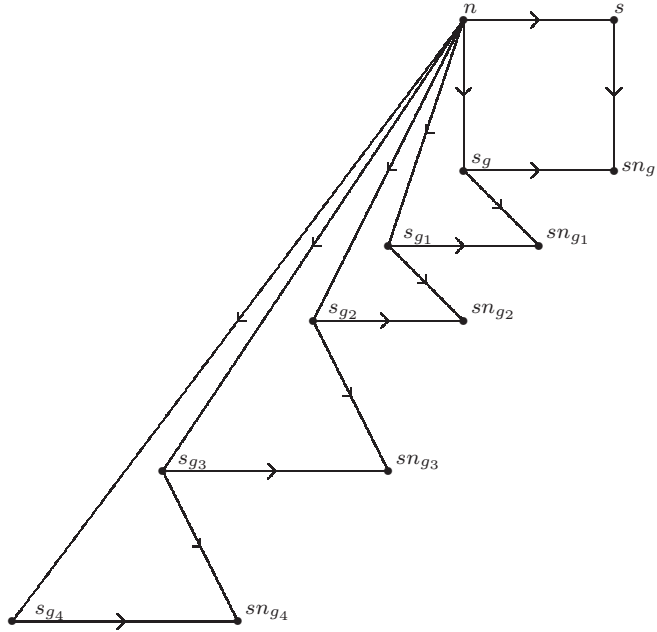
\begin{center}
\setlength{\unitlength}{1cm}
\begin{picture}(9,9)
\label{fig7}
\thinlines

\put(8,8){\circle*{.1}}
\put(6,8){\circle*{.1}}
\put(8,6){\circle*{.1}}
\put(6,6){\circle*{.1}}
\put(5,5){\circle*{.1}}
\put(7,5){\circle*{.1}}
\put(4,4){\circle*{.1}}
\put(6,4){\circle*{.1}}
\put(2,2){\circle*{.1}}
\put(5,2){\circle*{.1}}
\put(0,0){\circle*{.1}}
\put(3,0){\circle*{.1}}

\scriptsize
\put(8,8.1){$\displaystyle s$}
\put(6,8.1){$\displaystyle n$}
\put(8.1,6){$\displaystyle sn_g$}
\put(6.1,6.1){$\displaystyle s_g$}
\put(5.2,5.1){$\displaystyle s_{g_1}$}
\put(7.1,5.1){$\displaystyle sn_{g_1}$}
\put(4.2,4.1){$\displaystyle s_{g_2}$}
\put(6.1,4.1){$\displaystyle sn_{g_2}$}
\put(2.2,2.1){$\displaystyle s_{g_3}$}
\put(5.1,2.1){$\displaystyle sn_{g_3}$}
\put(3.1,.1){$\displaystyle sn_{g_4}$}
\put(.3,.1){$\displaystyle s_{g_4}$}

\normalsize

\qbezier(8,8)(7,8)(6,8)
\qbezier(6,8)(6,7)(6,6)
\qbezier(6,6)(7,6)(8,6)
\qbezier(8,8)(8,7)(8,6)

\qbezier(6,8)(5.5,6.5)(5,5)
\qbezier(7,5)(6,5)(5,5)
\qbezier(7,5)(6.5,5.5)(6,6)

\qbezier(6,8)(5,6)(4,4)
\qbezier(6,4)(5,4)(4,4)
\qbezier(6,4)(5.5,4.5)(5,5)

\qbezier(6,8)(4,5)(2,2)
\qbezier(5,2)(3.5,2)(2,2)
\qbezier(4,4)(4.5,3)(5,2)

\qbezier(6,8)(3,4)(0,0)
\qbezier(3,0)(1.5,0)(0,0)
\qbezier(2,2)(2.5,1)(3,0)


\qbezier(7,8)(6.95,8.05)(6.9,8.1)
\qbezier(7,8)(6.95,7.95)(6.9,7.9)

\qbezier(7,6)(6.95,6.05)(6.9,6.1)
\qbezier(7,6)(6.95,5.95)(6.9,5.9)

\qbezier(6,7)(6.05,7.05)(6.1,7.1)
\qbezier(6,7)(5.95,7.05)(5.9,7.1)

\qbezier(8,7)(8.05,7.05)(8.1,7.1)
\qbezier(8,7)(7.95,7.05)(7.9,7.1)

\qbezier(5.5,6.5)(5.5,6.55)(5.5,6.6)
\qbezier(5.5,6.5)(5.55,6.5)(5.6,6.5)

\qbezier(6,5)(5.95,5.05)(5.9,5.1)
\qbezier(6,5)(5.95,4.95)(5.9,4.9)

\qbezier(6.5,5.5)(6.5,5.55)(6.5,5.6)
\qbezier(6.5,5.5)(6.45,5.5)(6.4,5.5)

\qbezier(5,6)(5,6.05)(5,6.1)
\qbezier(5,6)(5.05,6)(5.1,6)

\qbezier(5,4)(4.95,4.05)(4.9,4.1)
\qbezier(5,4)(4.95,3.95)(4.9,3.9)

\qbezier(5.5,4.5)(5.5,4.55)(5.5,4.6)
\qbezier(5.5,4.5)(5.45,4.5)(5.4,4.5)

\qbezier(4,5)(4,5.05)(4,5.1)
\qbezier(4,5)(4.05,5)(4.1,5)

\qbezier(3.5,2)(3.45,2.05)(3.4,2.1)
\qbezier(3.5,2)(3.45,1.95)(3.4,1.9)

\qbezier(4.5,3)(4.5,3.05)(4.5,3.1)
\qbezier(4.5,3)(4.45,3)(4.4,3)

\qbezier(3,4)(3,4.05)(3,4.1)
\qbezier(3,4)(3.05,4)(3.1,4)

\qbezier(1.5,0)(1.45,.05)(1.4,.1)
\qbezier(1.5,0)(1.45,-.05)(1.4,-.1)

\qbezier(2.5,1)(2.5,1.05)(2.5,1.1)
\qbezier(2.5,1)(2.45,1)(2.4,1)

\end{picture}

$Fig.~\ref{fig1}.7:~s_g~is~connected~to~a~second~stable~node$
\end{center}

From figure \ref{fig1}.7, it turns out now that 
$$Ker\partial_1^{U_1}=i(Ker\partial_1^{U_2})\oplus
<s\pm s_g\pm s_{g_1}\pm...\pm s_{g_k}>$$ . 

\item[*****]
Consider now the cases
$s_{g_1}=s$ and $s_{g_1}=s_g$. If $s_{g_1}=s$, there are
two possible phase portraits (one of which, to be the phase portrait of a gradient
vector field, requires in $U_1$ a further unstable node and
yields in $U_2$ a degenerate square): however, in both cases, $Im\partial_0^{U_1}=
Ker\partial_1^{U_1}=\C$ and $Im\partial_0^{U_2}= Ker\partial_1^{U_2}=\{0\}$.
If $s_{g_1}=s_g$ then $W^u(n)\cap W^s(sn)$ is not bounded and the phase portrait
can not be extended to the compactification $S^2$
of $\R^2$ unless adding further critical points:
such phase portrait contradicts the compactness assumption \ref{assumption}
and so is not considered. 

\item[*****]
If $sn_g$ is connected to other saddles,
$dimKer\partial_1^{U_1}=dimKer\partial_1^{U_2}+1$ is still valid.
Indeed, observe that two cases may occur:
starting from $sn_g$, consider chains formed by a finite alternate 
sequence of saddles and stable nodes,
ending with a saddle, as $(sn_g,s_{g_1},sn_{g_1},s_{g_2},sn_{g_2},...,sn_{g_{k-1}},s_{g_k})$, 
or with a stable node, as $(sn_g,s_{g_1},sn_{g_1},s_{g_2},sn_{g_2},...,s_{g_j},sn_{g_j})$,
for some finite $k$ or $j$ (two chains may have elements in common). 
In figure \ref{fig1}.8 the two cases are represented by the chain $(sn_g,s_m)$
and by the chain $(sn_g,s_{g_1},sn_{g_1})$. Referring to this situation and
supposing for example that
$\gamma_{s,sn_g}$, $\gamma_{s_m,sn_g}$ and $\gamma_{s_{g_1},sn_g}$ 
have positive orientation,
and $\gamma_{s_g,sn_g}$ has negative orientation,
it follows that $s_m+s_g\in Ker\partial_1^{U_2}$, and that 
$s_m+s_g$, $s+s_g$, $s-s_m\in Ker\partial_1^{U_1}$, and
thus $Ker\partial_1^{U_1}$ is generated by two of these elements. 
The same argument works
when considering longer chains.
Note also that when a chain ends with a stable node,
a linear combination of the saddles of the chain with other
saddles in the phase portrait never belongs to $Ker\partial_1^{U_1}$ or to
$Ker\partial_1^{U_2}$.

\begin{center}
\setlength{\unitlength}{1cm}
\begin{picture}(6,4)
\label{fig8}
\thinlines

\put(1,2){\circle*{.1}}
\put(1,3){\circle*{.1}}
\put(3,2){\circle*{.1}}
\put(3,3){\circle*{.1}}
\put(3,1){\circle*{.1}}
\put(4,1){\circle*{.1}}
\put(5,1){\circle*{.1}}

\scriptsize
\put(1.1,2.1){$\displaystyle s_g$}
\put(.7,3){$\displaystyle n$}
\put(3.1,3.1){$\displaystyle s$}
\put(3.1,2.1){$\displaystyle sn_g$}
\put(3,.7){$\displaystyle s_m$}
\put(4,.7){$\displaystyle s_{g_1}$}
\put(5,.7){$\displaystyle sn_{g_1}$}

\normalsize

\qbezier(1,2)(1,2.5)(1,3)
\qbezier(3,2)(3,2.5)(3,3)
\qbezier(1,2)(2,2)(3,2)
\qbezier(1,3)(2,3)(3,3)
\qbezier(3,3)(4,3)(5,3)
\qbezier(3,1)(3,1.5)(3,2)
\qbezier(4,1)(3.5,1.5)(3,2)
\qbezier(4,1)(4.5,1)(5,1)


\qbezier(1,2.5)(1.05,2.55)(1.1,2.6)
\qbezier(1,2.5)(.95,2.55)(.9,2.6)
\qbezier(3,2.5)(3.05,2.55)(3.1,2.6)
\qbezier(3,2.5)(2.95,2.55)(2.9,2.6)
\qbezier(2,3)(1.95,3.05)(1.9,3.1)
\qbezier(2,3)(1.95,2.95)(1.9,2.9)
\qbezier(2,2)(1.95,2.05)(1.9,2.1)
\qbezier(2,2)(1.95,1.95)(1.9,1.9)
\qbezier(4,3)(4.05,3.05)(4.1,3.1)
\qbezier(4,3)(4.05,2.95)(4.1,2.9)
\qbezier(3,1.5)(3.05,1.45)(3.1,1,4)
\qbezier(3,1.5)(2.95,1.45)(2.9,1,4)
\qbezier(3.5,1.5)(3.5,1.45)(3.5,1.4)
\qbezier(3.5,1.5)(3.55,1.5)(3.6,1.5)
\qbezier(4.5,1)(4.45,1.05)(4.4,1.1)
\qbezier(4.5,1)(4.45,.95)(4.4,.9)

\end{picture}

$Fig.~\ref{fig1}.8:~s~is~connected~to~a~stable~node$
\end{center}

\item[*****]
It remains to consider the case where $s_g=s$. If $s$ has a double connection
to $sn$ (as explained, it can occur into two ways, in one of which, a second
unstable node is required in order to make the phase portrait that of
a gradient vector field: this 
anticipates the situation, which will be considered soon, where $s$ is connected
to two unstable nodes), then $Ker\partial_1^{U_1}=<s>$ and
$Ker\partial_1^{U_2}=\{0\}$. If, instead, $s$ has in $U_1$ a double connection to $n$,
a homoclinic orbit appears in the phase portrait
along $C$: this can not occur for a gradient vector field. 

%

\item[****]
Suppose now that $s$ is connected to a second unstable node $un$:
then by lemma \ref{php} there exists a saddle $s_j$ forming a square
together with $un$, $s$ and $sn_g$, as figure \ref{fig1}.9 
shows (suppose for the moment $s_j\neq s_g$
and $s_j\neq s$).

\begin{center}
\setlength{\unitlength}{1cm}
\begin{picture}(6,4)
\label{fig9}
\thinlines

\put(1,2){\circle*{.1}}
\put(1,3){\circle*{.1}}
\put(3,2){\circle*{.1}}
\put(3,3){\circle*{.1}}
\put(3,1){\circle*{.1}}
\put(4,1){\circle*{.1}}
\put(5,1){\circle*{.1}}
\put(5,3){\circle*{.1}}
\put(5,2){\circle*{.1}}

\scriptsize
\put(1.1,2.1){$\displaystyle s_g$}
\put(.7,3){$\displaystyle n$}
\put(3.1,3.1){$\displaystyle s$}
\put(3.1,2.1){$\displaystyle sn_g$}
\put(3,.7){$\displaystyle s_m$}
\put(4,.7){$\displaystyle s_{g_1}$}
\put(5,.7){$\displaystyle sn_{g_1}$}
\put(5.1,3){$\displaystyle un$}
\put(5,1.8){$\displaystyle s_{j}$}

\normalsize

\qbezier(1,2)(1,2.5)(1,3)
\qbezier(3,2)(3,2.5)(3,3)
\qbezier(1,2)(2,2)(3,2)
\qbezier(1,3)(2,3)(3,3)
\qbezier(3,3)(4,3)(5,3)
\qbezier(3,1)(3,1.5)(3,2)
\qbezier(4,1)(3.5,1.5)(3,2)
\qbezier(4,1)(4.5,1)(5,1)
\qbezier(5,2)(4,2)(3,2)
\qbezier(5,2)(5,2.5)(5,3)


\qbezier(1,2.5)(1.05,2.55)(1.1,2.6)
\qbezier(1,2.5)(.95,2.55)(.9,2.6)
\qbezier(3,2.5)(3.05,2.55)(3.1,2.6)
\qbezier(3,2.5)(2.95,2.55)(2.9,2.6)
\qbezier(2,3)(1.95,3.05)(1.9,3.1)
\qbezier(2,3)(1.95,2.95)(1.9,2.9)
\qbezier(2,2)(1.95,2.05)(1.9,2.1)
\qbezier(2,2)(1.95,1.95)(1.9,1.9)
\qbezier(4,3)(4.05,3.05)(4.1,3.1)
\qbezier(4,3)(4.05,2.95)(4.1,2.9)
\qbezier(3,1.5)(3.05,1.45)(3.1,1,4)
\qbezier(3,1.5)(2.95,1.45)(2.9,1,4)
\qbezier(3.5,1.5)(3.5,1.45)(3.5,1.4)
\qbezier(3.5,1.5)(3.55,1.5)(3.6,1.5)
\qbezier(4.5,1)(4.45,1.05)(4.4,1.1)
\qbezier(4.5,1)(4.45,.95)(4.4,.9)
\qbezier(4,2)(4.05,2.05)(4.1,2.1)
\qbezier(4,2)(4.05,1.95)(4.1,1.9)
\qbezier(5,2.5)(5.05,2.55)(5.1,2.6)
\qbezier(5,2.5)(4.95,2.55)(4.9,2.6)

\end{picture}

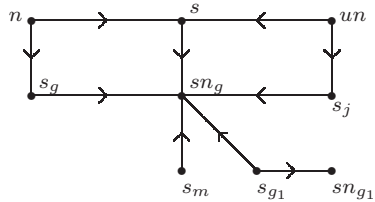
$Fig.~\ref{fig1}.9:~s~is~connected~to~a~stable~node~$
$and~to~a~second~unstable~node$
\end{center}

Observe that if the orientations of gradient lines are chosen as in
the case represented in figure \ref{fig1}.8, then the gradient line
$\gamma_{s_j,sn_g}$ has positive orientation. A computation proves
that 
$$Ker\partial_1^{U_2}=<s_m+s_g,s_j+s_g,s_j-s_m>=
<s_m+s_g,s_j+s_g>$$ 
and 
\begin{eqnarray}
Ker\partial_1^{U_1} & = &<s_m+s_g,s_j+s_g,s_j-s_m,
s-s_m,s+s_g,s-s_j>
\nonumber\\
& = & <s_m+s_g,s_j+s_g,s-s_g>
\nonumber
\end{eqnarray} 
thus
$$Ker\partial_1^{U_1}=i(Ker\partial_1^{U_2})\oplus<s-s_g>$$ 
which implies
$dimKer\partial^{U_1}=dimKer\partial^{U_2}+1$.

It remains to consider the degenerate cases. If $s_j=s_g$, 
a degenerate square appears in $U_2$, but the
relation just obtained above still holds. If $s_j=s$ and $s_g\neq s$,
then $s$ has a double connection only with $sn_g$ (and not with
$un$ because $s$ is already connected to $n$), which can be realized by
two different phase portraits, depending on how $\gamma_{s,sn_g}$
winds in the phase portrait: 
in one case $W^u(n)\cap W^s(sn_g)$ is not bounded and the phase portrait
can not be extended to the compactification $S^2$ of
$\R^2$ unless adding further critical points,
but this is not compatible with the compactness hypothesis \ref{assumption}; 
in the other case,
(note that the phase portrait, in order to be that of a gradient vector field, requires to 
be completed by the addition of further critical points)
the compactness hypothesis is fulfilled
and the expected relation between $Ker\partial_1^{U_1}$ and $Ker\partial_1^{U_2}$
is verified.
Finally, the case $s_g=s=s_j$ was already examined in 
the degenerate case of
$s$ connected to a unique unstable node $n$.

\item[***]
Suppose now that
$\partial_1^{U_2}(s)=sn_{g}+ sn_{k}$. Then, by lemma \ref{php}, there exist two saddles
$s_g$ and $s_k$ forming two distinct squares together with $n$, $s$
and, respectively, $sn_g$ and $sn_h$ (note that the squares can not be degenerate in $U_1$),
and suppose, for the moment, $s_g\neq s_k$. In figure \ref{fig1}.10 
it is presented the case where
$s$ is not connected to any other unstable node, except $n$.

\begin{center}
\setlength{\unitlength}{1cm}
\begin{picture}(6,6)
\label{fig10}
\thinlines

\put(1,2){\circle*{.1}}
\put(1,3){\circle*{.1}}
\put(1,4){\circle*{.1}}
\put(3,2){\circle*{.1}}
\put(3,4){\circle*{.1}}
\put(3,3){\circle*{.1}}
\put(3,5){\circle*{.1}}
\put(4,5){\circle*{.1}}
\put(5,5){\circle*{.1}}
\put(3,1){\circle*{.1}}
\put(4,1){\circle*{.1}}
\put(5,1){\circle*{.1}}

\scriptsize
\put(1.1,2.1){$\displaystyle s_g$}
\put(.7,3){$\displaystyle n$}
\put(1,4.1){$\displaystyle s_k$}
\put(3.2,3.1){$\displaystyle s$}
\put(3.1,2.1){$\displaystyle sn_g$}
\put(3.1,3.8){$\displaystyle sn_k$}
\put(3,5.1){$\displaystyle s_m$}
\put(4,5.1){$\displaystyle s_{k_1}$}
\put(5,5.1){$\displaystyle sn_{k_1}$}
\put(3,.7){$\displaystyle s_n$}
\put(4,.7){$\displaystyle s_{g_1}$}
\put(5,.7){$\displaystyle sn_{g_1}$}

\normalsize

\qbezier(1,2)(1,3)(1,4)
\qbezier(3,2)(3,3)(3,4)
\qbezier(1,2)(2,2)(3,2)
\qbezier(1,4)(2,4)(3,4)
\qbezier(1,3)(2,3)(3,3)
\qbezier(3,3)(4,3)(5,3)
\qbezier(3,5)(3,4.5)(3,4)
\qbezier(4,5)(3.5,4.5)(3,4)
\qbezier(4,5)(4.5,5)(5,5)
\qbezier(3,1)(3,1.5)(3,2)
\qbezier(4,1)(3.5,1.5)(3,2)
\qbezier(4,1)(4.5,1)(5,1)


\qbezier(1,2.5)(1.05,2.55)(1.1,2.6)
\qbezier(1,2.5)(.95,2.55)(.9,2.6)
\qbezier(1,3.5)(1.05,3.45)(1.1,3.4)
\qbezier(1,3.5)(.95,3.45)(.9,3.4)
\qbezier(3,2.5)(3.05,2.55)(3.1,2.6)
\qbezier(3,2.5)(2.95,2.55)(2.9,2.6)
\qbezier(3,3.5)(3.05,3.45)(3.1,3.4)
\qbezier(3,3.5)(2.95,3.45)(2.9,3.4)
\qbezier(2,4)(1.95,4.05)(1.9,4.1)
\qbezier(2,4)(1.95,3.95)(1.9,3.9)
\qbezier(2,3)(1.95,3.05)(1.9,3.1)
\qbezier(2,3)(1.95,2.95)(1.9,2.9)
\qbezier(2,2)(1.95,2.05)(1.9,2.1)
\qbezier(2,2)(1.95,1.95)(1.9,1.9)
\qbezier(4,3)(4.05,3.05)(4.1,3.1)
\qbezier(4,3)(4.05,2.95)(4.1,2.9)
\qbezier(3,4.5)(3.05,4.55)(3.1,4,6)
\qbezier(3,4.5)(2.95,4.55)(2.9,4,6)
\qbezier(3,1.5)(3.05,1.45)(3.1,1,4)
\qbezier(3,1.5)(2.95,1.45)(2.9,1,4)
\qbezier(3.5,4.5)(3.5,4.55)(3.5,4.6)
\qbezier(3.5,4.5)(3.55,4.5)(3.6,4.5)
\qbezier(4.5,5)(4.45,5.05)(4.4,5.1)
\qbezier(4.5,5)(4.45,4.95)(4.4,4.9)
\qbezier(3.5,1.5)(3.5,1.45)(3.5,1.4)
\qbezier(3.5,1.5)(3.55,1.5)(3.6,1.5)
\qbezier(4.5,1)(4.45,1.05)(4.4,1.1)
\qbezier(4.5,1)(4.45,.95)(4.4,.9)

\end{picture}
$Fig.~\ref{fig1}.10:~s~is~connected~to~two~stable~nodes$
\end{center}

If instead $s$
is connected to a second unstable node $un_1$, then, as figure \ref{fig1}.11 
shows, there exist two saddle
$s_{h}$ and $s_{j}$ forming a square with $un_1$, $s$ and respectively
$sn_k$ and $sn_g$. Suppose for the moment that all the saddles in the phase portrait
are distinct. 

\begin{center}
\setlength{\unitlength}{1cm}
\begin{picture}(6,6)
\label{fig11}
\thinlines

\put(1,2){\circle*{.1}}
\put(1,3){\circle*{.1}}
\put(1,4){\circle*{.1}}
\put(3,2){\circle*{.1}}
\put(3,4){\circle*{.1}}
\put(3,3){\circle*{.1}}
\put(3,5){\circle*{.1}}
\put(4,5){\circle*{.1}}
\put(5,5){\circle*{.1}}
\put(3,1){\circle*{.1}}
\put(4,1){\circle*{.1}}
\put(5,1){\circle*{.1}}
\put(5,4){\circle*{.1}}
\put(5,3){\circle*{.1}}
\put(5,2){\circle*{.1}}

\scriptsize
\put(1.1,2.1){$\displaystyle s_g$}
\put(.7,3){$\displaystyle n$}
\put(1,4.1){$\displaystyle s_k$}
\put(3.2,3.1){$\displaystyle s$}
\put(3.1,2.1){$\displaystyle sn_g$}
\put(3.1,3.8){$\displaystyle sn_k$}
\put(3,5.1){$\displaystyle s_m$}
\put(4,5.1){$\displaystyle s_{k_1}$}
\put(5,5.1){$\displaystyle sn_{k_1}$}
\put(3,.7){$\displaystyle s_n$}
\put(4,.7){$\displaystyle s_{g_1}$}
\put(5,.7){$\displaystyle sn_{g_1}$}
\put(5,4.1){$\displaystyle s_{h}$}
\put(5.1,3){$\displaystyle un_1$}
\put(5,1.8){$\displaystyle s_{j}$}

\normalsize

\qbezier(1,2)(1,3)(1,4)
\qbezier(3,2)(3,3)(3,4)
\qbezier(1,2)(2,2)(3,2)
\qbezier(1,4)(2,4)(3,4)
\qbezier(1,3)(2,3)(3,3)
\qbezier(3,3)(4,3)(5,3)
\qbezier(3,5)(3,4.5)(3,4)
\qbezier(4,5)(3.5,4.5)(3,4)
\qbezier(4,5)(4.5,5)(5,5)
\qbezier(3,1)(3,1.5)(3,2)
\qbezier(4,1)(3.5,1.5)(3,2)
\qbezier(4,1)(4.5,1)(5,1)
\qbezier(5,4)(4,4)(3,4)
\qbezier(5,4)(5,4.5)(5,3)
\qbezier(5,2)(4,2)(3,2)
\qbezier(5,2)(5,2.5)(5,3)


\qbezier(1,2.5)(1.05,2.55)(1.1,2.6)
\qbezier(1,2.5)(.95,2.55)(.9,2.6)
\qbezier(1,3.5)(1.05,3.45)(1.1,3.4)
\qbezier(1,3.5)(.95,3.45)(.9,3.4)
\qbezier(3,2.5)(3.05,2.55)(3.1,2.6)
\qbezier(3,2.5)(2.95,2.55)(2.9,2.6)
\qbezier(3,3.5)(3.05,3.45)(3.1,3.4)
\qbezier(3,3.5)(2.95,3.45)(2.9,3.4)
\qbezier(2,4)(1.95,4.05)(1.9,4.1)
\qbezier(2,4)(1.95,3.95)(1.9,3.9)
\qbezier(2,3)(1.95,3.05)(1.9,3.1)
\qbezier(2,3)(1.95,2.95)(1.9,2.9)
\qbezier(2,2)(1.95,2.05)(1.9,2.1)
\qbezier(2,2)(1.95,1.95)(1.9,1.9)
\qbezier(4,3)(4.05,3.05)(4.1,3.1)
\qbezier(4,3)(4.05,2.95)(4.1,2.9)
\qbezier(3,4.5)(3.05,4.55)(3.1,4,6)
\qbezier(3,4.5)(2.95,4.55)(2.9,4,6)
\qbezier(3,1.5)(3.05,1.45)(3.1,1,4)
\qbezier(3,1.5)(2.95,1.45)(2.9,1,4)
\qbezier(3.5,4.5)(3.5,4.55)(3.5,4.6)
\qbezier(3.5,4.5)(3.55,4.5)(3.6,4.5)
\qbezier(4.5,5)(4.45,5.05)(4.4,5.1)
\qbezier(4.5,5)(4.45,4.95)(4.4,4.9)
\qbezier(3.5,1.5)(3.5,1.45)(3.5,1.4)
\qbezier(3.5,1.5)(3.55,1.5)(3.6,1.5)
\qbezier(4.5,1)(4.45,1.05)(4.4,1.1)
\qbezier(4.5,1)(4.45,.95)(4.4,.9)
\qbezier(4,4)(4.05,4.05)(4.1,4.1)
\qbezier(4,4)(4.05,3.95)(4.1,3.9)
\qbezier(5,3.5)(5.05,3.45)(5.1,3.4)
\qbezier(5,3.5)(4.95,3.45)(4.9,3.4)
\qbezier(4,2)(4.05,2.05)(4.1,2.1)
\qbezier(4,2)(4.05,1.95)(4.1,1.9)
\qbezier(5,2.5)(5.05,2.55)(5.1,2.6)
\qbezier(5,2.5)(4.95,2.55)(4.9,2.6)

\end{picture}

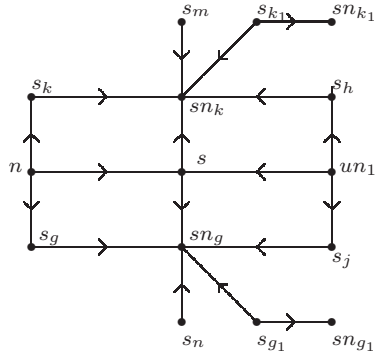
$Fig.~\ref{fig1}.11:~s~is~connected~to~two~stable~nodes~$
$and~to~a~second~unstable~node$
\end{center}

In all the cases, it can be checked that a choice of a coherent orientation of squares
allows to prove that $dimKer\partial_1^{U_1}=dimKer\partial_1^{U_2}+1$. For instance,
regarding the situation pictured in figure \ref{fig1}.9, choosing a coherent
orientation, for example one such that $\gamma_{s,sn_g}$, $\gamma_{s_n,sn_g}$,
$\gamma_{s_j,sn_g}$ and $\gamma_{s_k,sn_k}$ have positive orientation while
$\gamma_{s_g,sn_g}$, $\gamma_{s,sn_k}$, $\gamma_{s_h,sn_k}$ and $\gamma_{s_m,sn_k}$
have negative orientation, generators can be chosen such that
$$Ker\partial_1^{U_2}=<s_g+s_j,s_g+s_n,s_k+s_h,s_k+s_m>$$ 
and
$$Ker\partial_1^{U_1}=i(Ker\partial_1^{U_2})\oplus<s+s_k+s_g>$$     
This result is also achieved when, as already proved, sequences
of squares starting, as in figure \ref{fig1}.5, 
from $n$ and $s_g$ or from $n$ and $s_h$ are added in the phase portrait.

Finally, note that
$s_g$ may coincide with $s_j$ or $s_k$, and with $s_h$ provided, in this case, both $un_1$
and $sn_k$ are connected to $s_g$, that is, $s_g=s_h=s_j=s_k$. Depending
on which direction the gradient lines between saddles and nodes wind
in the phase portrait, further critical points must be included in order
the vector field to be gradient. Moreover, when $un_1$ is in the phase portrait, it may
happen that some moduli space of gradient lines
from unstable nodes to stable nodes is not bounded and the phase portrait can not
be extended to the comapctification $S^2$ of $\R^2$ unless adding further
critical points:
for example, this happens with $W^u(un_1)\cap W^s(sn_k)$ when $s_k=s_g=s_g$
(in this case, the extension of the phase portrait to $S^2$ 
yields a homoclinic orbit through $\infty$).
However, when the compactness hypothesis \ref{assumption} is fulfilled,
$dimKer\partial^{U_1}=dimKer\partial^{U_2}+1$.


\item
Suppose now the birth-death pair $(p_{k+1},p_{k+2})$ 
is represented by a stable node $n$ and a saddle $s$:
we are going to prove that $dimIm\partial_0^{U_1}=dimIm\partial_0^{U_2}$ and
$dimKer\partial_1^{U_1}=dimKer\partial_1^{U_2}$.

The existence in $U_1$, by lemma \ref{exgrl}, of 
a gradient line $\gamma_{n,s}$ from $n$ to $s$,
implies, as already explained for a birth-death pair given
by an unstable node and a saddle,
that if, in $U_1$, $n$ is connected to other
unstable nodes $un_j$ and saddles $s_i$ and if $s$ is connected at most to
two unstable nodes and to a second
stable node $sn$, then in $U_2$ gradient lines from $un_j$ and from $s_i$
to $sn$ appear in the phase portrait
(see figure \ref{fig1}.1 after reverting all the arrows).
Such gradient lines from $s_i$ to $sn$ implies a change in $Ker\partial_1^{U_2}$
with respect to $Ker\partial_1^{U_1}$.

The saddle $s$
can be connected to at most two stable nodes, one of which is $n$, and
two unstable nodes. 

\item[*]
Consider first the simplest case: $s$ is connected only to
$n$: 

\item[**]
clearly $dimIm\partial_0^{U_1}=dimIm\partial_0^{U_2}$.

\item[**]
as to kernels, 
in priciple
also other saddles $s_1$,...,$s_i$ might be connected to $n$, of these,
say $s_1$, ...,$s_j$, with $j\leq i$, might be each connected to a second
stable node $sn_1$, ..., $sn_j$, as in figure \ref{fig1}.12 (more generally, finite
chains of saddles and stable nodes, as $(n,s_{k_1}$, $sn_{k_1}$, $s_{k_2}$, $sn_{k_2}$,...),
starting from $n$ and ending with a saddle or a stable node, should be considered,
however the lenght of such chains does not affect the argument).

\begin{center}
\setlength{\unitlength}{1cm}
\begin{picture}(6,4)
\label{fig12}
\thinlines

\put(1,2){\circle*{.1}}
\put(3,3){\circle*{.1}}
\put(2,2){\circle*{.1}}
\put(3,1){\circle*{.1}}
\put(4,1){\circle*{.1}}

\scriptsize
\put(1,2.1){$\displaystyle s$}
\put(3.1,2.8){$\displaystyle s_2$}
\put(3,.7){$\displaystyle s_1$}
\put(4,.7){$\displaystyle sn_1$}
\put(1.9,2.1){$\displaystyle n$}

\normalsize

\qbezier(1,2)(1.5,2)(2,2)
\qbezier(2,2)(2.5,2.5)(3,3)
\qbezier(2,2)(2.5,1.5)(3,1)
\qbezier(3,1)(3.5,1)(4,1)


\qbezier(1.5,2)(1.45,2.05)(1.4,2.1)
\qbezier(1.5,2)(1.45,1.95)(1.4,1.9)

\qbezier(2.5,2.5)(2.5,2.55)(2.5,2.6)
\qbezier(2.5,2.5)(2.55,2.5)(2.6,2.5)

\qbezier(2.5,1.5)(2.5,1.45)(2.5,1.4)
\qbezier(2.5,1.5)(2.55,1.5)(2.6,1.5)

\qbezier(3.5,1)(3.45,1.05)(3.4,1.1)
\qbezier(3.5,1)(3.45,.95)(3.4,.9)

\end{picture}

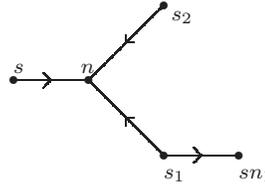
$Fig.~\ref{fig1}.12:~a~saddle~and~a~stable~node~forming~a~birth-death~pair$
\end{center}

Choosing the orientation for which
the gradient lines $\gamma_{s,n}$, $\gamma_{s_l,n}$, for $1\leq l\leq i$, are positive
and, as a consequence, $\gamma_{s_r,sn_r}$, for $1\leq r\leq j$, are negative,
it follows that 
$$Ker\partial_1^{U_2}=<s_{j+1},...,s_i>$$ 
and 
$$Ker\partial_1^{U_1}=<s-s_{j+1},...,s-s_i>$$ 
and so $dimKer\partial_1^{U_1}=dimKer\partial_1^{U_2}$. Note that the saddles
$s_1$, ...,$s_j$ were irrelevant in the computation of the kernel, so from now on, assume
every chain ends with a saddle.

\item[*]
Suppose also now that $s$ is connected to a 
second stable node $sn$, but no unstable nodes are connected to $s$. 

\item[**]
Clearly $dimIm\partial_0^{U_1}=dimIm\partial_0^{U_2}$.

\item[**]
As explained, in $U_2$ gradient
lines from the saddles $s_1$,...,$s_i$ to $sn$ appear in the phase portrait. 
Suppose that other saddles
$s^{'}_1$,...,$s^{'}_r$ are connected to $sn$ (see figure \ref{fig1}.13).

\begin{center}
\setlength{\unitlength}{1cm}
\begin{picture}(7,4)
\label{fig13}
\thinlines

\put(1,3){\circle*{.1}}
\put(1,1){\circle*{.1}}
\put(3,2){\circle*{.1}}
\put(2,2){\circle*{.1}}
\put(4,2){\circle*{.1}}
\put(5,1){\circle*{.1}}
\put(5,3){\circle*{.1}}
\put(5,2){\circle*{.1}}

\scriptsize
\put(3,2.1){$\displaystyle s$}
\put(5,.7){$\displaystyle s_1$}
\put(5,1.7){$\displaystyle s_2$}
\put(5,2.7){$\displaystyle s_3$}
\put(3.9,2.1){$\displaystyle n$}
\put(2,2.1){$\displaystyle sn$}
\put(.7,1.1){$\displaystyle s^{'}_1$}
\put(.7,2.9){$\displaystyle s^{'}_2$}

\normalsize

\qbezier(3,2)(2.5,2)(2,2)
\qbezier(4,2)(3.5,2)(3,2)
\qbezier(4,2)(4.5,2)(5,2)
\qbezier(4,2)(4.5,1.5)(5,1)
\qbezier(4,2)(4.5,2.5)(5,3)
\qbezier(2,2)(1.5,2.5)(1,3)
\qbezier(2,2)(1.5,1.5)(1,1)


\qbezier(2.5,2)(2.55,2.05)(2.6,2.1)
\qbezier(2.5,2)(2.55,1.95)(2.6,1.9)

\qbezier(3.5,2)(3.45,2.05)(3.4,2.1)
\qbezier(3.5,2)(3.45,1.95)(3.4,1.9)

\qbezier(4.5,2)(4.55,2.05)(4.6,2.1)
\qbezier(4.5,2)(4.55,1.95)(4.6,1.9)

\qbezier(4.5,2.5)(4.5,2.55)(4.5,2.6)
\qbezier(4.5,2.5)(4.55,2.5)(4.6,2.5)

\qbezier(4.5,1.5)(4.5,1.45)(4.5,1.4)
\qbezier(4.5,1.5)(4.55,1.5)(4.6,1.5)

\qbezier(1.5,2.5)(1.5,2.55)(1.5,2.6)
\qbezier(1.5,2.5)(1.45,2.5)(1.4,2.5)

\qbezier(1.5,1.5)(1.5,1.45)(1.5,1.4)
\qbezier(1.5,1.5)(1.45,1.5)(1.4,1.5)

\end{picture}

$Fig.~\ref{fig1}.13:~s~is~connected~to~a~second~stable~node$
\end{center}

Choose an orientation such that the gradient lines
$\gamma_{s^{'}_l,sn}$, for $1\leq l\leq r$, are positive. 
Observe that $\partial_1^{U_2}(s)=n-sn$. A computation
shows that 
$$Ker\partial_1^{U_1}=<s-s_1+s^{'}_1,...,s-s_1+s^{'}_r,s-s_2+s^{'}_1,...,s-s_i-s^{'}_1>$$
and
$$Ker\partial_1^{U_2}=<-s_1+s^{'}_1,...,-s_1+s^{'}_r,-s_2+s^{'}_1,...,-s_i-s^{'}_1>$$
and so $dimKer\partial_1^{U_1}=dimKer\partial_1^{U_2}$. 

Note that if $s_l=s^{'}_m$ for
some $1\leq l\leq i$ and $1\leq m\leq r$, then, as shown in figure \ref{fig1}.14,
there exists at least an unstable node to which $s$ is connected,
anticipating the case, where $s$ is connected to an unstable node, which will be considered soon.

\begin{center}
\setlength{\unitlength}{1cm}
\begin{picture}(4,4)
\label{fig14}
\thinlines

\put(3.5,2){\circle*{.1}}
\put(2,.5){\circle*{.1}}
\put(.5,2){\circle*{.1}}
\put(2,3.5){\circle*{.1}}
\put(2,2){\circle*{.1}}

\scriptsize
\put(2,3.6){$\displaystyle sn$}
\put(0.3,1.8){$\displaystyle s$}
\put(2,.2){$\displaystyle n$}
\put(3.6,2){$\displaystyle s_l=s^{'}_m$}
\put(2.1,2.1){$\displaystyle un$}
\normalsize

\qbezier(3.5,2)(2.75,1.25)(2,.5)
\qbezier(3.5,2)(2.75,2.75)(2,3.5)
\qbezier(.5,2)(1.25,1.25)(2,.5)
\qbezier(.5,2)(1.25,2.75)(2,3.5)

\qbezier(.5,2)(1,2)(2,2)
\qbezier(3.5,2)(3,2)(2,2)
\qbezier(2,.5)(2,1)(2,2)
\qbezier(2,3.5)(2,3)(2,2)

\qbezier(2.75,2.75)(2.76,2.75)(2.85,2.75)
\qbezier(2.75,2.75)(2.75,2.74)(2.75,2.65)

\qbezier(1.25,2.75)(1.2,2.75)(1.15,2.75)
\qbezier(1.25,2.75)(1.25,2.74)(1.25,2.65)

\qbezier(2.75,1.25)(2.76,1.25)(2.85,1.25)
\qbezier(2.75,1.25)(2.75,1.3)(2.75,1.35)

\qbezier(1.25,1.25)(1.2,1.25)(1.15,1.25)
\qbezier(1.25,1.25)(1.25,1.3)(1.25,1.35)

\qbezier(2.75,2)(2.7,2.05)(2.65,2.1)
\qbezier(2.75,2)(2.7,1.95)(2.65,1.9)

\qbezier(1.25,2)(1.3,2.05)(1.35,2.1)
\qbezier(1.25,2)(1.3,1.95)(1.35,1.9)

\qbezier(2,1.25)(2.05,1.3)(2.1,1.35)
\qbezier(2,1.25)(1.95,1.3)(1.9,1.35)

\qbezier(2,2.75)(2.05,2.7)(2.1,2.65)
\qbezier(2,2.75)(1.95,2.7)(1.9,2.65)


\end{picture}

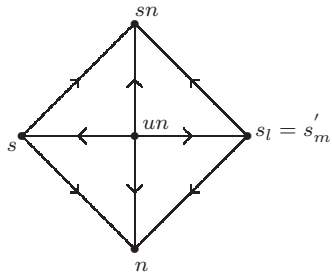
$Fig.~\ref{fig1}.14:~if~s_l=s^{'}_m~there~is~also~an~unstable~node$
\end{center}

\item[*]
Allow now $s$ to be connected also to one unstable node $un$. 
By lemma \ref{php} there exist saddles $s_1$ and $s^{'}_1$ forming
squares with $un$, $sn$ and respectively $s_1$ and $s^{'}_1$. Suppose, for the
moment, the two squares non-degenerate and $s_1\neq s^{'}_1$.
Suppose $un$ is connected, besides $s$, $s_1$ and $s^{'}_1$, to other saddles
$s_{g_1}$, ..., $s_{g_m}$, so that $\partial_0^{U_1}(un)=s+s_1+s^{'}_1+\sum^m_{i=1} s_{g_i}$
(as usual, signs depend on the orientation, however
fixing a choice does not modify the argument).
Suppose that, both in $U_1$ and in $U_2$, the unstable nodes are $un$ and 
$un_{h_1}$, ..., $un_{h_n}$.
Since the only unstable node connected to $s$ is $un$, 
$$Im\partial_0^{U_1}=i(<\sum^n_{j=1}\partial_0^{U_2}(un_{h_j})>)\oplus\partial_0^{U_1}(un)$$
moreover, being $\partial_0^{U_2}(un)=s_1+s^{'}_1+\sum^m_{i=1} s_{g_i}$ and if\\
$s_1+s^{'}_1+\sum^m_{i=1} s_{g_i}\notin<\sum^n_{j=1}\partial_0^{U_2}(un_{h_j})>$ then
$$Im\partial_0^{U_2}=<\sum^n_{j=1}\partial_0^{U_2}(un_{h_j})>\oplus\partial_0^{U_2}(un)$$
and so $dimIm\partial_0^{U_1}=dimIm\partial_0^{U_2}$.
If instead $s_1+s^{'}_1+\sum^m_{i=1} s_{g_i}\in\\<\sum^n_{j=1}\partial_0^{U_2}(un_{h_j})>$ 
then $\partial_0^{U_2}(un)$ is a linear combination of $\partial_0^{U_2}(un_{h_j})$: 
consider the situation represented in figure \ref{fig1}.15, where both $un$ and 
a second unstable node $un_1$ are connected to the same pair of saddles in such a way that
$\partial_0^{U_i}(un)=\partial_0^{U_i}(un)$, for $i=1,2$,

\begin{center}
\setlength{\unitlength}{1cm}
\begin{picture}(6,7)
\label{fig15}
\thinlines

\put(1,4){\circle*{.1}}
\put(1,5){\circle*{.1}}
\put(3,4){\circle*{.1}}
\put(5,4){\circle*{.1}}
\put(3,5){\circle*{.1}}
\put(5,5){\circle*{.1}}
\put(3,3){\circle*{.1}}
\put(3,2){\circle*{.1}}

\scriptsize
\put(1.1,4.1){$\displaystyle s_1$}
\put(.7,5.1){$\displaystyle sn$}
\put(3.1,5.1){$\displaystyle s$}
\put(3.1,4.1){$\displaystyle un$}
\put(4.9,5.1){$\displaystyle n$}
\put(5,3.8){$\displaystyle s_1$}
\put(3.1,2.8){$\displaystyle sn_1$}
\put(3.1,1.8){$\displaystyle un_1$}
\put(2.9,.7){$\displaystyle \infty$}
\put(2.9,6.1){$\displaystyle \infty$}

\normalsize

\qbezier(1,4)(1,4.5)(1,5)
\qbezier(3,4)(3,4.5)(3,5)
\qbezier(1,4)(2,4)(3,4)
\qbezier(1,5)(2,5)(3,5)
\qbezier(3,5)(4,5)(5,5)
\qbezier(5,4)(4,4)(3,4)
\qbezier(5,4)(5,4.5)(5,5)
\qbezier(1,4)(2,3.5)(3,3)
\qbezier(1,4)(2,3)(3,2)
\qbezier(5,4)(4,3.5)(3,3)
\qbezier(5,4)(4,3)(3,2)
\qbezier(3,1)(3,1.5)(3,2)
\qbezier(3,6)(3,5.5)(3,5)


\qbezier(1,4.5)(1.05,4.45)(1.1,4.4)
\qbezier(1,4.5)(.95,4.45)(.9,4.4)

\qbezier(3,4.5)(3.05,4.45)(3.1,4.4)
\qbezier(3,4.5)(2.95,4.45)(2.9,4.4)

\qbezier(2,5)(2.05,5.05)(2.1,5.1)
\qbezier(2,5)(2.05,4.95)(2.1,4.9)
\qbezier(2,4)(2.05,4.05)(2.1,4.1)
\qbezier(2,4)(2.05,3.95)(2.1,3.9)
\qbezier(4,5)(3.95,5.05)(3.9,5.1)
\qbezier(4,5)(3.95,4.95)(3.9,4.9)

\qbezier(4,4)(3.95,4.05)(3.9,4.1)
\qbezier(4,4)(3.95,3.95)(3.9,3.9)
\qbezier(5,4.5)(5.05,4.45)(5.1,4.4)
\qbezier(5,4.5)(4.95,4.45)(4.9,4.4)

\qbezier(3,5.5)(3.05,5.55)(3.1,5.6)
\qbezier(3,5.5)(2.95,5.55)(2.9,5.6)

\qbezier(3,1.5)(3.05,1.55)(3.1,1.6)
\qbezier(3,1.5)(2.95,1.55)(2.9,1.6)

\qbezier(4,3)(3.95,3)(3.9,3)
\qbezier(4,3)(4,2.95)(4,2.9)

\qbezier(2,3)(2.05,3)(2.1,3)
\qbezier(2,3)(2,2.95)(2,2.9)

\qbezier(4,3.5)(4.025,3.55)(4.05,3.6)
\qbezier(4,3.5)(4.05,3.475)(4.1,3.45)

\qbezier(2,3.5)(1.975,3.55)(1.95,3.6)
\qbezier(2,3.5)(1.95,3.475)(1.9,3.45)

\end{picture}

$Fig.~\ref{fig1}.15:~un~and~ 
a~second~unstable~node~un_1$
$are~connected~to~the~same~pair~of~saddles$
\end{center}

Observe
that $s$ has a separatrix from $\infty$, thus, if the compactness hypothesis \ref{assumption}
is fulfilled, there exists an unstable node $un_{h_p}$ with a gradient line from 
$un_{h_p}$ to $\infty$; this gradient line, in the compactification $S^2$, becomes
a gradient line from $un_{h_p}$ to $s$ and $(un_{h_p},s,s_1,n)$ and
$(un_{h_p},s,s^{'}_1,sn)$ are squares; thus $Im\partial_0^{U_1}(un)\in
i(<\sum^n_{j=1}\partial^{U_2}(un_{h_j})>)$, and so it follows again
$dimIm\partial_0^{U_1}=dimIm\partial_0^{U_2}$. Now, since an unstable node connected to $s$
does not affect the kernel of $\partial_1$, it also follows 
$dimKer\partial_1^{U_1}=dimKer\partial_1^{U_2}$.

It remains to consider few special cases. 

\item[**]
The situation where $s_1= s^{'}_1$
does not differ so much from that above: the same arguments can be suitably applied
and the same conclusion achieved. 

\item[**]
If, instead one or both the 
squares formed by $un$ with $s$ and 
respectively $s_1$ and $s^{'}_1$ are degenerate, the only possibility is
that $W^u(s)\subset W^u(un)$, that is, $s$ has two separatrices connecting
it to $un$ and having opposite orientations: this implies 
$$\partial_0^{U_1}(n)=i(\partial_0^{U_2}(n))$$
and so again $dimIm\partial_1^{U_1}=dimIm\partial_1^{U_2}$.

\item[*]
Suppose now that $s$ is connected to two unstable nodes $un_1$ and $un_2$.
Then, by lemma \ref{php}, there are saddles $s_i$ and $s^{'}_i$ forming squares
(which can not be degenerate) with $un_i$, $s$ and respectively $n$ and $sn$,
for $i=1,2$. The proof that $dimIm\partial_0^{U_1}=dimIm\partial_0^{U_2}$
and $dimKer\partial_1^{U_1}=dimKer\partial_1^{U_2}$ goes along the same lines
of the case where $s$ is connected to a unique unstable node, however
some attention must be paid to prove the relation between
$Ker\partial_1^{U_1}$ and $Ker\partial_1^{U_2}$ when 
$\partial_0^{U_1}(un_i)=\partial_0^{U_1}(\tilde{un}_i)$ for some unstable node $\tilde{un}_i$ 
(see figure \ref{fig1}.16): 
in this case, in $U_1$ there are
moduli spaces of gradient lines from unstable to stable nodes which
are not bounded (namely, ${\cal M}(\tilde{un}_i,n)$ and ${\cal M}(\tilde{un}_i,sn)$) and the phase portrait cannot be extended to the compactification
$S^2$ of $\R^2$ unless adding
further critical points. 
This case is not consistent with the compactness hypothesis
\ref{assumption} (however see also remark \ref{infinito}).

\begin{center}
\setlength{\unitlength}{1cm}
\begin{picture}(6,10)
\label{fig16}
\thinlines

\put(1,4){\circle*{.1}}
\put(1,5){\circle*{.1}}
\put(3,4){\circle*{.1}}
\put(5,4){\circle*{.1}}
\put(3,5){\circle*{.1}}
\put(5,5){\circle*{.1}}
\put(3,3){\circle*{.1}}
\put(3,2){\circle*{.1}}

\put(1,6){\circle*{.1}}
\put(3,6){\circle*{.1}}
\put(5,6){\circle*{.1}}
\put(3,7){\circle*{.1}}
\put(3,8){\circle*{.1}}

\scriptsize
\put(1.1,4.1){$\displaystyle s_1^{'}$}
\put(1.1,5.1){$\displaystyle sn$}
\put(3.1,5.1){$\displaystyle s$}
\put(3.1,4.1){$\displaystyle un_1$}
\put(4.7,5.1){$\displaystyle n$}
\put(5,3.8){$\displaystyle s_1$}
\put(3.1,2.8){$\displaystyle sn_1$}
\put(3.1,1.8){$\displaystyle \tilde{un}_1$}
\put(2.9,.7){$\displaystyle \infty$}
\put(2.9,9.1){$\displaystyle \infty$}

\put(.6,6.1){$\displaystyle s_2^{'}$}
\put(3.1,6.1){$\displaystyle un_2$}
\put(5.1,6){$\displaystyle s_2$}
\put(2.8,7.1){$\displaystyle sn_2$}
\put(3.1,8.1){$\displaystyle \tilde{un}_2$}
\put(0,5.15){$\displaystyle \infty$}
\put(5.7,5.15){$\displaystyle \infty$}

\normalsize

\qbezier(1,4)(1,4.5)(1,5)
\qbezier(3,4)(3,4.5)(3,5)
\qbezier(1,4)(2,4)(3,4)
\qbezier(1,5)(2,5)(3,5)
\qbezier(3,5)(4,5)(5,5)
\qbezier(5,4)(4,4)(3,4)
\qbezier(5,4)(5,4.5)(5,5)
\qbezier(1,4)(2,3.5)(3,3)
\qbezier(1,4)(2,3)(3,2)
\qbezier(5,4)(4,3.5)(3,3)
\qbezier(5,4)(4,3)(3,2)
\qbezier(3,1)(3,1.5)(3,2)

\qbezier(1,6)(1,5.5)(1,5)
\qbezier(3,6)(3,5.5)(3,5)
\qbezier(1,6)(2,6)(3,6)
\qbezier(5,6)(4,6)(3,6)
\qbezier(5,6)(5,5.5)(5,5)
\qbezier(1,6)(2,6.5)(3,7)
\qbezier(1,6)(2,7)(3,8)
\qbezier(5,6)(4,6.5)(3,7)
\qbezier(5,6)(4,7)(3,8)
\qbezier(3,9)(3,8.5)(3,8)
\qbezier(1,5)(.5,5)(0,5)
\qbezier(6,5)(4.5,5)(5,5)


\qbezier(1,4.5)(1.05,4.45)(1.1,4.4)
\qbezier(1,4.5)(.95,4.45)(.9,4.4)

\qbezier(3,4.5)(3.05,4.45)(3.1,4.4)
\qbezier(3,4.5)(2.95,4.45)(2.9,4.4)

\qbezier(2,5)(2.05,5.05)(2.1,5.1)
\qbezier(2,5)(2.05,4.95)(2.1,4.9)
\qbezier(2,4)(2.05,4.05)(2.1,4.1)
\qbezier(2,4)(2.05,3.95)(2.1,3.9)
\qbezier(4,5)(3.95,5.05)(3.9,5.1)
\qbezier(4,5)(3.95,4.95)(3.9,4.9)

\qbezier(2,6)(2.05,6.05)(2.1,6.1)
\qbezier(2,6)(2.05,5.95)(2.1,5.9)
\qbezier(4,6)(3.95,6.05)(3.9,6.1)
\qbezier(4,6)(3.95,5.95)(3.9,5.9)

\qbezier(4,4)(3.95,4.05)(3.9,4.1)
\qbezier(4,4)(3.95,3.95)(3.9,3.9)
\qbezier(5,4.5)(5.05,4.45)(5.1,4.4)
\qbezier(5,4.5)(4.95,4.45)(4.9,4.4)

\qbezier(3,5.5)(3.05,5.55)(3.1,5.6)
\qbezier(3,5.5)(2.95,5.55)(2.9,5.6)

\qbezier(3,1.5)(3.05,1.55)(3.1,1.6)
\qbezier(3,1.5)(2.95,1.55)(2.9,1.6)

\qbezier(4,3)(3.95,3)(3.9,3)
\qbezier(4,3)(4,2.95)(4,2.9)

\qbezier(2,3)(2.05,3)(2.1,3)
\qbezier(2,3)(2,2.95)(2,2.9)

\qbezier(4,3.5)(4.025,3.55)(4.05,3.6)
\qbezier(4,3.5)(4.05,3.475)(4.1,3.45)

\qbezier(2,3.5)(1.975,3.55)(1.95,3.6)
\qbezier(2,3.5)(1.95,3.475)(1.9,3.45)

\qbezier(4,7)(3.95,7)(3.9,7)
\qbezier(4,7)(4,7.05)(4,7.1)

\qbezier(2,7)(2.05,7)(2.1,7)
\qbezier(2,7)(2,7.05)(2,7.1)

\qbezier(2,6.5)(1.975,6.45)(1.95,6.4)
\qbezier(2,6.5)(1.95,6.525)(1.9,6.55)

\qbezier(4,6.5)(4.025,6.45)(4.05,6.4)
\qbezier(4,6.5)(4.05,6.525)(4.1,6.55)

\qbezier(1,5.5)(1.05,5.55)(1.1,5.6)
\qbezier(1,5.5)(.95,5.55)(.9,5.6)

\qbezier(5,5.5)(5.05,5.55)(5.1,5.6)
\qbezier(5,5.5)(4.95,5.55)(4.9,5.6)

\qbezier(3,8.5)(3.05,8.45)(3.1,8.4)
\qbezier(3,8.5)(2.95,8.45)(2.9,8.4)

\qbezier(.5,5)(.45,5.05)(.4,5.1)
\qbezier(.5,5)(.45,4.95)(.4,4.9)

\qbezier(5.5,5)(5.55,5.05)(5.6,5.1)
\qbezier(5.5,5)(5.55,4.95)(5.6,4.9)

\end{picture}

$Fig.~\ref{fig1}.16:\partial_0^{U_1}(un_i)=\partial_0^{U_1}(\tilde{un}_i)$
\end{center}

Finally, if some of the saddles forming the square coincide,
the situation can be analyzed
in a way similar 
to the case, already examined, where $n$ is an unstable node.
\end{itemize}

\item
It remains now to confront $HM_0(U_1)$ and $HM_2(U_2)$, which amounts to compute
respectively $Ker\partial_0$ and $Im\partial_1$.

\begin{itemize}

\item
Suppose first that the birth-death pair $(n,s)$ is given by an
unstable node $n$ and a saddle $s$. 

\item[*]
Consider $HM_0$. In $U_1$, let $un$ be an eventual second unstable node connected
to $s$. Consider chains $(n,s_{1i}^{'},un_{1i}^{'},...)$ and $(un,s_{1j},un_{1j},$ ...)
of saddles and unstable nodes starting respectively from $n$ and $un$. If at least one chain of each type 
ends with an unstable nodes, then they determine an element in $Ker\partial_0^{U_1}$ which
is not in $Ker\partial_0^{U_2}$. For simplicty, suppose
to have only one chain of each type and of shortest lenght: $(n,s_{1}^{'},un_{1}^{'})$ and $(n,s_{1},un_{1})$.
Then $n+un+un_1^{'}+un_1\in Ker\partial_0^{U_1}$. In $U_2$, $un$ becomes connected to $s_1^{'}$, giving
$un+un_1^{'}+un_1\in Ker\partial_0^{U_2}$. Since $un+un_1^{'}+un_1\notin Ker\partial_0^{U_1}$, it follows
that $dimKer\partial_0^{U_1}=dimKer\partial_0^{U_2}$ and so $HM_0(U_1)\cong HM_0(U_2)$.

\item[*]
Consider now $HM_2$. Since $Im\partial_1^{U_1}=i(Im\partial_1^{U_2})\oplus <\partial_1^{U_1}(s)>$ it remains 
to prove that $\partial_1^{U_1}(s)\in i(Im\partial_1^{U_2})$. Suppose $\partial_1^{U_1}(s)=sn_1+sn_2$,
then by lemma \ref{php} there are saddles $s_i$, with $i=1,2$, forming squares together with $n$, $s$ and $sn_i$. 
Suppose also $\partial_1^{U_1}(s_i)=sn_i+sn_i^{'}$ and note that $i(\partial_1^{U_2}(s_i))=sn_i+sn_i^{'}$ and that
by lemma \ref{php} there is a saddle $s^{'}$ connected to $n$ and such that $\partial_1^{U_1}(s^{'})=sn_1^{'}+sn_2^{'}=
i(\partial_1^{U_2}(s^{'}))$. It follows now that $\partial_1^{U_1}(s)=i(\partial_1^{U_1}(s_1-s_2+s^{'}))$. This proves
$HM_2(U_1)\cong HM_2(U_2)$.

\item
Suppose now that the birth-death pair $(n,s)$ is given by a
stable node $n$ and a saddle $s$.

\item[*]
Consider $HM_0$. Let $un_i$, for $i=1,2$, two unstable nodes eventually connected to $s$. By lemma
\ref{php} there are saddles $s_i$ forming squares respectively with $un_i$, $s$ and $n$. This implies
that any linear combination of $un_1$ and $un_2$ never belongs to $Ker\partial_0^{U_1}$, thus
$dimKer\partial_0^{U_1}=dimKer\partial_0^{U_2}$ and so $HM_0(U_1)\cong HM_0(U_2)$.

\item[*]
Consider now $HM_2$. As $n\in Ker\partial_2^{U_1}$, it follows that 
$dimKer\partial_2^{U_1}=dimKer\partial_2^{U_2}+1$. It remains to prove  
$dimIm\partial_1^{U_1}=dimIm\partial_1^{U_2}+1$. Suppose $sn$ is an eventual
second stable node connected to $s$ and consider chains
$(n,s_{1i}^{'},sn_{1i}^{'},...)$ and $(sn,s_{1j},sn_{1j},...)$ of saddles and
stable nodes starting respectively from $n$ and $sn$. Observe that 
$i(\partial_1^{U_2}(s^{'}_{ki}))$, $i(\partial_1^{U_2}(s_{kj}))\in Im\partial_1^{U_1}$,
and that $Im\partial_1^{U_1}$ contains also $\partial_1^{U_1}(s)=n-sn$ (for some
choice of signs). Hence $dimIm\partial_1^{U_1}=dimIm\partial_1^{U_2}+1$ and $HM_2(U_1)\cong HM_2(U_2)$.
\end{itemize}
\end{enumerate}
\end{proof}

\begin{remark}
\label{infinito}
\rm
Observe that, in the last case shown, for instance, in figure \ref{fig16}.16,
$\infty$ is the critical point to be added in order to extend the phase portrait
to $S^2$, particularly, is a saddle:
in this way on $S^2$ the relation $dimIm\partial_0^{U_1}=dimIm\partial_0^{U_2}$
still holds.
\end{remark}

\begin{remark}
\label{remarkcau}
\rm
Let $i$ denote the natural injection
$$\C[p_{1}^{U_2}]\oplus ... \oplus\C[p_{k}^{U_2}]\hookrightarrow
\C[p_{1}^{U_1}]\oplus ... \oplus\C[p_k^{U_1}]
\oplus\C[p_{k+1}^{U_1}]\oplus\C[p_{k+2}^{U_1}]$$
and observe that:
\begin{enumerate}
\item if the birth-death pair is given by an unstable node $n$ and a saddle $s$
then 
$$dimKer\partial_1^{U_1}=dimKer\partial_1^{U_2}+1$$  
$$dimIm\partial_0^{U_1}=dimIm\partial_0^{U_2}+1$$
and
$$i(Ker\partial_1^{U_2})\subsetneqq Ker\partial_1^{U_1}$$ 
$$Im\partial_0^{U_1}=i(Im\partial_0^{U_2})\oplus<\partial_0^{U_1}(n)>$$ 
moreover, the 1-dimensional complement which added to $i(Ker\partial_1^{U_2})$ yields
$Ker\partial_1^{U_1}$ contains the saddle $s$: a choice of a generator may be for example,
referring for notations to figure \ref{fig1}.11
$$s\pm s_g\pm\sum^n_{i=1}s_{g_i}\pm s_k\pm\sum^m_{j=1}s_{k_j}$$
where the signs depend on the choice of a coherent orientation,
$s_g$ is the saddle forming a square together with $n$, $s$ and $sn_1$,
$s_{g_i}$, with $1\leq i\leq N$ are the saddles in the chain of squares
of lenght $N$ attached to
$n$ and $s_g$ as in figure \ref{fig1}.7, $s_h$ is the saddle forming a square together with 
$n$, $s$ and $sn_2$, while
$s_{h_j}$, for $1\leq j\leq m$, are the saddles in the chain of squares
of lenght $m$ attached to
$n$ and $s_k$. The expression given for the generator includes all the cases:
indeed, if $s$ is not connected to any stable node, there are no squares, 
so the generator is simply $s$; if $s$ is connected to a unique stable node $sn_1$,
there is only the square with $s_g$, so the generator is $s\pm s_g$; if a chain
of squares of lenght $N$ is attached to $n$ and $s_g$, then the generator is
$s\pm s_g\pm\sum^N_{i=1}s_{g_i}$.

As to $HM_0$, observe that, considering only the unstable nodes $un$ and $n$, the element
$un$ in $Ker\partial_0^{U_2}$ is replaced by $un+n$ in $Ker\partial_0^{U_2}$.

As to $HM_2$, observe simply that $Im\partial_1^{U_1}=i(Im\partial_1^{U_2})$.

\item if the birth-death pair is given by a stable node and a saddle, then
$$dimKer\partial_1^{U_1}=dimKer\partial_1^{U_2}$$
$$dimIm\partial_0^{U_1}=dimIm\partial_0^{U_2}$$
moreover, if $un_1$ and $un_2$ are unstable nodes connected to $s$ in $U_1$,
writing $Im\partial_0^{U_2}=I\oplus<\partial_0^{U_2}(un_1),\partial_0^{U_2}(un_2)>$ for
some $I$, then
$$Im\partial_0^{U_1}=i(I)\oplus<\partial_0^{U_1}(un_1),\partial_0^{U_1}(un_2)>$$
instead, as to kernel of $\partial_1$,
suppose $s$ is connected, besides $n$, to a second stable node $sn$, assume that,
starting from $n$, there are finite 
chains of the form $(n,s_{i_1},sn_{i_1},s_{i_1i_2},sn_{i_1i_2},...,s_{i_1...i_p})$ ending
with a saddle $s_{i_1...i_p}$, where $1\leq i_r\leq i_r^{max}$ for $r=1,...,p$, 
and if from $sn$
there are finite 
chains of the form $(n,s^{'}_{j_1},sn^{'}_{j_1},s^{'}_{j_1j_2},sn^{'}_{j_1j_2},...,
s^{'}_{j_1...j_q})$ ending
with a saddle $s^{'}_{j_1...j_p}$, where $1\leq j_s\leq j_s^{max}$ for $s=1,...,q$, 
setting $\tilde{s}_{i_1...i_p}=\pm s_{i_1}\pm s_{i_1i_2}\pm
...\pm s_{i_1...i_p}$ and $\tilde{s}^{'}_{j_1...j_q}=\pm s_{j_1}\pm s_{j_1j_2}\pm
...\pm s_{j_1...j_q}$,
and writing 
\begin{eqnarray}
Ker\partial_1^{U_2} & = & K\oplus<\tilde{s}_{1...1}\pm \tilde{s}^{'}_{1...1},...,
\tilde{s}_{1...1}\pm 
\tilde{s}^{'}_{j_1^{max}...j_q^{max}},
\nonumber\\
 & & \tilde{s}_{21...1}\pm 
\tilde{s}^{'}_{1...1},...,\tilde{s}_{21...1}\pm 
\tilde{s}^{'}_{1...j_q^{max}},...,
\nonumber\\
 & & \tilde{s}_{i_1^{max}...i_p^{max}}\pm 
\tilde{s}^{'}_{1...1},...,
\tilde{s}_{i_1^{max}...i_p^{max}}\pm \tilde{s}^{'}_{j_1^{max}...j_q^{max}}>
\nonumber
\end{eqnarray}
for some $K$, then
\begin{eqnarray}
Ker\partial_1^{U_1} & = & i(K)\oplus<s\pm \tilde{s}_{1...1}\pm \tilde{s}^{'}_{1...1},...,
s\pm \tilde{s}_{1...1}\pm 
\tilde{s}^{'}_{j_1^{max}...j_q^{max}},
\nonumber\\
 & & s\pm \tilde{s}_{21...1}\pm \tilde{s}^{'}_{1...1},...,s\pm \tilde{s}_{21...1}\pm 
\tilde{s}^{'}_{1...j_q^{max}},...,
\nonumber\\
 & & s\pm \tilde{s}_{i_1^{max}...i_p^{max}}\pm \tilde{s}^{'}_{1...1},...,
s\pm \tilde{s}_{i_1^{max}...i_p^{max}}\pm \tilde{s}^{'}_{j_1^{max}...j_q^{max}}>
\nonumber
\end{eqnarray}
where signs depend, as usual, on a choice of a coherent orientation. If, as in the proof
of lemma \ref{causdiag}, we assume that all chains from $n$ and $sn$ are of lenght 1,
that is, $n$ is connected, besides $s$, to saddles $s_1$, ..., $s_n$ and these are not
connected to any other stable node, and, similarly, $n$ is connected, besides $s$, 
to saddles $s^{'}_1$, ..., $s^{'}_m$ and these too are not
connected to any other stable node, then the above expression simplifies considerably
\begin{eqnarray}
Ker\partial_1^{U_1} & = & i(K)\oplus<s\pm s_1\pm s^{'}_1,...,s\pm s_1\pm s^{'}_m,
\nonumber\\
 & & s\pm s_2\pm s^{'}_1,s\pm s_3\pm s^{'}_1,...,s\pm s_n\pm s^{'}_m>
\nonumber
\end{eqnarray}
Particularly, if $s$ is connected to a unique stable node, that is, only to $n$, we have
$$Ker\partial^{U_1}=i(K)\oplus<s\pm s_1,...,s\pm s_n>$$

As to $HM_0$, observe simply that $Ker\partial_0^{U_1}=i(Ker\partial_0^{U_2})$.

As to $HM_2$, observe that, though $i(\partial_1^{U_2}(s^{'}_{ki}))$, $i(\partial_1^{U_2}(s_{kj}))\in Im\partial_1^{U_1}$,
however $i(\partial_1^{U_2}(s^{'}_{ki}))\neq \partial_1^{U_1}(s^{'}_{ki})$ and 
$i(\partial_1^{U_2}(s_{kj}))=\partial_1^{U_1}(s_{kj})$.
\end{enumerate}
\end{remark}

Lemma \ref{causdiag} shows that $HM(U_1)\cong HM(U_2)$. The purpose now is to pick up
an isomorphism which will be used to glue the holomorphic structure of
the mirror object along the caustic $C$, providing the quantum corrections.

\begin{definition}
\label{qcc}
\rm
If $\partial U_1\cap\partial U_2\neq\varnothing$ 
contains only folds not limit points of
the bifurcation locus $B$, define an isomorphism $HM(U_1)\cong HM(U_2)$
as follows:
\begin{enumerate}
\item if the birth-death point is represented by an unstable node $n$ and a saddle $s$,
the isomorphism $M:HM(U_1)\cong HM(U_2)$ is the one induced in
homology by the map
$$\tilde{M}:\C[p_{1}^{U_2}]\oplus ... \oplus\C[p_{k}^{U_2}]\rightarrow
\C[p_{1}^{U_1}]\oplus ... \oplus\C[p_k^{U_1}]
\oplus\C[p_{k+1}^{U_1}]\oplus\C[p_{k+2}^{U_1}]$$
where $\tilde{M}$ is the natural injection $i$ on saddles and stable nodes
and, if $un_i,un$ are the unstable nodes appearing in $U_2$ with $un$ the
eventual second unstable node connected to $s$ in $U_1$ besides $n$, it is still 
the natural injection $i$
on $un_i$, while on $un$ it acts as
$$\tilde{M}(un)=un+n$$
\item If the birth-death point is represented by a stable node $n$ and a saddle $s$,
and if $s_i\neq s$, for $i=1,...,N$, are the saddles in the phase portraits 
over $U_1$ and $U_2$,
such that, for $i=1,...,m\leq N$, $s_i$ are the saddles connected to $n$, define
the isomorphism $M:HM(U_1)\cong HM(U_2)$ as the one induced in
homology by the map 
$$\tilde{M}:\C[p_{1}^{U_2}]\oplus ... \oplus\C[p_{k}^{U_2}]\rightarrow
\C[p_{1}^{U_1}]\oplus ... \oplus\C[p_k^{U_1}]
\oplus\C[p_{k+1}^{U_1}]\oplus\C[p_{k+2}^{U_1}]$$
such that
\begin{displaymath}
\tilde{M}(s_i)=\left\{ \begin{array}{lll}
s_i\pm s & , & 1\leq i\leq m\\
s_i & , & m<i\leq N
\end{array} \right.
\end{displaymath}
where the sign depends on orientation: 
$+$ if $n(\gamma_{s,n})=-n(\gamma_{s_i,n})$ and $-$ otherwise
($n(\gamma)$ denotes the sign of the gradient line $\gamma$),
while it is the natural injection $i$ on stable and unstable nodes.
\end{enumerate}
\end{definition}

To have a well-defined map $M$ in homology, it is necessary to check that  
$\tilde{M}$, defined on generators $p_i$, maps kernel and image of $\partial^{U_2}$
to, respectively, those of $\partial^{U_1}$.

\begin{lemma}
The map $\tilde{M}$ satisfies $\tilde{M}(Im\partial^{U_2})\subseteq Im\partial^{U_1}$
and $\tilde{M}(Ker\partial^{U_2})\subseteq Ker\partial^{U_1}$ and so it induces an
isomorphism $M:HM(U_1)\cong HM(U_2)$.
\end{lemma}
\begin{proof}
\begin{enumerate}
\item
When the birth-death pair is represented by an unstable node and a saddle,
the lemma follows from formulas in part 1 of remark \ref{remarkcau}.

\item
Suppose, instead, the birth-death pair is given by a stable node and a saddle.

\begin{itemize}
\item
The relations $\tilde{M}(Ker\partial_0^{U_2})\subseteq Ker\partial_0^{U_1}$ and
$\tilde{M}(Im\partial_1^{U_2})\subseteq Im\partial_1^{U_1}$ are easily verified. 

\item
So
consider the free $\C$-module generated by saddles.
Since $\tilde{M}$ is non trivial on those saddles $s_i$ connected to $n$, we are going
to compute $\tilde{M}(\partial_0^{U_2})$ 
on those unstable nodes connected to $s_i$.
If there are unstable nodes $un_1$ and $un_2$ connected to $s$,
let $s_i$ be the saddle forming squares with $s$, $n$ and $un_i$, for $i=1,2$.
Writing $\partial_0^{U_1}(un_i)=s\pm s_i\pm...$,
where, as usual, signs depend on orientation, 
then $\partial_0^{U_2}(un_i)=s_i\pm...$. 
So, in the special case where each node
$un_i$ is connected, in $U_1$, to only the saddle $s$ and $s_i$, 
since in a square 
$$n(\gamma_{s,n})=n(\gamma_{s_i,n})\Longleftrightarrow
n(\gamma_{un_i,s})\neq n(\gamma_{un_i,s_i})$$
it follows that $\tilde{M}(\partial_0^{U_2}(un_i))=\tilde{M}(s_i)=s_i\pm s$, and so
$\tilde{M}(Im\partial_0^{U_2})\subseteq Im\partial_0^{U_1}$.

If, instead, $un_i$ is connected to other saddles besides $s$ and $s_i$ (suppose, for
simplicity, only to the saddle $s_{k_i}$), then the formula
$\tilde{M}(Im\partial_0^{U_2})\subseteq Im\partial_0^{U_1}$ is still easily verified
provided that $s_{k_i}$ is not connected to $n$. 

In the opposite case, note,
first of all, that the phase portrait must exhibit other critical points,
at least a saddle, a stable node and an unstable node,
in order to be the phase portrait of a gradient vector field: this is shown in figure 
\ref{fig1}.17,
where it is represented only the unstable node $un_1$, and the critical points, added
in the phase portrait, are
denoted by $\tilde{s}$, $\tilde{sn}$, $\tilde{un}$.

\begin{center}
\setlength{\unitlength}{1cm}
\begin{picture}(8,6)
\label{fig17}
\thinlines

\put(1,1){\circle*{.1}}
\put(3,1){\circle*{.1}}
\put(1,2){\circle*{.1}}
\put(3,2){\circle*{.1}}
\put(3,3){\circle*{.1}}
\put(5,3){\circle*{.1}}
\put(5,2){\circle*{.1}}
\put(7,1){\circle*{.1}}

\scriptsize
\put(1,.7){$\displaystyle s$}
\put(3,.7){$\displaystyle n$}
\put(.4,2){$\displaystyle un_1$}
\put(3.1,2.1){$\displaystyle s_1$}
\put(3,3.1){$\displaystyle \tilde{sn}$}
\put(5.1,3.1){$\displaystyle \tilde{s}$}
\put(4.5,2.1){$\displaystyle \tilde{un}$}
\put(7,.7){$\displaystyle s_{1_1}$}
\put(1.9,1.4){$\displaystyle S_1$}
\put(5.9,3.8){$\displaystyle S_2$}
\normalsize

\qbezier(1,1)(3,1)(7,1)
\qbezier(1,5)(3,5)(7,5)
\qbezier(1,1)(1,2)(1,5)
\qbezier(7,1)(7,2)(7,5)
\qbezier(1,2)(3,2)(5,2)
\qbezier(3,3)(5,3)(6,3)
\qbezier(3,1)(3,2)(3,3)
\qbezier(6,1.5)(6,2)(6,3)
\qbezier(2,4)(3,4)(5,4)
\qbezier(2,4)(1.5,3)(1,2)
\qbezier(6,1.5)(4.5,1.25)(3,1)
\qbezier(5,2)(5,3)(5,4)

\qbezier(2,2)(1.95,2.05)(1.9,2.1)
\qbezier(2,2)(1.95,1.95)(1.9,1.9)

\qbezier(2,1)(1.95,1.05)(1.9,1.1)
\qbezier(2,1)(1.95,.95)(1.9,.9)

\qbezier(5,1)(5.05,1.05)(5.1,1.1)
\qbezier(5,1)(5.05,.95)(5.1,.9)

\qbezier(4,2)(4.05,2.05)(4.1,2.1)
\qbezier(4,2)(4.05,1.95)(4.1,1.9)

\qbezier(4,3)(4.05,3.05)(4.1,3.1)
\qbezier(4,3)(4.05,2.95)(4.1,2.9)

\qbezier(4,4)(3.95,4.05)(3.9,4.1)
\qbezier(4,4)(3.95,3.95)(3.9,3.9)

\qbezier(4,5)(3.95,5.05)(3.9,5.1)
\qbezier(4,5)(3.95,4.95)(3.9,4.9)

\qbezier(1,1.5)(1.05,1.55)(1.1,1.6)
\qbezier(1,1.5)(.95,1.55)(.9,1.6)

\qbezier(3,1.5)(3.05,1.55)(3.1,1.6)
\qbezier(3,1.5)(2.95,1.55)(2.9,1.6)

\qbezier(3,2.5)(3.05,2.55)(3.1,2.6)
\qbezier(3,2.5)(2.95,2.55)(2.9,2.6)

\qbezier(5,2.5)(5.05,2.55)(5.1,2.6)
\qbezier(5,2.5)(4.95,2.55)(4.9,2.6)

\qbezier(6,2.5)(6.05,2.55)(6.1,2.6)
\qbezier(6,2.5)(5.95,2.55)(5.9,2.6)


\end{picture}

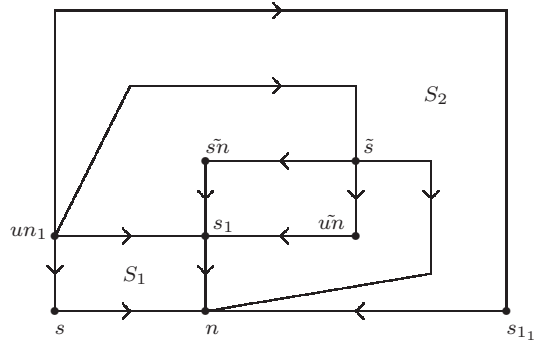
$Fig.~\ref{fig1}.17:~un_1~and~n~are~connected~to~the~same~saddle$
\end{center}

Considering the squares $S_1$ and $S_2$, it follows, respectively, that
$$n(\gamma_{s_1,n})=n(\gamma_{s,n})\Longleftrightarrow 
n(\gamma_{un_1,s_1})\neq n(\gamma_{un_1,s})$$
$$n(\gamma_{s_{1_1},n})=n(\gamma_{\tilde{s},n})\Longleftrightarrow 
n(\gamma_{un_1,s_{1_1}})\neq n(\gamma_{un_1,\tilde{s}})$$
so the quantum correction when crossing the caustic from $U_2$ to $U_1$ is given by
$$\tilde{M}(s_i)=s_i-\epsilon(s_i)s$$
for $i=1,3,4$, where
\begin{displaymath}
\epsilon(s_i) = \left\{ \begin{array}{lll}
1 & , & n(\gamma_{s_i,n})=n(\gamma_{s,n})\\
-1 & , & n(\gamma_{s_i,n})\neq n(\gamma_{s,n}) 
\end{array} \right.
\end{displaymath}

The Morse differential is defined as
$$\partial_0^{U_2}(un_1)=n(\gamma_{un_1,s_1})s_1+n(\gamma_{un_1,s_{1_1}})s_{1_1}+
n(\gamma_{un_1,\tilde{s}})\tilde{s}$$
$$\partial_0^{U_1}(un_1)=n(\gamma_{un_1,s})s+n(\gamma_{un_1,s_1})s_1
+n(\gamma_{un_1,s_{1_1}})s_{1_1}+
n(\gamma_{un_1,\tilde{s}})\tilde{s}$$
so it follows that
\begin{eqnarray}
\tilde{M}(\partial_0^{U_2}(un_1)) & = &n(\gamma_{un_1,s_1})s_1+n(\gamma_{un_1,s_{1_1}})s_{1_1}+
n(\gamma_{un_1,\tilde{s}})\tilde{s}+
\nonumber\\
 & & -[n(\gamma_{un_1,s_1})\epsilon(s_1)+n(\gamma_{un_1,s_{1_1}})\epsilon(s_{1_1})+
\nonumber\\
 & & +n(\gamma_{un_1,\tilde{s}})\epsilon(\tilde{s})]s
\nonumber
\end{eqnarray}

Since 
$$n(\gamma_{un_1,s_1})\epsilon(s_1)=-n(\gamma_{un_1,s})$$ 
while  
$$n(\gamma_{un_1,s_{1_1}})=n(\gamma_{un_1,\tilde{s}})\Leftrightarrow
n(\gamma_{s_{1_1},n})\neq n(\gamma_{\tilde{s},n})\Leftrightarrow
\epsilon(s_{1_1})\neq \epsilon(\tilde{s})$$
implies
$$n(\gamma_{un_1,s_{1_1}})\epsilon(s_{1_1})+
n(\gamma_{un_1,\tilde{s}})\epsilon(\tilde{s})=0$$
it follows that
$$\tilde{M}(\partial_0^{U_2}(un_1))=\partial_0^{U_1}(un_1)$$

Suppose now that another unstable node $\tilde{un}$ is connected to one of the saddle $s_i$,
say, for example, $s_1$
(the case considered in figure \ref{fig1}.17 is a particular case, more generally see figure 
\ref{fig1}.18). 
Suppose also that $\tilde{un}$ is connected to further saddles
$\tilde{s_h}$ for some parameter $h$.
It follows that $\tilde{M}(\partial_0^{U_2}(\tilde{un}))\in 
Im\partial_0^{U_1}$. By lemma \ref{php} there exists a saddle $s_q$ forming a square
with $\tilde{un}$, $s_1$ and $n$, and so $\tilde{M}(s_q)=s_q-\epsilon(s_q)s$

\begin{center}
\setlength{\unitlength}{1cm}
\begin{picture}(6,3)
\label{fig18}
\thinlines

\put(1,1){\circle*{.1}}
\put(1,2){\circle*{.1}}
\put(3,1){\circle*{.1}}
\put(3,2){\circle*{.1}}
\put(5,2){\circle*{.1}}
\put(5,1){\circle*{.1}}
\put(6,2){\circle*{.1}}

\scriptsize
\put(1.1,1.1){$\displaystyle s$}
\put(.7,2.1){$\displaystyle un_1$}
\put(3.1,2.1){$\displaystyle s_1$}
\put(3.1,1.1){$\displaystyle n$}
\put(4.9,2.1){$\displaystyle \tilde{un}$}
\put(5,.8){$\displaystyle s_q$}
\put(6.1,2.1){$\displaystyle \tilde{s_h}$}

\normalsize

\qbezier(1,1)(1,1.5)(1,2)
\qbezier(3,1)(3,1.5)(3,2)
\qbezier(1,1)(2,1)(3,1)
\qbezier(1,2)(2,2)(3,2)
\qbezier(3,2)(4,2)(5,2)
\qbezier(5,1)(4,1)(3,1)
\qbezier(5,1)(5,1.5)(5,2)
\qbezier(6,2)(5.5,2)(5,2)


\qbezier(1,1.5)(1.05,1.55)(1.1,1.6)
\qbezier(1,1.5)(.95,1.55)(.9,1.6)

\qbezier(3,1.5)(3.05,1.55)(3.1,1.6)
\qbezier(3,1.5)(2.95,1.55)(2.9,1.6)

\qbezier(2,2)(1.95,2.05)(1.9,2.1)
\qbezier(2,2)(1.95,1.95)(1.9,1.9)
\qbezier(2,1)(1.95,1.05)(1.9,1.1)
\qbezier(2,1)(1.95,.95)(1.9,.9)
\qbezier(4,2)(4.05,2.05)(4.1,2.1)
\qbezier(4,2)(4.05,1.95)(4.1,1.9)

\qbezier(4,1)(4.05,1.05)(4.1,1.1)
\qbezier(4,1)(4.05,.95)(4.1,.9)
\qbezier(5,1.5)(5.05,1.55)(5.1,1.6)
\qbezier(5,1.5)(4.95,1.55)(4.9,1.6)

\qbezier(5.5,2)(5.45,2.05)(5.4,2.1)
\qbezier(5.5,2)(5.45,1.95)(5.4,1.9)

\end{picture}

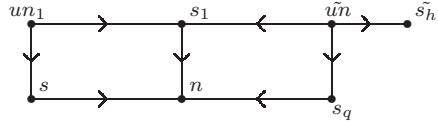
$Fig.~\ref{fig1}.18:~another~unstable~node~is~connected~to~s_1$
\end{center}

\begin{eqnarray}
\tilde{M}(\partial_0^{U_2}(\tilde{un})) & = & \tilde{M}(n(\gamma_{\tilde{un},s_1})s_1+
n(\gamma_{\tilde{un},s_q})s_q+
\nonumber\\
 & & +\sum_h n(\gamma_{\tilde{un},\tilde{s_h}})\tilde{s_h})=
\nonumber\\
 & = & n(\gamma_{\tilde{un},s_1})s_1+n(\gamma_{\tilde{un},s_q})s_q
+\sum_hn(\gamma_{\tilde{un},\tilde{s_h}})\tilde{s_h}+
\nonumber\\
 & & -[n(\gamma_{\tilde{un},s_1})\epsilon(s_1)+
n(\gamma_{\tilde{un},s_q})\epsilon(s_q)]s=
\nonumber\\
 & & n(\gamma_{\tilde{un},s_1})s_1+n(\gamma_{\tilde{un},s_q})s_q
+\sum_hn(\gamma_{\tilde{un},\tilde{s_h}})\tilde{s_h}=
\nonumber\\
 & = & \partial_0^{U_1}(\tilde{un})
\nonumber
\end{eqnarray}
because, as seen above, for a square always 
$$n(\gamma_{\tilde{un},s_1})\epsilon(s_1)+
n(\gamma_{\tilde{un},s_q})\epsilon(s_q)=0$$

As said at the beginning, 
the same conclusions are achieved if two unstable nodes $un_1$ and $un_2$ are connected
to $s$, and the proof does not change substantially.

It remains to check that $\tilde{M}(Ker\partial_1^{U_2})=Ker\partial_1^{U_1}$.
Observing that \\
$n(\gamma_{s_1,n})s_1-n(\gamma_{s_j,n})s_j$ is an element of both $Ker\partial_1^{U_1}$
and $Ker\partial_1^{U_2}$ for every $j=2,...,n$, and that
$n(\gamma_{s_1,n})s_1-n(\gamma_{s_j,n})s_j=\epsilon(s_1)s_1-\epsilon(s_j)s_j$,
it follows that
\begin{eqnarray}
Ker\partial_1^{U_1} & = & i(K)\oplus<s-\epsilon(s_1)s_1-\epsilon(s_1^{'})s_1^{'},...,
s-\epsilon(s_1)s_1+
\nonumber\\
 & & -\epsilon(s_m^{'})s_m^{'},\epsilon(s_1)s_1-\epsilon(s_2)s_2,...,
\epsilon(s_1)s_1-\epsilon(s_n)s_n>
\nonumber
\end{eqnarray}
and
\begin{eqnarray}
Ker\partial_1^{U_2} & = & K\oplus<-\epsilon(s_1)s_1-\epsilon(s_1^{'})s_1^{'},...,
-\epsilon(s_1)s_1-\epsilon(s_m^{'})s_m^{'},
\nonumber\\
 & & \epsilon(s_1)s_1-\epsilon(s_2)s_2,...,
\epsilon(s_1)s_1-\epsilon(s_n)s_n>
\nonumber
\end{eqnarray}
hence, since $\tilde{M}=i$ on all saddles but $s_j$,
\begin{eqnarray}
\tilde{M}(Ker\partial_1^{U_2}) & = & i(K)\oplus<-\epsilon(s_1)[s_1-\epsilon(s_1)s]
-\epsilon(s_1^{'})s_1^{'},...,
\nonumber\\
 & & -\epsilon(s_1)[s_1-\epsilon(s_1)s]-\epsilon(s_m^{'})s_m^{'},
\nonumber\\
 & & \epsilon(s_1)[s_1-\epsilon(s_1)s]-\epsilon(s_2)[s_2-\epsilon(s_2)s],...,
\nonumber\\
 & & \epsilon(s_1)[s_1-\epsilon(s_1)s]-\epsilon(s_n)[s_n-\epsilon(s_n)s]>=
\nonumber\\
 & = & Ker\partial_1^{U_1}
\nonumber
\end{eqnarray}
\end{itemize}
\end{enumerate}
\end{proof}

\section{The bifurcation locus}
\label{bifloc}
Points $x$ of the bifurcation locus $B$ 
are defined as those points at which a saddle-to-saddle separatrix 
appears in the phase portrait, so that to each bifurcation point $x$ it is associated
the pair of saddles $s_1(x)$ and $s_2(x)$ and the exceptional gradient 
line connecting them. Call them bifurcating saddles. 
In dimension 2, $B$ is an immersed submanifold 
of codimension 1. Bifurcation points characterized by the presence of two
saddle-to-saddle separatrixes form a codimension 2 subset. This endows $B$
with a stratification, where strata are given by codimension 1 and codimension 2
points (that is, points where one, respectively, two saddle-to-saddle separatrixes appear in
the phase portarit).

Consider two subset $U_1$ and $U_2$ such that $\partial U_1\cap \partial U_2\neq\varnothing$
is contained in a bifurcation line.
Note that if critical points are seen as sections of the fibration $\R^4\rightarrow\R^2$,
these are smooth sections over $B$.
The phase portraits over points of, respectively, $U_1$ and $U_2$ differ, in general, 
because of gradient lines appearing or disappearing between the bifurcating saddles
and some nodes (see figure \ref{fig19}.19). 
This changes the Morse differential $\partial$. However, even when $\partial$
is unchanged by the bifurcation, generically the phase portraits of $U_1$ and 
$U_2$ are not orbit isotopic. 
Consider, for simplicity, a 1-parameter family $\{f_t\}$ of functions such that $t=0$ is a 
bifurcation point.
Suppose there are two saddles $s_1$ and $s_2$ connected by a saddle-to-saddle separatrix 
$\gamma_{s_1,s_2}$.
Referring to figure \ref{fig19}.19,
if $W^u_j(s_i(t))$ and $W^s_j(s_i(t))$ denote, for $j=1,2$ the two components of the
unstable and respectively stable manifolds of $s_i$, $i=1,2$, then the family $\{f_t\}$
provides isotopies, for $t\leq0$, between $W^u_1(s_1(t))$ and $\gamma_{s_1,s_2}\cup
s_1\cup W^u_1(s_2(t))$ and between $W^s_1(s_2(t))$ and $W^s_2(s_1(t))\cup
s_1\cup\gamma_{s_1,s_2}$, for $t\geq0$, between $W^u_1(s_1(t))$ and $\gamma_{s_1,s_2}\cup
s_1\cup W^u_2(s_2(t))$ and between $W^s_1(s_2(t))$ and $W^s_1(s_1(t))\cup
s_1\cup\gamma_{s_1,s_2}$. Note that the above isotopy does not induce an orbit isotopy.
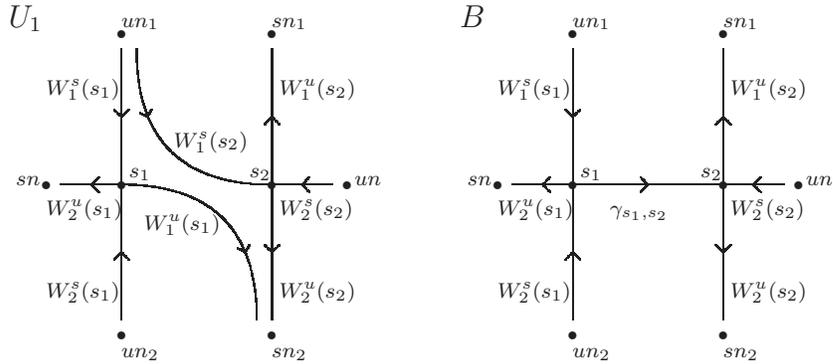
\begin{center}
\setlength{\unitlength}{1cm}
\begin{picture}(11,6)
\label{fig19}
\thinlines

\put(2,1){\circle*{.1}}
\put(2,3){\circle*{.1}}
\put(2,5){\circle*{.1}}
\put(4,1){\circle*{.1}}
\put(4,3){\circle*{.1}}
\put(4,5){\circle*{.1}}
\put(5,3){\circle*{.1}}
\put(1,3){\circle*{.1}}

\put(8,1){\circle*{.1}}
\put(8,3){\circle*{.1}}
\put(8,5){\circle*{.1}}
\put(10,1){\circle*{.1}}
\put(10,3){\circle*{.1}}
\put(10,5){\circle*{.1}}
\put(11,3){\circle*{.1}}
\put(7,3){\circle*{.1}}

\scriptsize
\put(2,.7){$\displaystyle un_2$}
\put(2.1,3.1){$\displaystyle s_1$}
\put(2,5.1){$\displaystyle un_1$}
\put(4,.7){$\displaystyle sn_2$}
\put(3.7,3.1){$\displaystyle s_2$}
\put(4,5.1){$\displaystyle sn_1$}
\put(.6,3){$\displaystyle sn$}
\put(5.1,3){$\displaystyle un$}
\put(1,1.5){$\displaystyle W^s_2(s_1)$}
\put(1,4.2){$\displaystyle W^s_1(s_1)$}
\put(4.1,1.5){$\displaystyle W^u_2(s_2)$}
\put(4.1,4.2){$\displaystyle W^u_1(s_2)$}
\put(1,2.6){$\displaystyle W^u_2(s_1)$}
\put(4.1,2.6){$\displaystyle W^s_2(s_2)$}
\put(2.7,3.5){$\displaystyle W^s_1(s_2)$}
\put(2.3,2.4){$\displaystyle W^u_1(s_1)$}

\put(8,.7){$\displaystyle un_2$}
\put(8.1,3.1){$\displaystyle s_1$}
\put(8,5.1){$\displaystyle un_1$}
\put(10,.7){$\displaystyle sn_2$}
\put(9.7,3.1){$\displaystyle s_2$}
\put(10,5.1){$\displaystyle sn_1$}
\put(6.6,3){$\displaystyle sn$}
\put(11.1,3){$\displaystyle un$}
\put(7,1.5){$\displaystyle W^s_2(s_1)$}
\put(7,4.2){$\displaystyle W^s_1(s_1)$}
\put(10.1,1.5){$\displaystyle W^u_2(s_2)$}
\put(10.1,4.2){$\displaystyle W^u_1(s_2)$}
\put(7,2.6){$\displaystyle W^u_2(s_1)$}
\put(10.1,2.6){$\displaystyle W^s_2(s_2)$}
\put(8.5,2.6){$\displaystyle \gamma_{s_1,s_2}$}

\normalsize

\put(.5,5.1){$\displaystyle U_1$}

\put(6.5,5.1){$\displaystyle B$}

\qbezier(1.2,3)(1.6,3)(2,3)
\qbezier(4.8,3)(4.4,3)(4,3)
\qbezier(2,1.2)(2,3)(2,4.8)
\qbezier(4,1.2)(4,3)(4,4.8)
\qbezier(2.2,4.8)(2.2,3)(4,3)
\qbezier(3.8,1.2)(3.8,3)(2,3)

\qbezier(7.2,3)(7.6,3)(8,3)
\qbezier(10.8,3)(10.4,3)(10,3)
\qbezier(8,1.2)(8,3)(8,4.8)
\qbezier(10,1.2)(10,3)(10,4.8)
\qbezier(10,3)(9,3)(8,3)


\qbezier(1.6,3)(1.65,3.05)(1.7,3.1)
\qbezier(1.6,3)(1.65,2.95)(1.7,2.9)

\qbezier(4.4,3)(4.45,3.05)(4.5,3.1)
\qbezier(4.4,3)(4.45,2.95)(4.5,2.9)

\qbezier(2,3.9)(2.05,3.95)(2.1,4)
\qbezier(2,3.9)(1.95,3.95)(1.9,4)

\qbezier(4,3.9)(4.05,3.85)(4.1,3.8)
\qbezier(4,3.9)(3.95,3.85)(3.9,3.8)

\qbezier(2,2.1)(2.05,2.05)(2.1,2)
\qbezier(2,2.1)(1.95,2.05)(1.9,2)

\qbezier(4,2.1)(4.05,2.15)(4.1,2.2)
\qbezier(4,2.1)(3.95,2.15)(3.9,2.2)

\qbezier(2.35,3.9)(2.3,3.925)(2.25,3.95)
\qbezier(2.35,3.9)(2.375,3.95)(2.4,4)

\qbezier(3.65,2.1)(3.675,2.15)(3.7,2.2)
\qbezier(3.65,2.1)(3.6,2.125)(3.55,2.15)

\qbezier(7.6,3)(7.65,3.05)(7.7,3.1)
\qbezier(7.6,3)(7.65,2.95)(7.7,2.9)

\qbezier(10.4,3)(10.45,3.05)(10.5,3.1)
\qbezier(10.4,3)(10.45,2.95)(10.5,2.9)

\qbezier(8,3.9)(8.05,3.95)(8.1,4)
\qbezier(8,3.9)(7.95,3.95)(7.9,4)

\qbezier(10,3.9)(10.05,3.85)(10.1,3.8)
\qbezier(10,3.9)(9.95,3.85)(9.9,3.8)

\qbezier(8,2.1)(8.05,2.05)(8.1,2)
\qbezier(8,2.1)(7.95,2.05)(7.9,2)

\qbezier(10,2.1)(10.05,2.15)(10.1,2.2)
\qbezier(10,2.1)(9.95,2.15)(9.9,2.2)

\qbezier(9,3)(8.95,3.05)(8.9,3.1)
\qbezier(9,3)(8.95,2.95)(8.9,2.9)

\end{picture}

\setlength{\unitlength}{1cm}
\begin{picture}(6,6)

\thinlines

\put(2,1){\circle*{.1}}
\put(2,3){\circle*{.1}}
\put(2,5){\circle*{.1}}
\put(4,1){\circle*{.1}}
\put(4,3){\circle*{.1}}
\put(4,5){\circle*{.1}}
\put(5,3){\circle*{.1}}
\put(1,3){\circle*{.1}}

\scriptsize
\put(2,.7){$\displaystyle un_2$}
\put(2.1,3.1){$\displaystyle s_1$}
\put(2,5.1){$\displaystyle un_1$}
\put(4,.7){$\displaystyle sn_2$}
\put(3.7,3.1){$\displaystyle s_2$}
\put(4,5.1){$\displaystyle sn_1$}
\put(.6,3){$\displaystyle sn$}
\put(5.1,3){$\displaystyle un$}
\put(1,1.5){$\displaystyle W^s_2(s_1)$}
\put(1,4.2){$\displaystyle W^s_1(s_1)$}
\put(4.1,1.5){$\displaystyle W^u_2(s_2)$}
\put(4.1,4.2){$\displaystyle W^u_1(s_2)$}
\put(1,2.6){$\displaystyle W^u_2(s_1)$}
\put(4.1,2.6){$\displaystyle W^s_2(s_2)$}
\put(2.3,3.5){$\displaystyle W^u_1(s_1)$}
\put(2.8,2.4){$\displaystyle W^s_1(s_2)$}

\normalsize

\put(.5,5.1){$\displaystyle U_2$}

\qbezier(1.2,3)(1.6,3)(2,3)
\qbezier(4.8,3)(4.4,3)(4,3)
\qbezier(2,1.2)(2,3)(2,4.8)
\qbezier(4,1.2)(4,3)(4,4.8)
\qbezier(2.2,1.2)(2.2,3)(4,3)
\qbezier(3.8,4.8)(3.8,3)(2,3)


\qbezier(1.6,3)(1.65,3.05)(1.7,3.1)
\qbezier(1.6,3)(1.65,2.95)(1.7,2.9)

\qbezier(4.4,3)(4.45,3.05)(4.5,3.1)
\qbezier(4.4,3)(4.45,2.95)(4.5,2.9)

\qbezier(2,3.9)(2.05,3.95)(2.1,4)
\qbezier(2,3.9)(1.95,3.95)(1.9,4)

\qbezier(4,3.9)(4.05,3.85)(4.1,3.8)
\qbezier(4,3.9)(3.95,3.85)(3.9,3.8)

\qbezier(2,2.1)(2.05,2.05)(2.1,2)
\qbezier(2,2.1)(1.95,2.05)(1.9,2)

\qbezier(4,2.1)(4.05,2.15)(4.1,2.2)
\qbezier(4,2.1)(3.95,2.15)(3.9,2.2)

\qbezier(2.35,2.1)(2.3,2.075)(2.25,2.05)
\qbezier(2.35,2.1)(2.375,2.05)(2.4,2)

\qbezier(3.65,3.9)(3.675,3.85)(3.7,3.8)
\qbezier(3.65,3.9)(3.6,3.875)(3.55,3.85)

\end{picture}
$Fig.~\ref{fig19}.19:~a~saddle-to-saddle~separatrix$
\end{center}

This implies that if there is an unstable node $un_1$ such that $\gamma_{un_1,s_1}=
W^s_1(s_1(t))$ for $t<0$, then for $t>0$, besides 
$\gamma_{un_1,s_1}=W^s_1(s_1(t))$, also the gradient line $\gamma_{un_1,s_2}=
W^s_1(s_2(t))$ appears in the phase portrait; 
if there exists an unstable node $un_2$ such that $\gamma_{un_2,s_1}=
W^s_2(s_1(t))$ for $t<0$, then there exists also the gradient line $\gamma_{un_2,s_2}=
W^s_1(s_2(t))$ for $t<0$, which breaks for $t>0$, while $\gamma_{un_2,s_1}=
W^s_2(s_1(t))$ persists in the phase portrait; if there is a stable node $sn_1$
such that $\gamma_{s_2,sn_2}=
W^u_1(s_2(t))$ for $t<0$, then there is also the gradient line $\gamma_{s_1,sn_2}=
W^u_1(s_1(t))$ for $t<0$, which breaks for $t>0$, while $\gamma_{s_2,sn_2}=
W^u_1(s_2(t))$ persists in the phase portrait; 
if there exists a stable node $sn_2$ such that $\gamma_{s_2,sn_2}=
W^s_2(s_2(t))$ for $t<0$, then for $t>0$ there is also, besides 
$\gamma_{s_2,sn_2}=W^s_2(s_2(t))$, the gradient line $\gamma_{s_1,sn_2}=
W^u_1(s_1(t))$; instead, if there are a stable node $sn$ such that $\gamma_{s_1,sn}=
W^u_2(s_1(t))$ or an unstable node $un$ such that $\gamma_{un,s_2}=
W^s_2(s_2(t))$, these lines appear in both the phase portrait for $t<0$ and $t>0$.
In this sense gradient lines between the bifurcating saddles
and some nodes appear or break in the phase portraits over points of $U_1$ and 
$U_2$ when crossing the bifurcation locus, changing the Morse differential $\partial$. The purpose is to prove that
the Morse complexes over $U_1$ and $U_2$ are isomorphic, and to pick up a suitable
isomorphism, providing the quantum correction.

\begin{lemma}
\label{signs}
In the situation described above of two bifurcating saddles,
suppose there exists, for instance in $U_1$, 
an unstable node $un_1$ with gradient lines $\gamma_{un_1,s_1}$
and $\gamma_{un_1,s_2}$, and 
assume that the signs given to the separatrixes of the saddles $s_1$ and $s_2$ coincide
in $U_1$ and $U_2$:
\begin{enumerate}
\item suppose there is also a stable node $sn_2$ with gradient lines $\gamma_{s_1,sn_2}$
and $\gamma_{s_2,sn_2}$, then
$$n(\gamma_{un_1,s_1})=n(\gamma_{un_1,s_2})\Longleftrightarrow 
n(\gamma_{s_1,sn_2})=-n(\gamma_{s_2,sn_2})$$
\item $n(\gamma_{un_1,s_1})=n(\gamma_{un_1,s_2})$ in $U_1$ if and only if
$n(\gamma_{un_2,s_1})=-n(\gamma_{un_2,s_2})$ and/or
$n(\gamma_{s_1,sn_1})=n(\gamma_{s_2,sn_1})$ in $U_2$ (provided these gradient lines do exist).
\end{enumerate}
\end{lemma}
\begin{proof}
\begin{enumerate}
\item It is a consequence of a choice of an orientation for a square
and already used along the proof of lemma \ref{causdiag}.
\item Suppose there exists in $U_1$ a second unstable node $un_2$ with the gradient line
$\gamma_{un_2,s_1}$, then in $U_2$, besides $\gamma_{un_1,s_1}$, 
both the lines $\gamma_{un_2,s_1}$ and
$\gamma_{un_2,s_2}$ appear in the phase portrait. Observe now that since
$n(W^u_1(s_1))=-n(W^u_2(s_1))$ and, by assumption, $n(\gamma_{un_1,s_2})=
n(\gamma_{un_2,s_2})$, it follows that $n(\gamma_{un_2,s_1})=-n(\gamma_{un_2,s_2})$.
Suppose there exist in $U_2$ a stable node $sn_1$ with the gradient line 
$\gamma_{s_2,sn_1}$, then in $U_2$, besides $\gamma_{un_1,s_1}$, also the lines 
$\gamma_{s_1,sn_1}$ and $\gamma_{s_2,sn_1}$ occur in the phase portrait.
By lemma \ref{php} there exist a saddle $s$ forming a square in $U_1$ with
$un_1$, $s_2$ and $sn_1$, and in $U_2$ with $un_1$, $s_1$ and $sn_1$.
If $n(\gamma_{un_1,s_2})=n(\gamma_{s_2,sn_1})$ in $U_1$ then 
$n(\gamma_{un_1,s})=-n(\gamma_{s,sn_1})$, which implies, assuming that
signs of separatrixes coincide in $U_1$ and $U_2$, that
$n(\gamma_{s_1,sn_1})=n(\gamma_{un_1,s_1})$ in $U_2$. As, by hypothesis,
$n(\gamma_{un_1,s_1})=n(\gamma_{un_1,s_2})$, it follows that
$n(\gamma_{s_1,sn_1})=n(\gamma_{s_2,sn_1})$. If, instead,
$n(\gamma_{un_1,s_2})=-n(\gamma_{s_2,sn_1})$ in $U_1$ then 
$n(\gamma_{un_1,s})=n(\gamma_{s,sn_1})$, which implies, that
$n(\gamma_{s_1,sn_1})=n(\gamma_{un_1,s_1})$ in $U_2$. As, by hypothesis,
$n(\gamma_{un_1,s_1})=n(\gamma_{un_1,s_2})$, it follows again that
$n(\gamma_{s_1,sn_1})=n(\gamma_{s_2,sn_1})$.
\end{enumerate}
\end{proof}

%
Let $s_1$, $s_2$, $s_3$, ..., $s_k$ be the saddles in the phase portrait, where
$s_1$ and $s_2$ form a pair of bifurcating saddles. Assume to choose signs according
to lemma \ref{signs} (this will be made more rigorous later in definition \ref{bifor}).

\begin{lemma}
\label{bifdiag}
Let $U_1$ and $U_2$ be open subsets such that $U_i\cap(B\cup C)=\varnothing$
and $\partial U_1\cap\partial U_2\subset B$ consists only of 
codiemension 1 bifurcation points.
Then the homology groups of the Morse complexes 
$$\C[p_{1}^{U_1}]\oplus ... \oplus\C[p_{k}^{U_1}]$$
$$\C[p_{1}^{U_2}]\oplus ... \oplus\C[p_k^{U_2}]$$
are isomorphic.
\end{lemma}
\begin{proof}
We are going to prove that $dimKer\partial^{U_1}=dimKer\partial^{U_2}$ and
$dimIm\partial^{U_1}\\
=dimIm\partial^{U_2}$. 

\begin{itemize}
\item
It is easily verified that $dimKer\partial_0^{U_1}=dimKer\partial_0^{U_2}$ and
$dimIm\partial_1^{U_1}
=dimIm\partial_1^{U_2}$.

\item
It remains to check the above relations at the level of 1-chains, that is, on saddles.
With the notation as in 
figure \ref{fig19}.19, what can change the Morse differential $\partial$
in $U_1$ and $U_2$ are the gradient lines between the saddle $s_1$ and $s_2$ on one side
and the nodes $un_1$, $un_2$, $sn_1$ and $sn_2$ on the other, thus, for simplicity,  
assume only these nodes in the phase portrait. 

\item[-]
Consider, first of all, how
$Im\partial_0$ can change when crossing the bifurcation line.

If neither $un_1$ nor $un_2$ are in the phase portrait then
$Im\partial_0^{U_1}=Im\partial_0^{U_2}=\{0\}$.

If just $un_1$ is in the phase portrait, then, as explained, $U_1$ exhibits the
gradient lines $\gamma_{un_1,s_1}$ and $\gamma_{un_1,s_2}$ while $U_2$ only
$\gamma_{un_1,s_1}$. So, for a choice of signs as lemma \ref{signs} suggets (see further on 
definition \ref{bifor}), we have that
$$Im\partial_0^{U_1}=<s_1+s_2\pm\sum_{j=1}^{m(un_1)} s_j>$$
$$Im\partial_0^{U_2}=<s_1\pm\sum_{j=1}^{m(un_1)} r_j>$$
where $r_j$, for $j=1,...,m(un_1)$, are further saddles connected to $un_1$,
hence 
$$dimIm\partial_0^{U_1}=dimIm\partial_0^{U_2}$$.

If just $un_2$ is in the phase portrait, then $U_1$ exhibits only the
gradient line $\gamma_{un_1,s_1}$ while $U_2$ both
$\gamma_{un_1,s_1}$ and $\gamma_{un_1,s_2}$. 
So, according to 
lemma \ref{signs}, 
$$Im\partial_0^{U_1}=<s_1\pm\sum_{l=1}^{m(un_2)} s_j>$$
$$Im\partial_0^{U_2}=<s_1-s_2\pm\sum_{l=1}^{m(un_2)} t_l>$$
where $t_l$, for $l=1,...,m(un_2)$, are further saddles connected to $un_2$.
Hence still $dimIm\partial_0^{U_1}=dimIm\partial_0^{U_2}$.

If both $un_1$ and $un_2$ are in the phase portrait, then in $U_1$ we have the
gradient lines $\gamma_{un_1,s_1}$, $\gamma_{un_2,s_1}$ and $\gamma_{un_1,s_2}$,
while in $U_2$ the lines $\gamma_{un_1,s_1}$, $\gamma_{un_2,s_1}$ and $\gamma_{un_2,s_2}$.
Lemma \ref{signs} implies that
$$Im\partial_0^{U_1}=<s_1+s_2\pm\sum_{j=1}^{m(un_1)} r_j,s_1\pm\sum_{l=1}^{m(un_2)} t_l>$$
$$Im\partial_0^{U_2}=<s_1\pm\sum_{j=1}^{m(un_1)} r_j,s_1-s_2\pm\sum_{l=1}^{m(un_2)} t_l>$$
so $dimIm\partial_0^{U_1}=dimIm\partial_0^{U_2}$.

\item[-]
Consider, now, how
$Ker\partial_1$ can change when crossing the bifurcation line.

If neither $sn_1$ nor $sn_2$ are in the phase portrait then
$$Ker\partial_1^{U_1}=Ker\partial_1^{U_2}=<s_1,s_2>$$ 
so $dimKer\partial_1^{U_1}=dimKer\partial_1^{U_2}$.

If just $sn_1$ is in the phase portrait, then in $U_1$ there is only the
gradient line $\gamma_{s_2,sn_1}$ while in $U_2$ we have $\gamma_{s_1,sn_1}$ and
$\gamma_{s_2,sn_1}$. 
Lemma \ref{signs} implies that
$$Ker\partial_1^{U_1}=<s_1,s_2\pm x_1,...,s_2\pm x_{m(sn_1)}>$$
$$Ker\partial_1^{U_2}=<s_1-s_2,s_2\pm x_1,...,s_2\pm x_{m(sn_1)}>$$
where $x_p$, for $p=1,...,m(sn_1)$, are further saddles connected to $sn_1$.
Therefore $dimKer\partial_1^{U_1}=dimKer\partial_1^{U_2}$.

If just $sn_2$ is in the phase portrait, then in $U_1$ there are the
gradient lines $\gamma_{s_1,sn_2}$ and $\gamma_{s_2,sn_2}$, while in $U_2$ only
$\gamma_{s_2,sn_2}$. 
Lemma \ref{signs} implies that 
$$Ker\partial_1^{U_1}=<s_1+s_2,s_2\pm y_1,...,s_2\pm y_{m(sn_2)}>$$
$$Ker\partial_1^{U_2}=<s_1,s_2\pm y_1,...,s_2\pm y_{m(sn_2)}>$$
where $y_q$, for $q=1,...,m(sn_2)$, are further saddles connected to $sn_2$,
in particular, it follows that $dimKer\partial_1^{U_1}=dimKer\partial_1^{U_2}$.

If both $sn_1$ and $sn_2$ are in the phase portrait, then this exhibits over $U_1$ the
gradient lines $\gamma_{s_1,sn_2}$, $\gamma_{s_2,sn_1}$ and $\gamma_{s_2,sn_2}$, while
over $U_2$ the gradient 
lines $\gamma_{s_1,sn_1}$, $\gamma_{s_2,sn_1}$ and $\gamma_{s_2,sn_2}$.
Lemma \ref{signs} implies that
\begin{eqnarray}
Ker\partial_1^{U_1} & = & <s_1+s_2\pm x_1,s_2\pm y_1\pm x_1,...,s_2\pm y_{m(sn_2)}+x_1,
\nonumber\\
 & & s_2\pm y_1\pm x_2,...,s_2\pm y_1\pm x_{m(sn_1)}>
\nonumber
\end{eqnarray}
\begin{eqnarray}
Ker\partial_1^{U_1} & = & <s_1\pm x_1,s_2\pm y_1\pm x_1,...,s_2\pm y_{m(sn_2)}+x_1,
\nonumber\\
 & & s_2\pm y_1\pm x_2,...,s_2\pm y_1\pm x_{m(sn_1)}>
\nonumber
\end{eqnarray}
so $dimKer\partial_1^{U_1}=dimKer\partial_1^{U_2}$.
\end{itemize}
\end{proof}

\begin{remark}
\label{bifrem}
\rm
The presence of a second stable node $sn$ connected to $s_1$ or of a second
unstable node $un$ connected to $s_2$ simply adds new terms in the expressions
of $Im\partial_0$ and $Ker\partial_1$, which, however, are not modified by the bifurcation
and so appears both in $U_1$ and in $U_2$.
\end{remark}

We now pick up an isomorphism between the Morse homologies in $U_1$ and $U_2$.

\begin{definition}
\label{qcb}
\rm
For a bifurcation characterized by the appearance of the saddle-to-saddle separatrix
$\gamma_{s_1,s_2}$ and whose locus is a line $B$, and, using notations as in figure 
\ref{fig19}.19,
if orientation is chosen in such a way that in $U_1$, as lemma \ref{signs} suggests,
$W^s_1(s_1)$ and $W^s_1(s_2)$ have same sign and $W^u_1(s_1)$ and $W^u_2(s_2)$ have
opposite sign,
define
an isomorphism $M:HM(U_1)\rightarrow HM(U_2)$ as the one induced in homology by the map
$\tilde{M}:\oplus^k_{i=1}\C[s_{i}^{U_1}]\rightarrow\oplus^k_{i=1}\C[s_{i}^{U_2}]$ such that
\begin{displaymath}
\tilde{M}(s_i)=\left\{ \begin{array}{lll}
s_1-s_2 & , & i=1\\
s_i & , & i\neq1
\end{array} \right.
\end{displaymath}
$\tilde{M}=Id$ on nodes
\end{definition}

\begin{remark}
Note that if in $U_1$ signs are not as definition \ref{qcb} requires, by lemma 
\ref{signs} this condition is fulfilled in $U_2$, so that $\tilde{M}$ defines
an isomorphism $M:HM(U_2)\rightarrow HM(U_1)$. The isomorphism $HM(U_2)\rightarrow HM(U_1)$
is provided by $\tilde{M}^{-1}$
\begin{displaymath}
\tilde{M}^{-1}(s_i)=\left\{ \begin{array}{lll}
s_1+s_2 & , & i=1\\
s_i & , & i\neq1
\end{array} \right.
\end{displaymath}
\end{remark}
%

\begin{lemma}
The map $\tilde{M}$ induces a map in homology.
\end{lemma}
\begin{proof}
\begin{itemize}
\item
Nothing to prove as to $HM_0$ and $HM_2$, being $\tilde{M}=Id$ on nodes. 

\item
So consider $HM_1$.
\item[-]
Lemma \ref{bifdiag} implies that $\tilde{M}(Ker\partial_1^{U_1})=Ker\partial_1^{U_2}$
and $\tilde{M}(Im\partial_0^{U_1})=Im\partial_0^{U_2}$,
when the only nodes in the phase portrait are those named $un_1$, $un_2$, $sn_1$ and
$sn_2$. However, as noted in remark \ref{bifrem}, there could be a further stable node $sn$
connected to $s_1$ and a further unstable node $un$ connected to $s_2$ In this case,
it is necessary to check that $\partial_0^{U_1}(un)$ is mapped by $\tilde{M}$ into
$Im\partial_0^{U_2}$:
since $\partial_0^{U_1}(un)=
\partial_0^{U_2}(un)=s_2+\sum_{j=1}^m r_j$ for some saddles $r_j$ with $r_j\neq s_1$ for
all $j=1,...,m$, and since $\tilde{M}$ acts as the identity on $s_2$ and $r_j$, it
follows that $\tilde{M}(\partial_0^{U_1}(un))=\partial_0^{U_2}(un)$. 
\item[-]
As to $Ker\partial_1$,
consider, first, the following case:
$s_2$ is not connected to any stable node, while to $sn$ it is connected, besides $s_1$, 
another saddle $s$; then 
$Ker\partial_1^{U_1}=Ker\partial_1^{U_2}=<s_1\pm s,s_2>$ and so 
$\tilde{M}(Ker\partial_1^{U_1})=Ker\partial_1^{U_2}$. If, instead, $s_2$ is connected to
two stable nodes $sn_1$ and $sn_2$, which, besides $s_2$, are connected, respectively,
to further saddles $r_1$ and $r_2$, we have that
$Ker\partial_1^{U_1}=<s_1+s-r_1,s_2+r_1+r_2>$ and
$Ker\partial_1^{U_2}=<s_1+s+r_2,s_2+r_1+r_2>$; thus, being
$\tilde{M}(s_1+s+r_1)=s_1+s+r_1-s_2=(s_1+s+r_2)-(s_2+r_1+r_2)$, it follows that
$\tilde{M}(Ker\partial_1^{U_1})=Ker\partial_1^{U_2}$. The argument is independent from
the chosen orientation. The same conclusion is achieved, modifying slightly
the proof, if just one stable node is connected to $s_2$.   
\end{itemize}
\end{proof}

To apply definition \ref{qcb} it is necessary that orientation is chosen in
a proper way. The following definition selects the class of orientations
for which quantum corrections can be constructed, in accordance with definition
\ref{qcb}.

\begin{definition}
\label{bifor}
\rm
An orientation in a phase portrait is said to satisfy
the ``signs convention'' if and only if it is chosen in such a way that,
for any pair of bifurcating
saddles $s_1$ and $s_2$ exhibiting the saddle-to-saddle separatrix $\gamma_{s_1,s_2}$
along a certain bifurcation line, 
$W^s_1(s_1)$ and $W^s_1(s_2)$ are given the same sign, while
$W^u_1(s_1)$ and $W^u_2(s_2)$ opposite sign (see figure \ref{fig19}.19 for notation), 
\end{definition}

The following proposition states that such class of orientations is not empty.

\begin{proposition}
For any given phase portrait, 
there exists a coherent orientation satisfying definition \ref{bifor}.
\end{proposition}
\begin{proof}
The coherent orientation of definition \ref{bifor} is the orientation 
corresponding to a phase portrait where
stable and unstable nodes $sn$ and $un$ are added in a such a way that, for each pair of
bifurcation saddles $s_1$ and $s_2$, $W^s_1(s_1)$ and $W^s_1(s_2)$ connect
$un$ with respectively $s_1$ and $s_2$, and $W^u_1(s_1)$ and $W^u_2(s_2)$ connect
respectively $s_1$ and $s_2$ to $sn$: in this case, indeed, $(un,s_1,s_2,sn)$, with
these separatrixes, form
a square, and this is just the orientation for a square. The existence of
coherent orientations for any phase portrait proves now the proposition.
\end{proof}

\section{Intersection of caustic and bifurcation lines}
\label{int}
As explained,
the caustic $C$ is an immersed submanifold of codimension 1
having, in dimension 2, two strata: the folds, forming the stratum of
codimension 1, and the cusps, the stratum of codimension 2.
Different branches of $C$ can intersect transversely one with another, 
generically at folds, to which
corresponds two birth-death pairs with no common points.

The bifurcation locus $B$ is as well an immersed submanifold of codimension 2
with two strata: codimension 1 and codimension 2 bifurcations, where the
corresponding phase portrait exhibits
one or respectively two saddle-to-saddle separatrixes. This means that
the intersection points of two bifurcation lines are codimension 2 bifurcations,
each line representing one of the two saddle-to-saddle separatrixes exhibited
by the codimension 2 bifurcation (see \cite{M1} and \cite{M2}).

A bifurcation line $B$ can intersect the caustic $C$, generically, at a fold:
indeed, if the intersection were a cusp, the exceptional gradient line appearing
at this point will break for any small perturbation, and, by a transversality
argument (see for example \cite{M2} for perturbations of the elliptic umbilic), 
will appear at a nearby point of the caustic, which, generically, is a fold.
Actually, to be precise, having defined bifurcation points away from
the caustic, it should be better to talk about points of the caustic being limit points
of the bifurcation locus rather than intersection points of the caustic and bifurcation
locus. 
If to $B$ it is associated the saddle-to-saddle separatrix $\gamma_{s_1,s_2}$
and to $C$ the pair of birth-death points $(s_i,n)$, where $i=1$ or $i=2$ and
$n$ is a node, then, in a neighbourhood of the intersection point, 
$B$ is an half-line lying in one of the two subsets determined by $C$,
precisely the one exhibiting $s_i$, and whose
origin is the intersection point. As already 
explained above, generically this point is a fold.
The two ways $B$ can meet $C$ are shown in figure \ref{fig20}.20:
case (a) was described just above, case (b) occurs when $i\neq1,2$.

\begin{center}
\setlength{\unitlength}{1cm}
\begin{picture}(8,4)
\label{fig20}
\thinlines

\put(2,2){\circle*{.1}}
\put(6,2){\circle*{.1}}

\scriptsize
\put(.7,1.9){$\displaystyle C$}
\put(1.9,.7){$\displaystyle B$}
\put(4.7,1.9){$\displaystyle C$}
\put(5.9,.7){$\displaystyle B$}

\normalsize
\put(.4,3){$\displaystyle (a)$}
\put(4.4,3){$\displaystyle (b)$}

\qbezier(1,2)(2,2)(3,2)
\qbezier(2,2)(2,1.5)(2,1)

\qbezier(5,2)(6,2)(7,2)
\qbezier(6,3)(6,2)(6,1)

\end{picture}
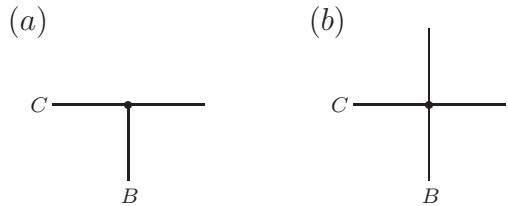
$Fig.~\ref{fig20}.20:~Intersections~between~the~caustic~and~the~bifurcation~locus$
\end{center}

Consider now the intersection of two bifurcation lines $B_1$ and $B_2$ (see \cite{M2}
for examples concerning the perturbed elliptic umbilic).
The phase portrait corresponding to the intersection point $z$ of $B_1$ and $B_2$
contains two saddle-to-saddle separatrixes $\gamma_1$ and $\gamma_2$, 
each appearing, respectively, 
along $B_1$ and $B_2$. The phase portrait associated to each subset
determined by $B_1$ and $B_2$ is obtained from the phase portrait over $z$ by breaking
$\gamma_1$ and $\gamma_2$: if $\gamma_1\neq\gamma_2$ there are at least four
of such phase portraits, which are not orbit
isotopic.

%
To study when the assignement of $B_1$ and $B_2$, 
together with the exceptional gradient lines $\gamma_1$ and $\gamma_2$ which they represent, 
can give rise to an admissible
$CB$-diagram, it is necessary to distinguish
between two cases. In fact, some attention must be paid when $B_1$ and $B_2$
represent bifurcations with saddle-to-saddle separatrix $\gamma_{s_1,s_2}$
and $\gamma_{s_2,s_3}$ respectively.

\begin{lemma}
\label{intbif}
Let $B_1$ and $B_2$ be bifurcation lines
representing bifurcations corresponding to saddle-to-saddle separatrixes $\gamma_{s_1,s_2}$
and $\gamma_{s_3,s_4}$, with either $s_2\neq s_3$ or $s_4\neq s_1$; 
if the phase portrait corresponding 
to the intersection point $B_1\cap B_2$ is the phase portrait of a gradient vector field,
and if to each subset determined by $B_1$ and $B_2$ it is associated a phase
portrait obtained, in the way described above, by breaking the exceptional
gradient lines, then
the resulting $CB$-diagram is admissible.
\end{lemma}
\begin{proof}
According to definition \ref{adm}
it is enough to prove that
there is a family of diffeomorphisms, providing orbit isotopies for $t<0$ and $t>0$,
and such that at $t=0$, that is along $B_1$ or $B_2$,
two separatrixes of the saddle $s_1$ and $s_2$, respectively of $s_3$ and $s_4$, form
$\gamma_{s_1,s_2}$
and $\gamma_{s_3,s_4}$,
in the way explained in section \ref{bifloc} and shown in figure \ref{fig19}.19.
If the
saddles $s_i$ are all distinct, for $i=1,...,4$, there are disjoint neighbourhood
$U_{12}$ and $U_{34}$ containing respectively $\gamma_{s_1,s_2}$
and $\gamma_{s_3,s_4}$ and such that
$U_{12}\cap(W^u(s_k)\cup W^s(s_k))=\varnothing$ for $k=3,4$ and
$U_{34}\cap(W^u(s_l)\cup W^s(s_l))=\varnothing$ for $l=1,2$.
Because of the way the phase portraits in each subset determined by $B_1$ and $B_2$
are constructed, and since $U_{12}$ and $U_{34}$ are disjoint,
it follows that it is possible to find a family of diffeomorphisms as above
which is the identity on the complement of $U_{34}$ or $U_{12}$,
and thus providing the required bifurcations along $B_1$ and $B_2$.

If $s_1=s_3$ (a similar argument if $s_2=s_4$) 
then $\gamma_{s_1,s_2}\cup\gamma_{s_1,s_4}=W^u(s_1)$, 
and so the two separatrixes of $s_1$ forming
$W^s(s_1)$ lie on different sides with respect to $\gamma_{s_1,s_2}\cup
s_3\cup\gamma_{s_3,s_4}$: this ensures that there are neighbourhood
$U_2$ and $U_4$ of respectively $\gamma_{s_1,s_2}$
and $\gamma_{s_1,s_4}$,
and containing the stable and unstable manifolds of respectively $s_2$ and $s_4$,
intersecting each one at most in $s_1=s_3$, such
that $U_2\cap(W^u(s_4)\cup W^s(s_4))=\varnothing$ and 
$U_4\cap(W^u(s_2)\cup W^s(s_2))=\varnothing$. Now the proof goes on as
in the first part for distict saddles.
\end{proof}

The case of two intersecting bifurcation lines $B_1$ and $B_2$, corresponding respectively 
to saddle-to-saddle separatrixes $\gamma_{s_1,s_2}$
and $\gamma_{s_2,s_3}$, performs a different behaviour: the reason is that
there are three different ways of breaking
the two non-generic gradient line appearing at $z$,
giving rise to a codimension 1 bifurcation: 
breaking $\gamma_{s_2,s_3}$ and leaving only
$\gamma_{s_1,s_2}$, as occurs along $B_1$; or breaking $\gamma_{s_1,s_2}$ and
leaving $\gamma_{s_2,s_3}$, as occurs along $B_2$; or forming 
the saddle-to-saddle separatrix
$\gamma_{s_1,s_3}$. The bifurcation line $B_3$ corresponding to $\gamma_{s_1,s_3}$
appears in the $CB$-diagram as an half-line with origin in z.
A case of this kind is considered in \cite{M2} for the perturbed elliptic
umbilic.

\begin{lemma}
\label{intbifb}
If, in a $CB$-diagram, 
the bifurcation lines $B_1$ and $B_2$ corresponding to the saddle-to-saddle separatrixes
$\gamma_{s_1,s_2}$ and $\gamma_{s_2,s_3}$ intersect each other in $z$ and if the phase
portrait at $z$ is that of a gradient vector field, then, to be admissible, the $CB$-diagram
must contain also a bifurcation half-line $B_3$, corresponding to 
the saddle-to-saddle separatrix $\gamma_{s_1,s_3}$, and whose origin is z.
\end{lemma}
\begin{proof}
Consider the phase portrait over $z$, and observe that the two separatrices
of $s_2$, not forming exceptional gradient lines, lie on the same side with respect
to $\gamma_{s_1,s_2}\cup s_2\cup\gamma_{s_2,s_3}$. There are two ways of breaking
a saddle-to-saddle separatrix, yielding two phase portraits which are not
orbit isotopic. So, while in the cases considered in lemma \ref{intbif}, the existence
of disjoint neighbourhood, each one containing one of the saddle-to-saddle separtrixes, 
provided four non-orbit equivalent phase portraits as a result of the breaking
the two exceptional gradient lines, and corresponding to the four
subsets determined by the intersection of $B_1$ and $B_2$, now, because of the relative
position of $\gamma_{s_1,s_2}$ and $\gamma_{s_2,s_3}$, and of the remaining separatrixes
of $s_2$, there are five non-orbit isotopic phase portraits, shown
in figure \ref{fig23}.23. Two of these
differs by a bifurcation, associated to the saddle-to-saddle separatrix $\gamma_{s_1,s_3}$,
and represented by a bifurcation line $B_3$ lying in one of the subsets, determined 
by $B_1$ and $B_2$, and dividing it into two disjoint subset. In other words, $B_3$ is
an half-line with origin in $z$. In one of those two phase portraits 
$\gamma_{s_1,s_2}$ can occur but not $\gamma_{s_2,s_3}$, the opposite happens
in the second. The bifurcation corresponding to $\gamma_{s_1,s_3}$ allows to
switch from one to the other. This shows the admissibility of the $CB$-diagram containing
$B_1$, $B_2$ and $B_3$.
\end{proof}

Figure \ref{fig21}.21 represents the two possibilities of intersection of bifurcation
lines in an admissible $CB$-diagram.

\begin{center}
\setlength{\unitlength}{1cm}
\begin{picture}(9,4)
\label{fig21}
\thinlines

\put(2,2){\circle*{.1}}
\put(6,2){\circle*{.1}}

\scriptsize
\put(1,.7){$\displaystyle B_2$}
\put(3,.7){$\displaystyle B_1$}
\put(5,.7){$\displaystyle B_2$}
\put(7,.7){$\displaystyle B_1$}
\put(7.1,1.9){$\displaystyle B_3$}

\normalsize
\put(0,3){$\displaystyle (a)$}
\put(4,3){$\displaystyle (b)$}

\qbezier(1,1)(2,2)(3,3)
\qbezier(3,1)(2,2)(1,3)

\qbezier(5,1)(6,2)(7,3)
\qbezier(7,1)(6,2)(5,3)
\qbezier(6,2)(6.5,2)(7,2)

\end{picture}
$Fig.~\ref{fig21}.21:~Intersections~of~bifurcation~lines$
\end{center}

The structure of the phase portraits in the subsets determined
by intersection of bifurcation lines as in case (b) of figure \ref{fig21}.21 is as follows. 
The phase portrait
at the intersection point contains the exceptional gradient lines 
$\gamma_{s_1,s_2}$ and $\gamma_{s_2,s_3}$; the remaining separatrices
of $s_2$ lie on the same side with respect
to $\gamma_{s_1,s_2}\cup s_2\cup\gamma_{s_2,s_3}$ (see figure \ref{fig22}.22). 

\begin{center}
\setlength{\unitlength}{1cm}
\begin{picture}(8,6)
\label{fig22}
\thinlines

\put(2,3){\circle*{.1}}
\put(4,3){\circle*{.1}}
\put(6,3){\circle*{.1}}

\scriptsize
\put(1.7,3.1){$\displaystyle s_1$}
\put(6.1,3.1){$\displaystyle s_3$}
\put(3.9,2.7){$\displaystyle s_2$}
\put(2.5,2.7){$\displaystyle \gamma_{s_1,s_2}$}
\put(4.5,2.7){$\displaystyle \gamma_{s_2,s_3}$}

\normalsize

\qbezier(1,3)(4,3)(7,3)
\qbezier(2,1)(2,3)(2,5)
\qbezier(6,1)(6,3)(6,5)
\qbezier(4,3)(3.5,4)(3,5)
\qbezier(4,3)(4.5,4)(5,5)

\qbezier(3,3)(2.95,3.05)(2.9,3.1)
\qbezier(3,3)(2.95,2.95)(2.9,2.9)
\qbezier(5,3)(4.95,3.05)(4.9,3.1)
\qbezier(5,3)(4.95,2.95)(4.9,2.9)
\qbezier(1.5,3)(1.55,3.05)(1.6,3.1)
\qbezier(1.5,3)(1.55,2.95)(1.6,2.9)
\qbezier(6.5,3)(6.55,3.05)(6.6,3.1)
\qbezier(6.5,3)(6.55,2.95)(6.6,2.9)

\qbezier(2,2)(2.05,1.95)(2.1,1.9)
\qbezier(2,2)(1.95,1.95)(1.9,1.9)
\qbezier(2,4)(2.05,4.05)(2.1,4.1)
\qbezier(2,4)(1.95,4.05)(1.9,4.1)
\qbezier(6,2)(6.05,2.05)(6.1,2.1)
\qbezier(6,2)(5.95,2.05)(5.9,2.1)
\qbezier(6,4)(6.05,3.95)(6.1,3.9)
\qbezier(6,4)(5.95,3.95)(5.9,3.9)

\qbezier(3.5,4)(3.55,3.975)(3.6,3.95)
\qbezier(3.5,4)(3.475,3.95)(3.45,3.9)

\qbezier(4.5,4)(4.55,4.025)(4.6,4.05)
\qbezier(4.5,4)(4.475,4.05)(4.45,4.1)

\end{picture}
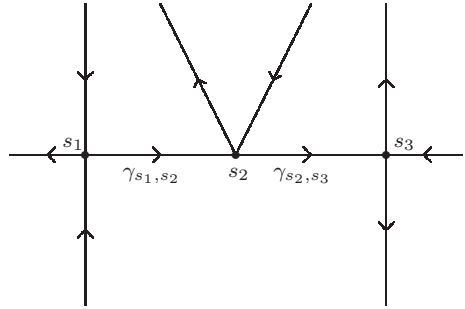
$Fig.~\ref{fig22}.22:~the~phase~portrait~over~the~intersection~of~B_1,~B_2~and~B_3$
\end{center}

Breaking the two exceptional gradient lines provides five phase portrait, as explained 
in lemma \ref{intbifb}, represented in figure \ref{fig23}.23.

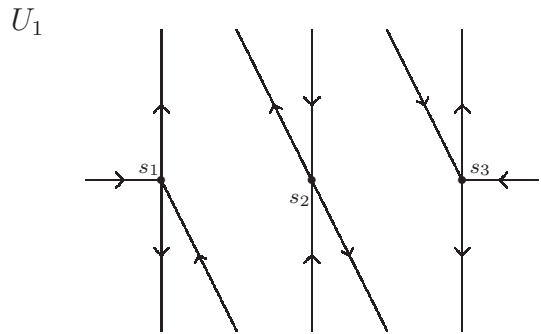
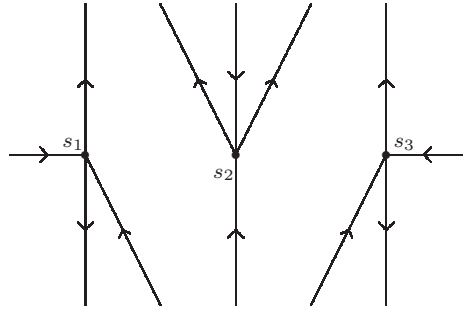
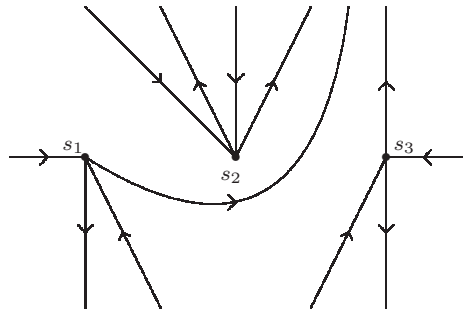
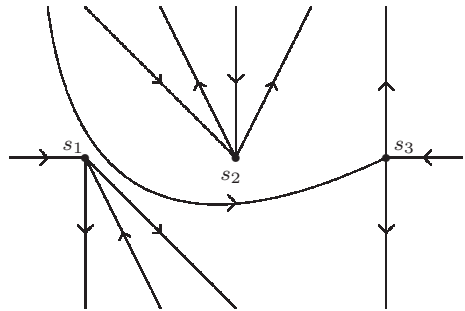
\begin{center}
\setlength{\unitlength}{1cm}
\begin{picture}(8,6)
\label{fig23}
\thinlines

\put(2,3){\circle*{.1}}
\put(4,3){\circle*{.1}}
\put(6,3){\circle*{.1}}

\scriptsize
\put(1.7,3.1){$\displaystyle s_1$}
\put(6.1,3.1){$\displaystyle s_3$}
\put(3.7,2.7){$\displaystyle s_2$}

\normalsize
\put(0,5){$\displaystyle U_1$}

\qbezier(1,3)(1.5,3)(2,3)
\qbezier(6,3)(6.5,3)(7,3)
\qbezier(2,1)(2,3)(2,5)
\qbezier(6,1)(6,3)(6,5)
\qbezier(4,1)(4,3)(4,5)
\qbezier(5,1)(4,3)(3,5)
\qbezier(6,3)(5.5,4)(5,5)
\qbezier(2,3)(2.5,2)(3,1)

\qbezier(1.5,3)(1.45,3.05)(1.4,3.1)
\qbezier(1.5,3)(1.45,2.95)(1.4,2.9)
\qbezier(6.5,3)(6.55,3.05)(6.6,3.1)
\qbezier(6.5,3)(6.55,2.95)(6.6,2.9)

\qbezier(2,2)(2.05,2.05)(2.1,2.1)
\qbezier(2,2)(1.95,2.05)(1.9,2.1)
\qbezier(2,4)(2.05,3.95)(2.1,3.9)
\qbezier(2,4)(1.95,3.95)(1.9,3.9)
\qbezier(4,2)(4.05,1.95)(4.1,1.9)
\qbezier(4,2)(3.95,1.95)(3.9,1.9)
\qbezier(4,4)(4.05,4.05)(4.1,4.1)
\qbezier(4,4)(3.95,4.05)(3.9,4.1)
\qbezier(6,2)(6.05,2.05)(6.1,2.1)
\qbezier(6,2)(5.95,2.05)(5.9,2.1)
\qbezier(6,4)(6.05,3.95)(6.1,3.9)
\qbezier(6,4)(5.95,3.95)(5.9,3.9)

\qbezier(3.5,4)(3.55,3.975)(3.6,3.95)
\qbezier(3.5,4)(3.475,3.95)(3.45,3.9)
\qbezier(2.5,2)(2.55,1.975)(2.6,1.95)
\qbezier(2.5,2)(2.475,1.95)(2.45,1.9)

\qbezier(4.5,2)(4.525,2.05)(4.55,2.1)
\qbezier(4.5,2)(4.45,2.025)(4.4,2.05)
\qbezier(5.5,4)(5.525,4.05)(5.55,4.1)
\qbezier(5.5,4)(5.45,4.025)(5.4,4.05)

\end{picture}

\setlength{\unitlength}{1cm}
\begin{picture}(8,5)

\thinlines

\put(2,3){\circle*{.1}}
\put(4,3){\circle*{.1}}
\put(6,3){\circle*{.1}}

\scriptsize
\put(1.7,3.1){$\displaystyle s_1$}
\put(6.1,3.1){$\displaystyle s_3$}
\put(3.7,2.7){$\displaystyle s_2$}

\normalsize
\put(0,5){$\displaystyle U_2$}

\qbezier(1,3)(1.5,3)(2,3)
\qbezier(6,3)(6.5,3)(7,3)
\qbezier(2,1)(2,3)(2,5)
\qbezier(6,1)(6,3)(6,5)
\qbezier(4,1)(4,3)(4,5)
\qbezier(4,3)(3.5,4)(3,5)
\qbezier(5,5)(4.5,4)(4,3)
\qbezier(6,3)(5.5,2)(5,1)
\qbezier(2,3)(2.5,2)(3,1)

\qbezier(1.5,3)(1.45,3.05)(1.4,3.1)
\qbezier(1.5,3)(1.45,2.95)(1.4,2.9)
\qbezier(6.5,3)(6.55,3.05)(6.6,3.1)
\qbezier(6.5,3)(6.55,2.95)(6.6,2.9)

\qbezier(2,2)(2.05,2.05)(2.1,2.1)
\qbezier(2,2)(1.95,2.05)(1.9,2.1)
\qbezier(2,4)(2.05,3.95)(2.1,3.9)
\qbezier(2,4)(1.95,3.95)(1.9,3.9)
\qbezier(4,2)(4.05,1.95)(4.1,1.9)
\qbezier(4,2)(3.95,1.95)(3.9,1.9)
\qbezier(4,4)(4.05,4.05)(4.1,4.1)
\qbezier(4,4)(3.95,4.05)(3.9,4.1)
\qbezier(6,2)(6.05,2.05)(6.1,2.1)
\qbezier(6,2)(5.95,2.05)(5.9,2.1)
\qbezier(6,4)(6.05,3.95)(6.1,3.9)
\qbezier(6,4)(5.95,3.95)(5.9,3.9)

\qbezier(3.5,4)(3.55,3.975)(3.6,3.95)
\qbezier(3.5,4)(3.475,3.95)(3.45,3.9)
\qbezier(2.5,2)(2.55,1.975)(2.6,1.95)
\qbezier(2.5,2)(2.475,1.95)(2.45,1.9)

\qbezier(4.5,4)(4.45,3.975)(4.4,3.95)
\qbezier(4.5,4)(4.525,3.95)(4.55,3.9)
\qbezier(5.5,2)(5.525,1.95)(5.55,1.9)
\qbezier(5.5,2)(5.45,1.975)(5.4,1.95)

\end{picture}

\setlength{\unitlength}{1cm}
\begin{picture}(8,5)

\thinlines

\put(2,3){\circle*{.1}}
\put(4,3){\circle*{.1}}
\put(6,3){\circle*{.1}}

\scriptsize
\put(1.7,3.1){$\displaystyle s_1$}
\put(6.1,3.1){$\displaystyle s_3$}
\put(3.8,2.7){$\displaystyle s_2$}

\normalsize
\put(0,5){$\displaystyle U_3$}

\qbezier(1,3)(1.5,3)(2,3)
\qbezier(6,3)(6.5,3)(7,3)
\qbezier(2,1)(2,2)(2,3)
\qbezier(5.5,5)(5,1.1)(2,3)
\qbezier(6,1)(6,3)(6,5)
\qbezier(4,3)(4,4)(4,5)
\qbezier(4,3)(3.5,4)(3,5)
\qbezier(4,3)(3,4)(2,5)
\qbezier(5,5)(4.5,4)(4,3)
\qbezier(6,3)(5.5,2)(5,1)
\qbezier(2,3)(2.5,2)(3,1)

\qbezier(1.5,3)(1.45,3.05)(1.4,3.1)
\qbezier(1.5,3)(1.45,2.95)(1.4,2.9)
\qbezier(6.5,3)(6.55,3.05)(6.6,3.1)
\qbezier(6.5,3)(6.55,2.95)(6.6,2.9)

\qbezier(4,4)(4.05,4.05)(4.1,4.1)
\qbezier(4,4)(3.95,4.05)(3.9,4.1)
\qbezier(2,2)(2.05,2.05)(2.1,2.1)
\qbezier(2,2)(1.95,2.05)(1.9,2.1)
\qbezier(6,2)(6.05,2.05)(6.1,2.1)
\qbezier(6,2)(5.95,2.05)(5.9,2.1)
\qbezier(6,4)(6.05,3.95)(6.1,3.9)
\qbezier(6,4)(5.95,3.95)(5.9,3.9)
\qbezier(3,4)(2.95,4)(2.9,4)
\qbezier(3,4)(3,4.05)(3,4.1)

\qbezier(3.5,4)(3.55,3.975)(3.6,3.95)
\qbezier(3.5,4)(3.475,3.95)(3.45,3.9)
\qbezier(2.5,2)(2.55,1.975)(2.6,1.95)
\qbezier(2.5,2)(2.475,1.95)(2.45,1.9)

\qbezier(4.5,4)(4.45,3.975)(4.4,3.95)
\qbezier(4.5,4)(4.525,3.95)(4.55,3.9)
\qbezier(5.5,2)(5.525,1.95)(5.55,1.9)
\qbezier(5.5,2)(5.45,1.975)(5.4,1.95)

\qbezier(4,2.4)(3.95,2.45)(3.9,2.5)
\qbezier(4,2.4)(3.95,2.35)(3.9,2.3)

\end{picture}

\setlength{\unitlength}{1cm}
\begin{picture}(8,5)

\thinlines

\put(2,3){\circle*{.1}}
\put(4,3){\circle*{.1}}
\put(6,3){\circle*{.1}}

\scriptsize
\put(1.7,3.1){$\displaystyle s_1$}
\put(6.1,3.1){$\displaystyle s_3$}
\put(3.8,2.7){$\displaystyle s_2$}

\normalsize
\put(0,5){$\displaystyle U_4$}

\qbezier(1,3)(1.5,3)(2,3)
\qbezier(6,3)(6.5,3)(7,3)
\qbezier(2,1)(2,2)(2,3)
\qbezier(4,1)(3,2)(2,3)
\qbezier(6,1)(6,3)(6,5)
\qbezier(4,3)(4,4)(4,5)
\qbezier(4,3)(3.5,4)(3,5)
\qbezier(4,3)(3,4)(2,5)
\qbezier(5,5)(4.5,4)(4,3)
\qbezier(6,3)(2,1.1)(1.5,5)
\qbezier(2,3)(2.5,2)(3,1)

\qbezier(1.5,3)(1.45,3.05)(1.4,3.1)
\qbezier(1.5,3)(1.45,2.95)(1.4,2.9)
\qbezier(6.5,3)(6.55,3.05)(6.6,3.1)
\qbezier(6.5,3)(6.55,2.95)(6.6,2.9)

\qbezier(4,4)(4.05,4.05)(4.1,4.1)
\qbezier(4,4)(3.95,4.05)(3.9,4.1)
\qbezier(2,2)(2.05,2.05)(2.1,2.1)
\qbezier(2,2)(1.95,2.05)(1.9,2.1)
\qbezier(6,2)(6.05,2.05)(6.1,2.1)
\qbezier(6,2)(5.95,2.05)(5.9,2.1)
\qbezier(6,4)(6.05,3.95)(6.1,3.9)
\qbezier(6,4)(5.95,3.95)(5.9,3.9)
\qbezier(3,4)(2.95,4)(2.9,4)
\qbezier(3,4)(3,4.05)(3,4.1)

\qbezier(3.5,4)(3.55,3.975)(3.6,3.95)
\qbezier(3.5,4)(3.475,3.95)(3.45,3.9)
\qbezier(2.5,2)(2.55,1.975)(2.6,1.95)
\qbezier(2.5,2)(2.475,1.95)(2.45,1.9)

\qbezier(4.5,4)(4.45,3.975)(4.4,3.95)
\qbezier(4.5,4)(4.525,3.95)(4.55,3.9)
\qbezier(3,2)(2.95,2)(2.9,2)
\qbezier(3,2)(3,2.05)(3,2.1)

\qbezier(4,2.4)(3.95,2.45)(3.9,2.5)
\qbezier(4,2.4)(3.95,2.35)(3.9,2.3)

\end{picture}

\setlength{\unitlength}{1cm}
\begin{picture}(8,5)

\thinlines

\put(2,3){\circle*{.1}}
\put(4,3){\circle*{.1}}
\put(6,3){\circle*{.1}}

\scriptsize
\put(1.7,3.1){$\displaystyle s_1$}
\put(6.1,3.1){$\displaystyle s_3$}
\put(3.8,2.7){$\displaystyle s_2$}

\normalsize
\put(0,5){$\displaystyle U_5$}

\qbezier(1,3)(1.5,3)(2,3)
\qbezier(6,3)(6.5,3)(7,3)
\qbezier(2,1)(2,2)(2,3)
\qbezier(4,1)(3,2)(2,3)
\qbezier(6,1)(6,3)(6,5)
\qbezier(4,3)(4,4)(4,5)
\qbezier(4,3)(3.5,4)(3,5)
\qbezier(4,3)(3,4)(2,5)
\qbezier(6,3)(5.5,4)(5,5)
\qbezier(4,3)(4.5,2)(5,1)
\qbezier(2,3)(2.5,2)(3,1)

\qbezier(1.5,3)(1.45,3.05)(1.4,3.1)
\qbezier(1.5,3)(1.45,2.95)(1.4,2.9)
\qbezier(6.5,3)(6.55,3.05)(6.6,3.1)
\qbezier(6.5,3)(6.55,2.95)(6.6,2.9)

\qbezier(4,4)(4.05,4.05)(4.1,4.1)
\qbezier(4,4)(3.95,4.05)(3.9,4.1)
\qbezier(2,2)(2.05,2.05)(2.1,2.1)
\qbezier(2,2)(1.95,2.05)(1.9,2.1)
\qbezier(6,2)(6.05,2.05)(6.1,2.1)
\qbezier(6,2)(5.95,2.05)(5.9,2.1)
\qbezier(6,4)(6.05,3.95)(6.1,3.9)
\qbezier(6,4)(5.95,3.95)(5.9,3.9)
\qbezier(3,4)(2.95,4)(2.9,4)
\qbezier(3,4)(3,4.05)(3,4.1)

\qbezier(3.5,4)(3.55,3.975)(3.6,3.95)
\qbezier(3.5,4)(3.475,3.95)(3.45,3.9)
\qbezier(2.5,2)(2.55,1.975)(2.6,1.95)
\qbezier(2.5,2)(2.475,1.95)(2.45,1.9)

\qbezier(3,2)(2.95,2)(2.9,2)
\qbezier(3,2)(3,2.05)(3,2.1)


\qbezier(4.5,2)(4.525,2.05)(4.55,2.1)
\qbezier(4.5,2)(4.45,2.025)(4.4,2.05)
\qbezier(5.5,4)(5.525,4.05)(5.55,4.1)
\qbezier(5.5,4)(5.45,4.025)(5.4,4.05)

\end{picture}
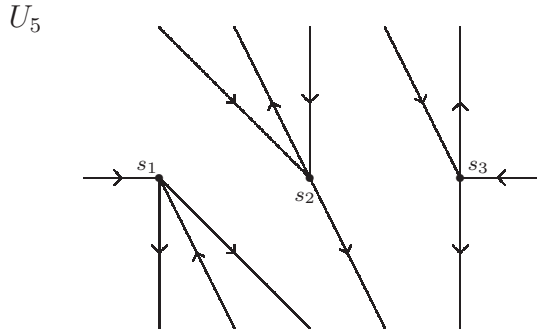
$Fig.~\ref{fig23}.23:~phase~portraits~in~the~subsets~determined~by~B_1,~B_2~and~B_3$
\end{center}

Note that the bifurcation line $B_3$, corresponding to the saddle-to-saddle separatrix 
$\gamma_{s_1,s_3}$, bounds $U_3$ and $U_4$. The bifurcation line
$B_1$, corresponding to $\gamma_{s_1,s_2}$, separates $U_1$ from $U_5$ and
$U_2$ from $U_3$, while $B_2$, corresponding to $\gamma_{s_2,s_3}$, 
separates $U_1$ from $U_2$ and
$U_4$ from $U_5$.

\section{Extension of quantum corrections}
The quantum corrections in definitions \ref{qcc} and \ref{qcb} 
allow to glue the holomorphic objects,
defined by means of Morse homology, on $U_1$ and $U_2$ along their common
boundary, when this is either a subset of the caustic $C$, consisting of folds not limit
points of the bifurcation locus $B$, or a subset of $B$ consisting of codimension 1
bifurcation points. This is a codimension 1 submanifold of $\R^2$. 
It remains to check that such holomorphic structure can be extended through
the codimension 2 subset of $\R^2$ formed by folds which are limit points of $B$
(that is, the intersection points of $C$ and $B$), by codimension 2 bifurcation points
(that is, the intersections of bifurcation lines) and cusps. To this purpose
it will be computed the monodromy of the holomorphic structure given by
quantum corrections and check that it is the identity.
We are not going to analyze 
in this paper the cusps: this was considered, though only 
for the particular case of the perturbed
elliptic umbilic, in \cite{M3}.

%
%

%
%
%
%
%
%
%
%

\begin{proposition}
Suppose that the caustic $C$ and a bifurcation line $B$ intersect as shown in figure 
\ref{fig24}.24., 
then the holomorphic structure of the mirror object can be extended through the
intersection point.
\end{proposition}

\begin{center}
\setlength{\unitlength}{1cm}
\begin{picture}(4,4)
\label{fig24}
\thinlines

\put(2,2){\circle*{.1}}

\scriptsize
\put(.7,1.9){$\displaystyle C$}
\put(1.9,.7){$\displaystyle B$}
\put(1.3,2.5){$\displaystyle U_1$}
\put(1.3,1.3){$\displaystyle U_1^{'}$}
\put(2.5,1.3){$\displaystyle U_2^{'}$}
\put(2.5,2.5){$\displaystyle U_2$}

\normalsize
\qbezier(1,2)(2,2)(3,2)
\qbezier(2,3)(2,1.5)(2,1)

\end{picture}

$Fig.~\ref{fig24}.24:~Intersection~of~C~and~B$
\end{center}

\begin{proof}
Suppose that $(n,s)$ is the birth-death pair associated to $C$, appearing in $U_1^{'}$
and $U_2^{'}$,
and $(s_1,s_2)$ the pair of bifurcating saddles associated to $B$. 
Note that $s\neq s_i$, for $i=1,2$. From
definition \ref{qcb} it follows that: the quantum corrections, 
glueing the mirror object over $U_2$, $U_1$ along $B$ and over $U_2^{'}$, $U_1^{'}$
along $B$, and induced respectively by $\tilde{M}_B^{U_2U_1}$, 
$\tilde{M}_B^{U_2^{'}U_1^{'}}=(\tilde{M}_B^{U_1^{'}U_2^{'}})^{-1}$,
coincide on all saddles except $s$, on which the former is not defined;
$\tilde{M}_B^{U_1^{'}U_2^{'}}$ is the
identity on all saddles, 
except on $s_1$: in fact, it acts as
$s_1\rightarrow s_1\pm s_2$, where the sign depends on orientation. 
On the other hand, by definition \ref{qcc}, the quantum corrections,
glueing along $C$
the mirror object over $U_1$, $U_1^{'}$, and over $U_2$, $U_2^{'}$, 
and denoted respectively by $\tilde{M}_C^{U_1U_1^{'}}$ and 
$\tilde{M}_C^{U_2U_2^{'}}$,
are both given, if $n$ is an unstable node, by the natural injection on saddles, if $n$ is a stable node,
by the natural injection on saddles not connected to $n$, respectively, in $U_1^{'}$
and $U_2^{'}$, and by
a shift by $s$ on the remaining saddles. Hence, if $n$ is unstable or 
if $n$ is stable but $s_1$ and $s_2$ are not
connected to $n$ both in $U_1^{'}$
and $U_2^{'}$, then clearly
$$\tilde{M}_B^{U_2U_1}\tilde{M}_C^{U_2^{'}U_2}
\tilde{M}_B^{U_1^{'}U_2^{'}}\tilde{M}_C^{U_1U_1^{'}}=Id$$
and so the holomorphic structure can be extended through the intersection point. 

This equality, and the same conclusion, 
holds also in the remaining cases, though more care must be paid: indeed,
if $n$ is stable node and only $s_1$ is connected to $n$ both in $U_1^{'}$
and $U_2^{'}$, then
$\tilde{M}_C^{U_1U_1^{'}}$  and $\tilde{M}_C^{U_2U_2^{'}}$ are the identity on $s_2$ and 
a shift by $s$ on $s_1$, that is, $s_1\rightarrow s_1\pm s$;
if $n$ is stable node, $s_1$ is connected to $n$ only in $U_1^{'}$ and $s_2$ is connected 
to $n$ both in $U_1^{'}$
and $U_2^{'}$, then, making for simplicity a choice of signs, although the argument
is independent from it, 
$\tilde{M}_C^{U_1U_1^{'}}$ is a shift by $s$ both on $s_1$ and $s_2$,
while $\tilde{M}_C^{U_2U_2^{'}}$ is the identity on $s_1$ and a shift by $s$ on $s_2$, 
therefore
the composition $\tilde{M}_C^{U_2^{'}U_2}\tilde{M}_B^{U_1^{'}U_2^{'}}
\tilde{M}_C^{U_1U_1^{'}}$
acts as follows
\begin{displaymath}
\begin{array}{ccccccc}
s_1 & \rightarrow & s_1+s & \rightarrow & s_1+s-s_2 & \rightarrow & s_1-s_2\\
s_2 & \rightarrow & s_2 & \rightarrow & s_2 & \rightarrow & s_2
\end{array}
\end{displaymath}

Finally, 
if $n$ is stable, $s_1$ is connected to $n$ only in $U_2^{'}$, and $s_2$ is connected 
to $n$ both in $U_1^{'}$
and $U_2^{'}$, then
$\tilde{M}_C^{U_1U_1^{'}}$ is the identity on $s_1$ and a shift by $s$ on $s_2$,
while $\tilde{M}_C^{U_2U_2^{'}}$ is a shift by $s$ on both $s_1$ and $s_2$, 
therefore
the composition $\tilde{M}_C^{U_2^{'}U_2}\tilde{M}_B^{U_1^{'}U_2^{'}}\tilde{M}_C^{U_1U_1^{'}}$
acts as follows
\begin{displaymath}
\begin{array}{ccccccc}
s_1 & \rightarrow & s_1 & \rightarrow & s_1-s_2-s & \rightarrow & s_1-s_2\\
s_2 & \rightarrow & s_2+s & \rightarrow & s_2+s & \rightarrow & s_2
\end{array}
\end{displaymath}

That the monodromy is the identity on nodes is easily verified, since $\tilde{M}_B^{U_2^{'}U_1^{'}}$ and 
$\tilde{M}_B^{U_2U_1}$ act as the identity on nodes.
\end{proof}

\begin{proposition}
Suppose that the caustic $C$ and a bifurcation line intersect as shown in figure 
\ref{fig25}.25,
then the holomorphic structure of the mirror object can be extended through the
intersection point.
\end{proposition}

\begin{center}
\setlength{\unitlength}{1cm}
\begin{picture}(4,4)
\label{fig25}
\thinlines

\put(2,2){\circle*{.1}}

\scriptsize
\put(.7,1.9){$\displaystyle C$}
\put(1.9,.7){$\displaystyle B$}

\put(1.3,2.5){$\displaystyle U_1$}
\put(1.3,1.3){$\displaystyle U_1^{'}$}
\put(2.5,1.3){$\displaystyle U_2^{'}$}

\normalsize

\qbezier(1,2)(2,2)(3,2)
\qbezier(2,2)(2,1.5)(2,1)

\end{picture}

$Fig.~\ref{fig25}.25:~intersection~of~C~and~B$
\end{center}

\begin{proof}
Suppose that $(s_1,s_2)$ is the pair of bifurcating saddles at $B$
and $(n,s_i)$ is the birth-death pair associated to $C$, for $i=1$ or $i=2$.
Let $s_3$, $s_4$, ..., $s_m$ be the saddles exhibited, besides $s_1$ and $s_2$, 
by the phase portrait in
$U_1^{'}$ and $U_2^{'}$. 
Denote by $M_C^{U_1U_1^{'}}$ and $M_C^{U_2^{'}U_1}$ the quantum corrections
glueing along $C$ the mirror object respectively over $U_1$, $U_1^{'}$, and over
$U_2^{'}$, $U_1$, and by $M_B^{U_1^{'}U_2^{'}}$ the quantum correction
glueing along $B$ the mirror object in $U_1^{'}$ and $U_2^{'}$.

Note that, for $j\geq3$, $\tilde{M}_B^{U_1^{'}U_2^{'}}(s_j)=s_j$ and since 
$s_j$ is connected to
$n$ in $U_1^{'}$ if and only if $s_j$ is connected to
$n$ in $U_2^{'}$, it follows $\tilde{M}_C^{U_1U_1^{'}}(s_j)=\tilde{M}_C^{U_1U_2^{'}}(s_j)$:
this implies that 
$$\tilde{M}_C^{U_2^{'}U_1}\tilde{M}_B^{U_1^{'}U_2^{'}}\tilde{M}_C^{U_1U_1^{'}}(s_j)=s_j$$
It remains to check the action of the above
composition of quantum corrections on $s_1$ and $s_2$.

Suppose first $n$ is an unstable node and $s_i=s_1$: by definition \ref{qcc}, 
$\tilde{M}_C^{U_1U_1^{'}}$ and $\tilde{M}_C^{U_1U_2^{'}}$ are the natural injection,
in particular, they map $s_2\rightarrow s_2$ (on $s_1$ they are not defined);
by definition \ref{qcb},
$\tilde{M}_B^{U_1^{'}U_2^{'}}(s_1)=s_1-s_2$ and $\tilde{M}_B^{U_1^{'}U_2^{'}}(s_2)=s_2$; 
hence 
$$\tilde{M}_C^{U_2^{'}U_1}\tilde{M}_B^{U_1^{'}U_2^{'}}\tilde{M}_C^{U_1U_1^{'}}=Id$$
thus the mirror object can be extended through the intersection point.

If $n$ is an unstable node but $s_i=s_2$, the same conclusion is achieved, the only
difference being that $\tilde{M}_C^{U_1U_1^{'}}$ and $\tilde{M}_C^{U_1U_2^{'}}$ 
are the natural injection on $s_1$ (and not defined on $s_2$).

Suppose now $n$ is a stable node and $s_i=s_2$. If $s_1$ is connected to $n$ in $U_1^{'}$
(but the same result is achieved also if $s_1$ is connected to $n$ in $U_2^{'}$),
and assuming $n(\gamma_{s_2,n})=-n(\gamma_{s_2,n})$
in $U_1^{'}$ (the argument works as well, up to signs, for the opposite choice), 
then $\tilde{M}_C^{U_1U_1^{'}}(s_1)=s_1+s_2$, 
$\tilde{M}_B^{U_1^{'}U_2^{'}}(s_1)=s_1-s_2$ and $\tilde{M}_B^{U_1^{'}U_2^{'}}(s_2)=s_2$,
while $\tilde{M}_C^{U_1U_2^{'}}(s_1)=s_1$  
since $s_1$ is not connected to $n$ in $U_2^{'}$.
Hence
$$\tilde{M}_C^{U_2^{'}U_1}\tilde{M}_B^{U_1^{'}U_2^{'}}\tilde{M}_C^{U_1U_1^{'}}=Id$$
and the mirror object can be extended through the intersection point.

Also when $n$ is a stable node but $s_i=s_1$ there is no monodromy given
by quantum corrections and the mirror object can be extended through the intersection
point: indeed, the maps $\tilde{M}_C^{U_2^{'}U_1}$, $\tilde{M}_B^{U_1^{'}U_2^{'}}$
and $\tilde{M}_C^{U_1U_1^{'}}$ are the identity on $s_2$.

That the monodromy is the identity on nodes is easily verified, since $\tilde{M}_B^{U_2^{'}U_1^{'}}$ and 
$\tilde{M}_B^{U_2U_1}$ act as the identity on nodes.
\end{proof}

Consider now intersection points of bifurcation lines: as seen, there are two cases.

\begin{proposition}
Suppose that two bifurcation lines $B_1$ and $B_2$ intersect as shown in figure 
\ref{fig26}.26,
then the holomorphic structure of the mirror object can be extended through the
intersection point.
\end{proposition}

\begin{center}
\setlength{\unitlength}{1cm}
\begin{picture}(4,4)
\label{fig26}
\thinlines

\put(2,2){\circle*{.1}}

\scriptsize
\put(1,.7){$\displaystyle B_1$}

\put(1.9,1.3){$\displaystyle U_1$}
\put(1.9,2.5){$\displaystyle U_3$}
\put(2.6,1.9){$\displaystyle U_2$}
\put(1.2,1.9){$\displaystyle U_4$}
\put(3,.7){$\displaystyle B_2$}

\normalsize

\qbezier(1,1)(2,2)(3,3)
\qbezier(3,1)(2,2)(1,3)

\end{picture}

$Fig.~\ref{fig26}.26:~Intersection~of~B_1~and~B_2$
\end{center}

\begin{proof}
Let $(s_1,s_1^{'})$ and $(s_2,s_2^{'})$ be the pair of bifurcating saddles 
corresponding respectively to
$B_1$ and $B_2$. By lemma \ref{intbif} at most either $s_1=s_2$
or $s_1^{'}=s_2^{'}$. Let $M^{U_1U_2}$, ..., $M^{U_4U_1}$ be the quantum corrections
glueing the mirror object over $U_1$ and $U_2$, ..., $U_4$ and $U_1$, 
along the common bifurcation line bounding these domains. 
Then, by definition \ref{qcb},
$M^{U_3U_4}=(M^{U_1U_2})^{-1}$ and $M^{U_4U_1}=(M^{U_2U_3})^{-1}$ and the
quantum corrections along $B_1$ commute with those along $B_2$. Therefore
$$M^{U_4U_1}M^{U_3U_4}M^{U_2U_3}M^{U_1U_2}=Id$$
and so the mirror object can be extended through the intersection point.
\end{proof}

\begin{proposition}
Suppose that two bifurcation lines $B_1$ and $B_2$ intersect as shown in figure 
\ref{fig27}.27,
then the holomorphic structure of the mirror object can be extended through the
intersection point.
\end{proposition}

\begin{center}
\setlength{\unitlength}{1cm}
\begin{picture}(5,4)
\label{fig27}
\thinlines

\put(2,2){\circle*{.1}}

\scriptsize
\put(1,.7){$\displaystyle B_1$}
\put(3,.7){$\displaystyle B_2$}
\put(3.1,1.9){$\displaystyle B_3$}

\put(1.9,1.3){$\displaystyle U_2$}
\put(1.9,2.5){$\displaystyle U_4$}
\put(2.6,1.6){$\displaystyle U_3$}
\put(2.6,2.2){$\displaystyle U_4$}
\put(1.2,1.9){$\displaystyle U_1$}

\normalsize

\qbezier(1,1)(2,2)(3,3)
\qbezier(3,1)(2,2)(1,3)
\qbezier(2,2)(2.5,2)(3,2)

\end{picture}

$Fig.~\ref{fig27}.27:~Intersection~of~B_1,~B_2~and~B_3$
\end{center}

\begin{proof}
By lemma \ref{intbifb}, assume that the phase portrait of $U_i$ is as represented in figure 
\ref{fig23}.23. 
We write the quantum corrections for each bifurcation line,
showing their action on generators, when non-trivial, and their associated matrix,
and then compute their composition:
\begin{enumerate}
\item
from $U_1$ to $U_2$ the quantum correction $\Psi_{12}$ is non-trivial on $s_2$
$$s_2\rightarrow s_2+\psi_{12}(s_3)$$
where $\psi_{12}(s_3)\in\{-1,1\}$, and its matrix is
\begin{displaymath}
\Psi_{12}=
\left[ \begin{array}{ccc}
1 & 0 & 0\\
0 & 1 & 0\\
0 & \psi_{12}(s_3) & 1
\end{array} \right]
\end{displaymath}
\item
from $U_2$ to $U_3$ the quantum correction $\Psi_{23}$ is non-trivial on $s_1$
$$s_1\rightarrow s_1+\psi_{23}(s_2)$$
where $\psi_{23}(s_2)\in\{-1,1\}$, and its matrix is
\begin{displaymath}
\Psi_{12}=
\left[ \begin{array}{ccc}
1 & 0 & 0\\
\psi_{23}(s_2) & 1 & 0\\
0 & 0 & 1
\end{array} \right]
\end{displaymath}
\item
from $U_3$ to $U_4$ the quantum correction $\Psi_{34}$ is non-trivial on $s_1$ 
$$s_1\rightarrow s_1+\psi_{34}(s_3)$$
where $\psi_{34}(s_3)\in\{-1,1\}$, and its matrix is
\begin{displaymath}
\Psi_{12}=
\left[ \begin{array}{ccc}
1 & 0 & 0\\
0 & 1 & 0\\
\psi_{34}(s_3) & 0 & 1
\end{array} \right]
\end{displaymath}
\item
from $U_4$ to $U_5$ the quantum correction $\Psi_{45}$ is non-trivial on $s_2$ 
$$s_2\rightarrow s_2+\psi_{45}(s_3)$$
where $\psi_{45}(s_3)\in\{-1,1\}$, and its matrix is
\begin{displaymath}
\Psi_{12}=
\left[ \begin{array}{ccc}
1 & 0 & 0\\
0 & 1 & 0\\
0 & \psi_{45}(s_3) & 1
\end{array} \right]
\end{displaymath}
\item
from $U_5$ to $U_1$ the quantum correction $\Psi_{51}$ is non-trivial on $s_1$ 
$$s_1\rightarrow s_1+\psi_{51}(s_2)$$
where $\psi_{51}(s_2)\in\{-1,1\}$, and its matrix is
\begin{displaymath}
\Psi_{51}=
\left[ \begin{array}{ccc}
1 & 0 & 0\\
\psi_{54}(s_2) & 1 & 0\\
0 & 0 & 1
\end{array} \right]
\end{displaymath}
\end{enumerate}
The composition $\Psi=\Psi_{51}\Psi_{45}\Psi_{34}\Psi_{23}\Psi_{12}$ has the
following action on generators:
\begin{eqnarray}
\Psi(s_1) & = & s_1+(\psi_{23}(s_2)+\psi_{54}(s_2))
+\psi_{34}(s_3)+\psi_{45}(s_3)\psi_{23}(s_2)
\nonumber\\
\Psi(s_2) & = & s_2+\psi_{12}(s_3)+\psi_{45}(s_3)
\nonumber\\
\Psi(s_3) & = & s_3
\nonumber 
\end{eqnarray}
Note now that $\psi_{23}(s_2)=-\psi_{54}(s_2)$ and 
$\psi_{12}(s_3)=-\psi_{45}(s_3)$: in fact, obeserving the phase portraits in figure 
\ref{fig23}.23, 
once an orientation for the separatrixes
of each saddle is chosen, the signs convention \ref{assumption} 
is satisfied by $\psi_{23}(s_2)$ if and only if
it is not satisfied by $\psi_{54}(s_2)$, and the same can be stated for
$\psi_{12}(s_3)$ and $\psi_{45}(s_3)$. In this way $\Psi$ can be simplified as:
\begin{eqnarray}
\Psi(s_1) & = & s_1+\psi_{34}(s_3)+\psi_{45}(s_3)\psi_{23}(s_2)
\nonumber\\
\Psi(s_2) & = & s_2
\nonumber\\
\Psi(s_3) & = & s_3
\nonumber 
\end{eqnarray}
It remains to prove that $\psi_{34}(s_3)+\psi_{45}(s_3)\psi_{23}(s_2)=0$.
Consider the phase portrait over $U_1$, shown in figure \ref{fig28}.28: choose an orientation
and compute the terms in the above equation.

\begin{center}
\setlength{\unitlength}{1cm}
\begin{picture}(8,6)
\label{fig28}
\thinlines

\put(2,3){\circle*{.1}}
\put(4,3){\circle*{.1}}
\put(6,3){\circle*{.1}}

\scriptsize
\put(1.7,3.1){$\displaystyle s_1$}
\put(6.1,3.1){$\displaystyle s_3$}
\put(3.7,2.7){$\displaystyle s_2$}

\put(.8,4){$\displaystyle W^u_1(s_1)$}
\put(.8,2){$\displaystyle W^u_2(s_1)$}
\put(2.5,2.2){$\displaystyle W^s_1(s_1)$}
\put(.5,3.2){$\displaystyle W^s_2(s_1)$}

\put(2.4,4){$\displaystyle W^u_1(s_2)$}
\put(4.6,2){$\displaystyle W^u_2(s_2)$}
\put(4.2,4){$\displaystyle W^s_1(s_2)$}
\put(3,1.6){$\displaystyle W^s_2(s_2)$}

\put(6.2,4){$\displaystyle W^u_1(s_3)$}
\put(6.2,2){$\displaystyle W^u_2(s_3)$}
\put(4.8,3.3){$\displaystyle W^s_1(s_3)$}
\put(6.5,3.2){$\displaystyle W^s_2(s_3)$}

\normalsize
\put(0,5){$\displaystyle U_1$}

\qbezier(1,3)(1.5,3)(2,3)
\qbezier(6,3)(6.5,3)(7,3)
\qbezier(2,1)(2,3)(2,5)
\qbezier(6,1)(6,3)(6,5)
\qbezier(4,1)(4,3)(4,5)
\qbezier(5,1)(4,3)(3,5)
\qbezier(6,3)(5.5,4)(5,5)
\qbezier(2,3)(2.5,2)(3,1)

\qbezier(1.5,3)(1.45,3.05)(1.4,3.1)
\qbezier(1.5,3)(1.45,2.95)(1.4,2.9)
\qbezier(6.5,3)(6.55,3.05)(6.6,3.1)
\qbezier(6.5,3)(6.55,2.95)(6.6,2.9)

\qbezier(2,2)(2.05,2.05)(2.1,2.1)
\qbezier(2,2)(1.95,2.05)(1.9,2.1)
\qbezier(2,4)(2.05,3.95)(2.1,3.9)
\qbezier(2,4)(1.95,3.95)(1.9,3.9)
\qbezier(4,2)(4.05,1.95)(4.1,1.9)
\qbezier(4,2)(3.95,1.95)(3.9,1.9)
\qbezier(4,4)(4.05,4.05)(4.1,4.1)
\qbezier(4,4)(3.95,4.05)(3.9,4.1)
\qbezier(6,2)(6.05,2.05)(6.1,2.1)
\qbezier(6,2)(5.95,2.05)(5.9,2.1)
\qbezier(6,4)(6.05,3.95)(6.1,3.9)
\qbezier(6,4)(5.95,3.95)(5.9,3.9)

\qbezier(3.5,4)(3.55,3.975)(3.6,3.95)
\qbezier(3.5,4)(3.475,3.95)(3.45,3.9)
\qbezier(2.5,2)(2.55,1.975)(2.6,1.95)
\qbezier(2.5,2)(2.475,1.95)(2.45,1.9)

\qbezier(4.5,2)(4.525,2.05)(4.55,2.1)
\qbezier(4.5,2)(4.45,2.025)(4.4,2.05)
\qbezier(5.5,4)(5.525,4.05)(5.55,4.1)
\qbezier(5.5,4)(5.45,4.025)(5.4,4.05)

\end{picture}

$Fig.~\ref{fig28}.28:~the~phase~portrait~over~U_1$
\end{center}

Observe that $\psi_{45}(s_3)$ is determined by the signs of $W^s_1(s_2)$ and $W^s_1(s_3)$,
$\psi_{23}(s_2)$ by the signs of $W^s_2(s_2)$ and $W^s_1(s_1)$, and
$\psi_{34}(s_3)$ by $W^s_1(s_1)$ and $W^s_1(s_3)$.
Suppose $n(W^s_1(s_2))=n(W^s_1(s_3))$: since $\psi_{45}(s_3)=-\psi_{12}(s_3)$,
then $\psi_{45}(s_3)=1$, moreover $n(W^s_2(s_2))=-n(W^s_1(s_2))$. Choose now
$n(W^s_1(s_1))$: if $n(W^s_1(s_1))=n(W^s_2(s_2))$ then $\psi_{23}(s_3)=-1$, on the
other hand, being $n(W^s_1(s_1))=-n(W^s_1(s_3))$, it also follows $\psi_{34}(s_3)=1$;
if $n(W^s_1(s_1))=-n(W^s_2(s_2))$ then $\psi_{23}(s_3)=1$ and $\psi_{34}(s_3)=-1$.
In both cases $\psi_{34}(s_3)+\psi_{45}(s_3)\psi_{23}(s_2)=0$.
The same conclusion is achieved supposing $n(W^s_1(s_2))=-n(W^s_1(s_3))$.
\end{proof}

We can sum up all the results in the following theorem:

\begin{theorem}
\label{theo}
For an admissible $CB$-diagram, such that for each $x\notin C\cup B$ Morse homology
is defined, that is, both the hypothesis of theorem \ref{jost} and assumption \ref{assumption} are fulfilled,
and where the orientation of gradient lines satisfies the signs convention \ref{bifor},
the quantum corrections introduced in definitions \ref{qcc} and \ref{qcb} allow to
extend the holomorphic structure of the mirror object to all folds of the caustic
and to all bifurcation points.
\end{theorem}

Another step is necessary to provide a globally defined holomorphic object
on the mirror fibration: to study its extensibility to cusps. This problem was analyzed in
\cite{M3} (but see also \cite{F2}) for the perturbed elliptic umbilic: in
that case, a quantum correction was defined, however, since it was related
to the existence of a spin strucure on the Lagrangian submaniold $L$ defined
by the generating function $f$ (and to the orientation problem of a family
in Morse theory and Floer theory), rather than to the bifurcations of the family
$f_x$, we prefer to postpone the analysis of cusps to another time.

\vskip 1cm
\rm
\noindent
GIOVANNI MARELLI\\
marelli@kusm.kyoto-u.ac.jp

\end{document}